\documentclass[amsfonts]{amsart}

\usepackage{amsmath} 

\theoremstyle{plain}
\newtheorem{theorem}{Theorem}[section]
\newtheorem{lemma}[theorem]{Lemma}
\newtheorem{proposition}[theorem]{Proposition}
\newtheorem{cor}[theorem]{Corollary}

\theoremstyle{definition}
\newtheorem{definition}[theorem]{Definition}
\newtheorem*{example}{Example}

\theoremstyle{remark}
\newtheorem*{rmrk}{Remark}

\usepackage[all]{xy}

\usepackage{enumitem}

\usepackage{float}
\usepackage{etoolbox}
\patchcmd{\subsection}{\bfseries}{\bf}{}{}
\patchcmd{\subsection}{-.5em}{.5em}{}{}

\patchcmd{\subsubsection}{\bfseries}{\bf}{}{}

\usepackage{graphicx} 
\newcommand{\iso}{\stackrel{\sim}{\to}}
\newcommand{\cR}{\mathcal R}
\newcommand{\cA}{\mathcal A}
\newcommand{\cB}{\mathcal B}

\newcommand{\cW}{\mathcal W}

\newcommand{\BLA}{\mathbf{\Lambda}}

\newcommand{\Span}{\mathrm{Span}}

\newcommand{\ed}{x}
\newcommand{\II}{\mathrm I}

\newcommand{\lc}{locally convex}

\newcommand{\eqns}{(\ref{eqns})}

\newcommand{\eqnsss}{(\ref{eqnsss})}

\newcommand{\Q}{Q^5}
\newcommand{\M}{M^4}
\newcommand{\N}{N^5}
\newcommand{\zm}{\ast}

\newcommand{\eto}{\xrightarrow{\sim}}

\newcommand{\tG}{\tilde \Gamma}
\newcommand{\PC}{\mathcal{PC}}

\renewcommand{\>}{\rangle}
\newcommand{\<}{\langle} 
\newcommand{\w}{{\mathbf w}}
\newcommand{\bv}{{\mathbf v}} 
\newcommand{\m}{\zeta}

\newcommand{\tO}{\widetilde{\mathbb O}}

\newcommand{\X}{\times}
\newcommand{\G}{{\rm G}_2}
\newcommand{\Aut}{{\rm Aut}}

\newcommand{\tr}{\mbox{tr}}
 
\renewcommand{\u}{{\bf u}} 
\newcommand{\g}{{\mathfrak{g}}}
\newcommand{\gl}{{\mathfrak{gl}}}
 
\renewcommand{\a}{{\bf a}}
\renewcommand{\b}{{\bf b}}
\newcommand{\bfc}{{\bf c}}

\renewcommand{\o}{ \boldsymbol{\omega}}
\renewcommand{\O}{ \boldsymbol{\Omega}}
\renewcommand{\P}{\mathbb{P}}

\newcommand{\pf}{{\n\em Proof. }}

\newcommand{\tQ}{\overline{Q}^5}

\renewcommand{\Im}{{\rm Im}}
\renewcommand{\Re}{{\rm Re}}

 \newcommand{\wnab}{\widetilde\nabla}
 \newcommand{\rh}{\rho_h}

 \newcommand{\CC}{\mathcal{C}}
 
\newcommand{\RP}{\R\P}
\newcommand{\RPt}{\RP^2}
\newcommand{\RPts}{\RP^{2*}}
\newcommand{\Rts}{(\R^3)^*}
\newcommand{\Rt}{\R^3}

\newcommand{\ent}{\Longrightarrow}
\newcommand{\st}{\, | \,}

\newcommand{\R}{\mathbb{R}}
\newcommand{\Rtt}{\R^{3,3}} 

\newcommand{\C}{\mathbb{C}} 

\newcommand{\T}{\mathbb{T}}

\renewcommand{\qed}{\hfill\mbox{$\Box$}}

\newcommand{\D}{{\mathcal D}} 
\newcommand{\tD}{\overline{\mathcal D}} 
\newcommand{\SL}{\mathrm{SL}}

\newcommand{\GL}{\mathrm{GL}}

\newcommand{\SLt}{{\SL_2(\R)}} 
\newcommand{\SLth}{{\SL_3(\R)}} 
\newcommand{\SO}{\mathrm{SO}}

\newcommand{\slt}{\mathfrak{sl}_2(\R)} 
\newcommand{\slth}{\mathfrak{sl}_3(\R)} 
\renewcommand{\sl}{\mathfrak{sl}}

\newcommand{\h}{{\mathfrak h}} 
\newcommand{\I}{{\mathbb I^3}}

\newcommand{\so}{\mathfrak{so}} 
\newcommand{\End}{{\rm End}} 
\newcommand{\Ker}{{\rm Ker}}
\newcommand{\Ad}{{\rm Ad}} 
\renewcommand{\mod}{\hbox{mod }}

\newcommand{\n}{\noindent} 
\newcommand{\bs}{\bigskip}
 
\newcommand{\sn}{\smallskip\n} 

\newcommand{\mn}{\medskip\noindent}
\newcommand{\bn}{\bs\n}

\newcommand{\p}{{\mathbf p} }
 
\newcommand{\tj}{{\tilde j} } 
\newcommand{\q}{\mathbf q }

\usepackage{etoolbox}
\makeatletter
\let\old@tocline\@tocline
\let\section@tocline\@tocline
\newcommand{\subsection@dotsep}{4.5}
\newcommand{\subsubsection@dotsep}{4.5}
\patchcmd{\@tocline}
  {\hfil}
  {\nobreak
     \leaders\hbox{$\m@th
        \mkern \subsection@dotsep mu\hbox{.}\mkern \subsection@dotsep mu$}\hfill
     \nobreak}{}{}
\let\subsection@tocline\@tocline
\let\@tocline\old@tocline

\patchcmd{\@tocline}
  {\hfil}
  {\nobreak
     \leaders\hbox{$\m@th
        \mkern \subsubsection@dotsep mu\hbox{.}\mkern \subsubsection@dotsep mu$}\hfill
     \nobreak}{}{}
\let\subsubsection@tocline\@tocline
\let\@tocline\old@tocline

\let\old@l@subsection\l@subsection
\let\old@l@subsubsection\l@subsubsection

\def\@tocwriteb#1#2#3{%
  \begingroup
    \@xp\def\csname #2@tocline\endcsname##1##2##3##4##5##6{%
      \ifnum##1>\c@tocdepth
      \else \sbox\z@{##5\let\indentlabel\@tochangmeasure##6}\fi}%
    \csname l@#2\endcsname{#1{\csname#2name\endcsname}{\@secnumber}{}}%
  \endgroup
  \addcontentsline{toc}{#2}%
    {\protect#1{\csname#2name\endcsname}{\@secnumber}{#3}}}%

\newlength{\@tocsectionindent}
\newlength{\@tocsubsectionindent}
\newlength{\@tocsubsubsectionindent}
\newlength{\@tocsectionnumwidth}
\newlength{\@tocsubsectionnumwidth}
\newlength{\@tocsubsubsectionnumwidth}
\newcommand{\settocsectionnumwidth}[1]{\setlength{\@tocsectionnumwidth}{#1}}
\newcommand{\settocsubsectionnumwidth}[1]{\setlength{\@tocsubsectionnumwidth}{#1}}
\newcommand{\settocsubsubsectionnumwidth}[1]{\setlength{\@tocsubsubsectionnumwidth}{#1}}
\newcommand{\settocsectionindent}[1]{\setlength{\@tocsectionindent}{#1}}
\newcommand{\settocsubsectionindent}[1]{\setlength{\@tocsubsectionindent}{#1}}
\newcommand{\settocsubsubsectionindent}[1]{\setlength{\@tocsubsubsectionindent}{#1}}

\renewcommand{\l@section}{\section@tocline{1}{\@tocsectionvskip}{\@tocsectionindent}{}{\@tocsectionformat}}%
\renewcommand{\l@subsection}{\subsection@tocline{1}{\@tocsubsectionvskip}{\@tocsubsectionindent}{}{\@tocsubsectionformat}}%
\renewcommand{\l@subsubsection}{\subsubsection@tocline{1}{\@tocsubsubsectionvskip}{\@tocsubsubsectionindent}{}{\@tocsubsubsectionformat}}%
\newcommand{\@tocsectionformat}{}
\newcommand{\@tocsubsectionformat}{}
\newcommand{\@tocsubsubsectionformat}{}
\expandafter\def\csname toc@1format\endcsname{\@tocsectionformat}
\expandafter\def\csname toc@2format\endcsname{\@tocsubsectionformat}
\expandafter\def\csname toc@3format\endcsname{\@tocsubsubsectionformat}
\newcommand{\settocsectionformat}[1]{\renewcommand{\@tocsectionformat}{#1}}
\newcommand{\settocsubsectionformat}[1]{\renewcommand{\@tocsubsectionformat}{#1}}
\newcommand{\settocsubsubsectionformat}[1]{\renewcommand{\@tocsubsubsectionformat}{#1}}
\newlength{\@tocsectionvskip}
\newcommand{\settocsectionvskip}[1]{\setlength{\@tocsectionvskip}{#1}}
\newlength{\@tocsubsectionvskip}
\newcommand{\settocsubsectionvskip}[1]{\setlength{\@tocsubsectionvskip}{#1}}
\newlength{\@tocsubsubsectionvskip}
\newcommand{\settocsubsubsectionvskip}[1]{\setlength{\@tocsubsubsectionvskip}{#1}}

\patchcmd{\tocsection}{\indentlabel}{\makebox[\@tocsectionnumwidth][l]}{}{}
\patchcmd{\tocsubsection}{\indentlabel}{\makebox[\@tocsubsectionnumwidth][l]}{}{}
\patchcmd{\tocsubsubsection}{\indentlabel}{\makebox[\@tocsubsubsectionnumwidth][l]}{}{}

\newcommand{\@sectypepnumformat}{}
\renewcommand{\contentsline}[1]{%
  \expandafter\let\expandafter\@sectypepnumformat\csname @toc#1pnumformat\endcsname%
  \csname l@#1\endcsname}
\newcommand{\@tocsectionpnumformat}{}
\newcommand{\@tocsubsectionpnumformat}{}
\newcommand{\@tocsubsubsectionpnumformat}{}
\newcommand{\setsectionpnumformat}[1]{\renewcommand{\@tocsectionpnumformat}{#1}}
\newcommand{\setsubsectionpnumformat}[1]{\renewcommand{\@tocsubsectionpnumformat}{#1}}
\newcommand{\setsubsubsectionpnumformat}[1]{\renewcommand{\@tocsubsubsectionpnumformat}{#1}}
\renewcommand{\@tocpagenum}[1]{%
  \hfill {\mdseries\@sectypepnumformat #1}}

\let\oldappendix\appendix
\renewcommand{\appendix}{%
  \leavevmode\oldappendix%
  \addtocontents{toc}{%
    \protect\settowidth{\protect\@tocsectionnumwidth}{\protect\@tocsectionformat\sectionname\space}%
    \protect\addtolength{\protect\@tocsectionnumwidth}{2em}}%
}
\makeatother

\makeatletter
\settocsectionnumwidth{2em}
\settocsubsectionnumwidth{2.5em}
\settocsubsubsectionnumwidth{3em}
\settocsectionindent{1pc}%
\settocsubsectionindent{\dimexpr\@tocsectionindent+\@tocsectionnumwidth}%
\settocsubsubsectionindent{\dimexpr\@tocsubsectionindent+\@tocsubsectionnumwidth}%
\makeatother

\settocsectionvskip{10pt}
\settocsubsectionvskip{0pt}
\settocsubsubsectionvskip{0pt}

\settocsectionformat{\bfseries}
\settocsubsectionformat{\mdseries}
\settocsubsubsectionformat{\mdseries}
\setsectionpnumformat{\bfseries}
\setsubsectionpnumformat{\mdseries}
\setsubsubsectionpnumformat{\mdseries}

\let\oldtableofcontents\tableofcontents
\renewcommand{\tableofcontents}{%
  \vspace*{-\linespacing}
  \oldtableofcontents}

\newcommand{\weta}{\hat\eta}

\newcommand{\om}[2]{\omega^{#1}_{\;#2}}

\newcommand{\be}{\begin{equation}}
\newcommand{\ee}{\end{equation}}

\newcommand{\met}{ \mathbf{ g}}

\title {The dancing metric, $\G$-symmetry  and  projective  rolling}

\author{Gil Bor, Luis Hern\'andez Lamoneda \& Pawel Nurowski} 

\AtEndDocument{\bigskip{\footnotesize%
  \textsc{CIMAT, A.P. 402, Guanajuato, Gto. 36000, Mexico} \par  
  \textit{E-mail address}, G. Bor: \texttt{gil@cimat.mx} \par
  \textit{E-mail address}, L.~Hern\'andez Lamoneda: \texttt{lamoneda@cimat.mx} \par
  \addvspace{\medskipamount}
 \textsc{Centrum Fizyki Teoretycznej, Polska Akademia Nauk, Al. Lotnik\'ow 32/46, 02-668 Warszawa, Poland} \par
  \textit{E-mail address}, P.~Nurowski: \texttt{nurowski@cft.edu.pl}
}}

\date{\today}

\begin{document}

\begin{abstract} The ``dancing metric" is a pseudo-riemannian metric $\met$ of signature (2,2) on the space $\M$ of  non-incident point-line pairs  in the real projective plane $\RPt$. The null-curves  of $(\M,\met)$ are given by the ``dancing condition": {\em the  point  is moving  towards a point  on the line, about which the line  is turning}. We  establish  a  dictionary between classical projective geometry (incidence, cross ratio, projective duality, projective invariants of plane curves\ldots) and pseudo-riemannian 4-dimensional conformal geometry (null-curves and geodesics, parallel transport, self-dual null 2-planes, the Weyl curvature,\ldots). There is also an unexpected bonus: by applying a  twistor construction to $(\M,\met)$, a $\G$-symmetry emerges, hidden deep in classical projective geometry. To uncover this symmetry, one needs to refine the ``dancing condition" by a higher-order condition, expressed in terms of the osculating conic along a plane curve. The outcome is a correspondence between curves in the projective plane and its dual, a projective geometry analog of the more  familiar ``rolling without slipping and twisting" for a pair of   riemannian surfaces. 

\end{abstract}

\maketitle

\tableofcontents

\section{Introduction} Let us consider the following system of 4
ordinary differential equations for 6 unknown functions $p_1,
p_2,p_3,$ $q^1,q^2,q^3$ of the variable $t$ 
$$ p_i
{dq^i\over dt}=0, \quad {dp_i\over dt}=\epsilon_{ijk}q^j{dq^k\over
dt},\quad i=1,2,3$$
 (we are using the summation
convention for repeated indices and the  symbol $\epsilon_{ijk}$, equal to $1$ for an even
permutation $ijk$ of $123$, -1 for an odd permutation, and 0
otherwise.)

It is convenient to recast these equations in  vector form by
introducing the notation
 $$\q=\left(\begin{matrix}q^1\\ q^2\\ q^3\end{matrix}\right)\in\R^3,\quad
\p=(p_1, p_2, p_3)\in\Rts.$$
Then, using the standard scalar and cross product
of vector calculus (and omitting the dot product symbol), the above
system can be written more compactly as
\begin{equation}\label{eqns} 
\p\q'=0, \quad \p'=\q\times \q'.
\end{equation}

This simple-looking  system of 4 ordinary differential equations for 6 unknown functions  enjoys a number of remarkable  properties and
interpretations, linking together old and new subjects, some of which we are going to explore in this article. The main themes  are 
\begin{itemize}
\item  generic  rank 2 distributions on 5-manifolds and their symmetries;
\item  4-dimensional pseudo-riemannian conformal geometry of split-signature;

\item  projective differential  geometry of plane curves. \end{itemize}
The relation between the first theme   and last  two is the main thrust  of this article.

\subsection{Summary of main results}

\subsubsection{A $(2,3,5)$-distribution  and its symmetries}\label{nonint}
Geometrically,  Eqns.~ \eqns\ define at each point  $(\q,\p)\in\R^6$ (away from some ``small" subset)  a $2$-dimensional subspace $\D_{(\q,\p)}\subset T_{(\q,\p)}\R^6$. Put together, these subspaces define (generically) a rank 2 distribution  $\D\subset T\R^6$, a field of tangent 2-planes,  so that the solutions to our system of equations are precisely the {\em integral curves} of $\D$: the  parametrized curves $(\q(t), \p(t))$ whose velocity vector $(\q'(t), \p'(t))$   lies in  $\D_{(\q(t),\p(t))}$ at each moment $t$.
 
 Furthermore,  we see readily from  Eqns.~\eqref{eqns} that the function $\p\q=p^iq_i:\R^6\to\R$
  is a ``conserved quantity" (constant along solutions), so $\D$
   is  tangent everywhere to the level surfaces of $\p\q$. 
   By a simple rescaling argument (Sect.\,\ref{TwoTwo}),  
   it suffices to consider one of its non-zero  level surfaces,  
   say  $\Q:=\{\p\q=1\}$.  Restricted to $\Q$, the equation 
   $\p\q'=0$ is a consequence of $\p'=\q\times\q',$ 
   hence our system of Eqns.~\eqns\ reduces to  
\begin{equation}\label{eqnsss}
\p\q=1, \quad \p'=\q\times\q'.
\end{equation}

The system $(\Q,\D)$ given by Eqns.~\eqref{eqnsss} does not have any more conserved quantities, since  $\D$ {\em bracket-generates} $TQ$,  in two steps: $\D^{(2)}=[\D,\D]$ is a rank 3 distribution and $\D^{(3)}=[\D, \D^{(2)}]=T\Q.$ Such a distribution is called  $(2,3,5)$-distribution. 

The study of  $(2,3,5)$-distributions  has a rich and fascinating history. Their local geometry was studied by \'Elie Cartan in his   celebrated ``5-variable paper" of 1910 \cite{C_5var}, where  he showed that the symmetry algebra of such a distribution (vector fields  whose flow preserves $\D$) is at most 14-dimensional. The most symmetric case  is realized, locally  uniquely, on a certain compact homogeneous 5-manifold $\tQ$ for the 14-dimensional simple exceptional non-compact Lie group $\G$ (Sect.\,\ref{Three}) equipped with a $\G$-invariant $(2,3,5)$-distribution $\tD$. This maximally-symmetric $(2,3,5)$-distribution $\tD$  is called the Cartan-Engel distribution, and was in fact  used  by \'E.~Cartan and F.~Engel  in 1893 \cite{C1, Eng} to  {\em define} $\g_2$ as the automorphism algebra of this distribution (the modern definition of  $\G$ as the automorphism group of  the octonions did not appear until 1908 \cite{Ca3}). 

Using \'E. Cartan's theory of $(2,3,5)$ distributions -- in particular, his {\em submaximality} result  -- we show (Thm.~\ref{sym}) that our distribution $(\Q,\D)$, as given by Eqns.~\eqnsss, is maximally-symmetric, i.e. admits a 14-dimensional symmetry algebra isomorphic to $\g_2$, and hence is {\em locally} diffeomorphic to the  Cartan-Engel distribution $(\tQ, \tD)$. Eqns.~\eqnsss\ thus provide an explicit model, apparently new,  for the Cartan-Engel distribution.

\begin{theorem} The system $(\Q,\D)$ given by  equations \eqnsss\ is  a $(2,3,5)$-distribution  with a 14-dimensional symmetry algebra, isomorphic to $\g_2$, the maximum possible for a $(2,3,5)$-distribution,  and is thus locally diffeomorphic to the Cartan-Engel $\G$-homogeneous distribution $(\tQ, \tD)$.
\end{theorem}

 Most of the symmetries of Eqns.~\eqref{eqns} implied by  this theorem are not obvious at all (``hidden"). There is however an 8-dimensional subalgebra  $\slth\subset \g_2$  of  ``obvious" symmetries, generated by 
$$(\q,\p)\mapsto (g\q, \p g^{-1}), \quad g \in \SLth$$
 (the cross-product  in  Eqns.~\eqnsss\ can be defined via  the standard volume form on $\R^3$,
hence the occurrence of $\SLth$; see Sect.~\ref{TwoFour}). The group $\SLth$ then acts transitively and effectively on $\Q$, preserving $\D$, and will be our main tool for studying the system \eqns. 

To explain the  appearance of  $\g_2$ as the symmetry algebra of $(\Q,\D)$ 
we construct in Sect.~\ref{Three} an embedding of $(\Q,\D)$  in the  ``standard model'' $(\tQ, \tD)$ of  the Cartan-Engel distribution, defined in terms of the split-octonions $\tO$. Using Zorn's ``vector matrices" to represent split-octonions (usually it is done with pairs of ``split-quaternions") we get explicit formulas for the symmetry algebra of Eqns.~\eqns.    

\begin{theorem} There is an embedding $\SLth\hookrightarrow \G=\Aut(\tO)$ and an $\SLth$-equivariant 
   embedding  $\R^6\hookrightarrow \RP^6=\P(\Im(\tO))$ (an affine chart), identifying  $\Q$  with  
   the open dense orbit  of $\SLth$ in $\tQ$ and mapping $\D$ over to $\tD$. 
   The $\G$-action on $\tQ$ defines a Lie subalgebra of vector fields on $\Q$ isomorphic to $\g_2$ (a 14-dimensional simple Lie algebra), forming the symmetry algebra of $(\Q, \D)$. 
\end{theorem}

\begin{cor}\label{cor1}For  each  $A=(a^i_j)\in\slth$, $\b=(b^i)\in \R^3$ and $\bfc=(c_i) \in \Rts$  the vector field   on $\R^6$ 
\begin{eqnarray*}
X_{A,\b, \bfc}&=&[ 2b^i + a^i_jq^j + \epsilon_{ijk}p^j c^k -(p_jb^j+ c_j q^j)q^i]\partial_{q^i} 
\\ &&
\qquad+[2c_i-a_i^jp_j  + \epsilon^{ijk}q_j b_k- (p_j b^j +c_j q^j )p_i]\partial_{p_i }
\end{eqnarray*}
is tangent to $\Q$ and preserves $\D$. The collection of these vector fields defines a 14-dimensional subalgebra of the Lie algebra of vector fields on $\Q$, isomorphic to $\g_2$, and forming the symmetry algebra of the system $(\Q, \D)$ defined by Eqns.~\eqns.
\end{cor}

\subsubsection{4-dimensional conformal  geometry in split signature}
Let $\M\subset\RPt\times\RPts$ be the (open dense) subset of {\em non-incident point-line pairs} $(q,p).$ There is a principal fibration 
 $\R^*\to\Q\to\M$ (the ``pseudo-Hopf-fibration") defined by regarding $(\q,\p)\in\Q$ as  homogeneous coordinates  of the pair $(q,p)=([\q], [\p])\in \M$. The fibration $\Q\to\M$ defines naturally a  pseudo-riemannian metric on  $\M$  by  a  standard procedure: restrict the flat $(3,3)$-signature metric on $\R^6$ given by $\p\q$  to $\Q$, then project to $\M$, using the fact that the principal $\R^*$-action on $\Q$ is by isometries. We call the resulting metric $\met$ on $\M$  the {\em dancing metric}. A similar procedure defines an orientation on $\M$. The dancing metric is a self-dual pseudo-riemannian metric of signature $(2,2)$, non-flat, irreducible,  $\SLth$-symmetric (as well as many other remarkable properties, see Thm.~\ref{Ten}). 

\begin{rmrk} Although only the conformal class $[\met]$ of the dancing metric is eventually used in this article, it is natural to consider the dancing metric $\met$ itself, as it is the unique (up to a constant)  $\SLth$-invariant metric in its class.
\end{rmrk}

 The main result in Sect.~\ref{pr} is  a correspondence between the geometries   of  $(\Q,\D)$ and $(\M,[\met])$. 

\begin{theorem}
The  above defined ``pseudo-Hopf-fibration"  $\Q\to \M$  establishes a bijection between  integral curves in   $(\Q,\D)$   and non-degenerate  null curves in $(\M,[\met])$ with   {\em parallel self-dual tangent null 2-plane}.\end{theorem}

The   condition ``parallel  self-dual tangent null 2-plane" on a null curve in an oriented  split-signature conformal 4-manifold can be regarded as ``one-half" of  the  geodesic equations. Every null-direction  is the unique intersection of two null 2-planes, one self-dual and the other anti-self-dual. It follows that given a null-curve in such a manifold  there are  two fields of tangent null 2-planes defined along it, one self-dual and the other anti-self-dual, intersecting in the tangent line. A null  curve is a geodesic if and only if its tangent line is parallel, which  is equivalent to the two  tangent fields of null 2-planes being  parallel; for our curves, only the self-dual field   is required to be parallel, hence ``half-geodesics".  

We derive various explicit formulas for the dancing metric. Perhaps the most elementary expression is the following: use the local coordinates $(x,y,a,b)$ on $\M$ where  $(x,y)$ are the Cartesian coordinates of a point $q\in\RPt$ (in some affine chart) and  $(a,b)$ the coordinates of a line $p\in\RPts$ given by $y=ax+b$. Then 
$$\met\sim  da[(y - b)dx - xdy] + db[adx - dy],$$
where $\sim$ denotes conformal equivalence. See 
Sect.~\ref{simple} for a quick derivation of this formula using the {\em dancing condition} (appearing also in the abstract to this paper, after which we name  the metric $\met$). An  explicit formula for the dancing {\em metric} $\met$ in  homogeneous coordinates  is given  in  Sect.~\ref{first} (Prop.~\ref{homo}). In Sect.~\ref{ss_cr} we give another formula for $\met$ in terms of  the cross-ratio.

Following the twistor construction in  \cite{AN}, we show that  $(\Q, \D)$ can  be naturally identified with the {\em non-integrable locus} of  the  total space of the {\em self-dual twistor fibration} $\RP^1\to \T^+(\M)\to \M$ associated with $(\M,[\met])$, equipped with its twistor distribution $\D^+$. The non-integrability of $\D$ is then seen to be equivalent to the non-vanishing of the self-dual Weyl tensor of $\met$.

This explains also why we do not look at the ``other-half" of the null-geodesic equations on $\M$. They correspond to integral curves of the twistor distribution $\D^-$ on the anti-self-dual twistor space $\T^-(\M)$, which turns out to be integrable, due to the vanishing of the anti-self-dual Weyl tensor of $(\M,[\met])$. The resulting ``anti-self-dual-half-geodesics" can be easily described and are rather uninteresting from the point of view of this article (see Thm.~\ref{Ten}, Sect.~\ref{proper}).

 \subsubsection{Projective geometry: dancing pairs and projective rolling.}\label{ProjGeom}

 Every integral curve of $(\Q,\D)$ projects, via $\Q\to\M\subset \RPt\times\RPts$, to a pair of curves $q(t), p(t)$ in $\RPt, \RPts$  (resp.). We offer two  interrelated projective geometric interpretations of the class of pairs of curves thus obtained:  ``dancing pairs" and ``projective rolling".   

By ``dancing" we refer to the interpretation of  $q(t), p(t)$   as the coordinated  motion of a non-incident point-line pair  in $\RPt$. We ask: what ``rules of choreography" should the  pair  follow so as to define (1) a null-curve in  $\M$ (2)  with a parallel self-dual tangent plane?  We call a pair of curves $q(t), p(t)$ satisfying these conditions a {\em dancing pair}. 

 The nullity condition on the pair  turns out to have a rather simple ``dancing" description. Consider a moving point  tracing a curve $q(t)$ in $\RPt$  with an associated  tangent line along it $q^*(t)\in\RPts$, the {\em dual curve} of  $q(t)$.  Likewise, a  moving line   in $ \RPt$ traces a curve $p(t)$ in $\RPts$,  whose dual curve $p^*(t)$ is a curve in $\RPt$, the   {\em envelope} of the family of lines in $\RPt$ represented by  $p(t)$, or  the curve in $\RPt$ traced out by the ``turning points" of the moving line $p(t)$.

\begin{theorem}\label{intro_dancing}A   non-degenerate parametrized curve  in $(\M,[\met])$ is null if and only if the corresponding pair of curves  $(q(t), p(t))$  satisfies the ``dancing condition": {\em at each moment $t$, the point $q(t)$ is moving towards the turning  point $p^*(t)$ of  the line $p(t)$}. 
\end{theorem}
\begin{figure}[h]\centering
\includegraphics[width=0.5\textwidth]{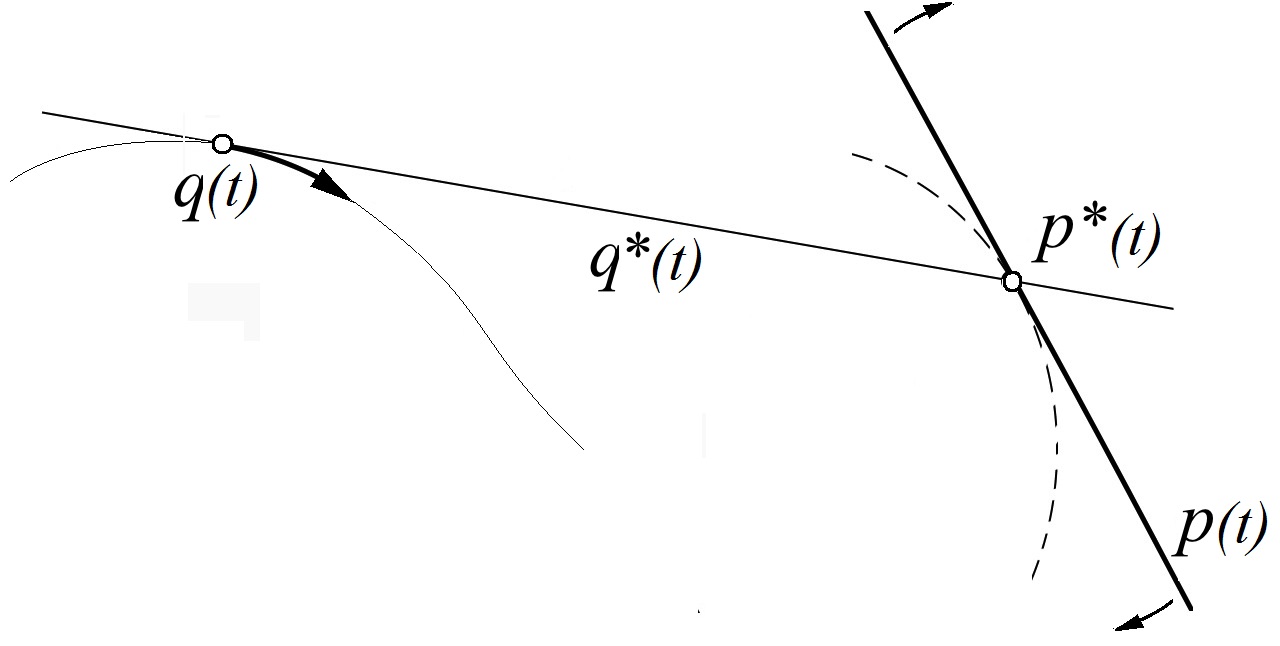}
 \caption{The  dancing condition}\label{fig:dance}
\end{figure}
In Sect.~\ref{sec:mate} we pose the following ``dancing mate" problem: fix  an arbitrary  curve $q(t)$ in $\RPt$ (with some mild non-degeneracy condition) and find its ``dancing mates" $p(t)$; that is, curves $p(t)$ in $\RPts$ such that   $(q(t),p(t))$ 
is a null curve in $(\M, [\met])$ with parallel self-dual tangent plane.  Abstractly, it is clear that there is a 3-parameter family of dancing mates for a given $q(t)$, corresponding to ``horizontal lifts" of $q(t)$ to integral curves of $(\Q, \D)$ via $\Q\to \M\to \RPt$, followed by the projection $\Q\to \M\to \RPts$. 

We study the resulting correspondence of curves $q(t)\mapsto p(t)$ from the point of view of classical projective differential geometry. We find that this correspondence preserves the natural {\em projective structures} on the curves $q(t),p(t)$,  but in general does not preserve the {\em projective arc length} nor the {\em projective curvature} (these are the basic projective invariants of a plane curve; any two of the  three invariants form a complete set of projective  invariants for  plane curves). We use the existence of a common projective parameter $t$ on $q(t),p(t)$ and the dancing condition to derive the  ``dancing mate equation":
\begin{equation}\label{mate}
y^{(4)}+2{y'''y'\over y}+3ry'+r'y=0.\end{equation}
Here,  $q(t)$ is given in homogeneous coordinates by a ``lift" $A(t)\in\R^3\setminus 0$ satisfying $A'''+rA=0$ for some function $r(t)$  (this is called the Laguerre-Forsyth form of the tautological equation for a plane curve) and the dual curve to $p(t)$ is given in homogeneous coordinates by $B=-y'A+yA'$. 

We study the special case  of  {\em the dancing mates of the circle}. That is, we look for dancing pairs $(q(t), p(t))$ where $q(t)$ parametrizes a fixed  circle $\CC\subset \RPt$ (or conic, projectively they are all equivalent). We show how the  above dancing mate equation  (\ref{mate}) reduces in this case to the 3rd order ODE $y'''y^2=1$. The dual dancing mates  $p^*(t)$ form a 3-parameter family of curves in the exterior of the circle $\CC$. We show here a computer generated image of a 1-parameter family of solutions (all other curves can be obtained from this family by the subgroup  $\SLt\subset\SLth$ preserving $\CC$).
\begin{figure}[h]\centering
\includegraphics[width=0.7\textwidth]{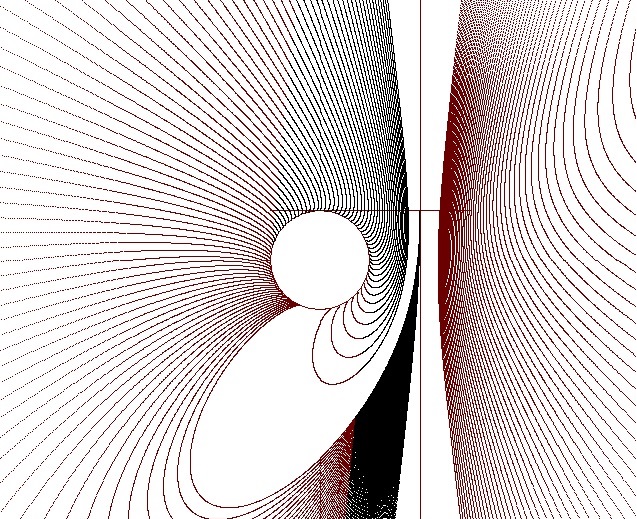}
\caption{Dancing mates around  the circle: the point-dancer moves along  the central circle, starting at its ``north pole", moving clockwise. 
The line-dancer starts in the vertical position (``$y$-axis''), keeping always tangent to one of the curves that spiral around the circle (the envelope of the line's motion). At all moments they  comply with the dancing condition; the figure shows the tangent direction of the point-dancer at the moment it passes through the north pole (horizontal line segment)  and its incidence  with the ``turning point" of the line-dancer at that moment.}\label{circ}
\end{figure}

\sn 

 As another illustration we give in Sect.~\ref{const}  examples of dancing pairs with {\em constant projective curvature} (logarithmic spirals, ``generalized parabolas", and exponential curves).  

Finally, in  Sect.~\ref{RT} we turn to   the  ``projective rolling" interpretation of Eqns.~\eqref{eqns}: imagine the curves $q(t)$ and $p(t)$ as the contact points of the two projective planes $ \RPt ,  \RPts$ as they ``roll" along each other. When rolling two surfaces along each other, one needs to pick at each moment, in addition to a pair of contact points $(q,p)$ on the two surfaces, an identification of   the tangent spaces $T_q\RPt$, $T_p\RPts$ at these points. 
In the case of usual rolling of riemannian surfaces, the identification is required to be an {\em isometry}. 
Here, we introduce the notion of ``projective contact" between the corresponding tangent spaces: 
it is an identification $\psi:T_q\RPt\to T_p\RPts$ (linear isomorphism) which  sends each line through 
$q$ to its intersection point with the line $p$ (thought of as a line in the tangent space to $\RPts$ at $p$).

 Now a simple calculation shows that this ``projective contact" condition is equivalent   to the condition that the 
 {\em graph of $\psi$  is a self-dual null 2-plane} in $T_q\RPt\oplus T_q\RPts\simeq T_{(q,p)}M.$ 
 The configuration space for projective rolling is thus the space $\PC$ of such projective contact 
 elements  $(q,p,\psi)$. Continuing the 
 analogy with the rolling of riemannian surfaces, we define {\em projective rolling without slipping} 
 as a curve    $(q(t), p(t), \psi(t))$ in $\PC$ satisfying $\psi(t)q'(t)=p'(t)$ for all $t$. 
\begin{theorem}
A curve $(q(t),p(t),\psi(t))$ in $\PC$ satisfies the no-slip condition $\psi(t)q'(t)=p'(t)$ if and only if $(q(t),p(t))$ is  a null curve in $(\M,[\met])$ (equivalently, it satisfies the dancing condition of Thm.~\ref{intro_dancing}). 
\end{theorem} 

Our next task is to  translate the ``half-geodesic" condition (parallel self-dual tangent plane) to rolling language. We use a notion of parallel transport of lines along (non-degenerate) curves in the projective plane, formulated in terms of Cartan's {\em development of the osculating conic along the curve}  (the unique  conic that touches a given point on the curve to 4th order; see Sect. \ref{osc}). We then define the ``no-twist" condition on a curve of projective contact elements  $(q(t),p(t),\psi(t))$ as follows: if  $\ell(t)$ is a parallel family of lines along $q(t)$ then $\psi(t)\ell(t)$ is a parallel family along $p(t)$. 

\begin{theorem}
A projective rolling curve $(q(t),p(t),\psi(t))$ satisfies the no-slip and no-twist condition if and only if  $(q(t),p(t))$ is a null-curve in $\M$ with parallel self-dual tangent plane. Equivalently, $(q(t),p(t))$ is the projection via $\Q\to \M$ of an integral curve of $(\Q,\D)$. 
\end{theorem}

The no-twist condition can be thought of as a  ``2nd dancing condition" for the dancing pair $(q(t),p(t))$; admittedly, it is a rather demanding  one: the dancers should be aware of the 5th order derivative of their motion in order to comply with  it\dots We believe there should be a simpler dancing rule that captures the no-twist condition but could not find it.

 

  \subsection{Background}
 Our original motivation for this article stems from  the article of the third author  with Daniel An
 \cite{AN}, where the twistor construction for split-signature 4-dimensional conformal metrics was introduced, raising  the following natural question: for which split-signature conformal 4-manifolds  $\M$ the associated twistor distribution $\D^+$ on $\T^+\M$ is a {\em flat}  $(2,3,5)$-distribution? (That is, with  $\g_2$-symmetry, the maximum possible). 

This is a hard problem, even when $\M$ is a product of riemannian surfaces $(\Sigma_i, g_i)$, $i=1,2,$  equipped with the difference metric   $g=g_1\ominus g_2.$
 In this case, the integral curves of the twistor distribution  can be interpreted as modeling rolling without slipping or twisting of the two surfaces along each other. It was known for a while to  R. Bryant (communicated in various places, like \cite{BM,Zel}) that the only case of pairs of {\em constant curvature}   surfaces that gives rises  to a flat $(2,3,5)$-distribution is that of curvature ratio 9:1 (or spheres of radius ratio 3:1, in the positive curvature case), but An-Nurowski found  in \cite{AN} a   new family of  examples, and it is still unknown if more examples exist. 

These new examples   of An-Nurowski motivated us looking for {\em irreducible} split signature 4-dimensional conformal metrics with flat twistor distribution. A natural place to start are homogeneous manifolds $\M=G/H$, with $G\subset \G$. We know of a few such examples, but we found the case of  $\SL_3(\R)/\GL_2(\R)$ studied in this article  the most attractive, due to its  projective geometric flavor (``dancing"  and ``projective rolling" interpretations), so we decided to dedicate an article to this example alone. 

\bn{\bf Acknowledgments.}  We are indebted to Richard Montgomery for many useful discussions,  revising  part of this work and making useful suggestions. We thank Serge Tabachnikov for suggesting we look for the formula appearing in Prop.~\ref{cr}. The first two authors   acknowledge support from grant 222870 of CONACyT. The third author acknowledges partial support from the Polish National Science Center (NCN) via DEC-2013/09/B/ ST1/01799.

\section{Geometric reformulation: a $(2,3,5)$-distribution} 
\subsection{First integral and reduction to the 5-manifold $\Q\subset \Rtt$}\label{TwoOne}
Let $\Rtt:=\R^3\oplus\Rts$, equipped with the quadratic form $(\q,\p)\mapsto \p\q=p^iq_i$. 
One can easily check that  $\p\q$ is  a first integral of equations \eqns\ (a
conserved quantity). That is,  for each $c\in\R$, a solution to \eqns\ that starts on the   level surface $$Q_c=\{(\q,\p)|\p\q=c\}$$ remains on
$Q_c$ for all times.  

Furthermore,  the map
$(\q,\p)\mapsto (\lambda \q, \lambda^2 \p)$, $\lambda\in\R$, maps
solutions on $Q_c$ to solutions on $Q_{\lambda^3c}$, hence it is
enough to study solutions of the system restricted to one of the 
(non-zero) level surfaces, say $$\Q:=\{(\q,\p)|\p\q=1\},$$ a 5-dimensional 
affine quadric
of signature
$ (3,3).$

\begin{rmrk} An {\em affine 
quadric}  is the non-zero  level set of a non-degenerate quadratic form on $\R^n$; its {\em signature}   is the signature of the defining quadratic form.
\end{rmrk}

\begin{rmrk}  We leave out the less interesting case of the zero level surface $Q_0$. 
\end{rmrk}

Now  restricted to $\Q$,  the equation $\p\q'=0$ is a consequence of $\p'=\q\times\q'$, hence we can replace equations \eqns\ with 
the somewhat simpler system 
\begin{equation}\label{eqnss}
\p\q=1, \quad \p'=\q\times \q'.
\end{equation}

\subsection{A rank 2 distribution $\D$ on $\Q$}\label{TwoTwo}
A geometric reformulation of  Eqns.~\eqref{eqnss} is the following: let us introduce the three 1-forms
$$\omega_i:=dp_i-\epsilon_{ijk}q^jdq^k\in\Omega^1(\Q),\quad i=1,2,3,$$
or in vector notation,
$$\o=d\p-\q\times d\q\in\Omega^1(\Q)\otimes\Rts.$$
Then the kernel of the 1-form $\o$ (the common kernel of its 3 components)  
defines  at each point $(\q,\p)\in \Q$ a 2-dimensional linear subspace  
$\D_{(\q,\p)}\subset T_{(\q,\p)}\Q$. Put together, these subspaces define 
a rank 2 distribution  $\D\subset T\Q$  (a field of tangent 2-planes on $\Q$),  
so that the solutions to our system of Eqns.~\eqref{eqnss} are precisely the 
{\em integral curves} of $\D$: the  parametrized curves $(\q(t) , \p(t) )$ 
whose velocity vector at each moment $t$  belongs to $\D_{(\q(t) ,\p(t) )}$.

\begin{proposition} The kernel of $\o=d\p-\q\times d\q$ defines on $\Q$ a rank 2 distribution $\D\subset T\Q$, whose integral curves are given by solutions to Eqns.~\eqref{eqnss}. 
\end{proposition}
\begin{proof} One checks easily that the 3 components of $\o$ are  linearly independent.
\end{proof}

\subsection{$\D$ is a $(2,3,5)$-distribution}  
We recall first   some standard terminology  from the general theory of distributions.
A distribution $\D$ on a manifold  is {\em integrable} if $[\D,\D]\subset\D$. It  is  {\em bracket generating} if one can obtain any tangent vector on the manifold by  successive Lie brackets of vector fields tangent to $\D$. Let $r_i=rank(\D^{(i)})$, where $\D^{(i)}:= [\D, \D^{(i-1)}]$  and $\D^{(1)}:=\D$. Then $(r_1, r_2, \ldots)$,  is the {\em growth
 vector} of $\D$. In general, the growth vector of a distribution may vary from point to point of the manifold, although not in our case, since our distribution is homogeneous, as we  shall soon see.  A distribution with constant growth vector is {\em regular}.  It can be shown that for  a regular  bracket-generating rank 2 distribution on  a 5-manifold there are only two possible growth vectors: $(2,3,4,5)$ (called {\em Goursat distributions}) or $(2,3,5)$, which is the generic case  (see \cite{BrHs}). 

\begin{definition} A {\em $(2,3,5)$-distribution} is  a bracket-generating  rank 2 distribution $\D$  on a 5-manifold $\Q$ with growth vector $(2,3,5)$ everywhere. That is, $\D^{(2)}=[\D,\D]$ is a rank 3 distribution, and 
$\D^{(3)}=[\D,\D^{(2)}]=T\Q$. 
\end{definition}

\begin{proposition} \label{prop235} $\D=\Ker(\o)\subset T\Q$, defined by Eqns.~\eqref{eqnss} above,  is a  $(2,3,5)$-distribution. 
\end{proposition}

 This is a calculation done most  easily  using the symmetries of the equations, so is postponed to the next subsection.

\subsection{$\SLth$-symmetry} \label{TwoFour}
A {\em symmetry} of a distribution $\D$ on a manifold $\Q$  is a diffeomorphism of $\Q$
 which preserves $\D$. An {\em infinitesimal symmetry} of $\D$  is a  vector field on $\Q$ whose flow preserves $\D$.

The use of the vector and scalar product on  $\R^3$ in Eqns.~\eqref{eqnss} may give the  impression that $\D$ depends on the euclidean structure on $\R^3$, so  $(\Q,\D)$ only admits $\SO(3)$ as  an obvious  group of symmetries (a 3-dimensional group). In fact, it is quite easy to see, as we will show now, that   $(\Q,\D)$   admits $\SLth$  as a symmetry group (8-dimensional). In the next section we will show the less obvious fact that  the symmetry {\em algebra} of $(\Q,\D)$ is $\g_2$ (14-dimensional).

Fix a volume form on $\R^3$, say $$vol:=dq^1\wedge dq^2\wedge dq^3,$$ and define the associated  covector-valued ``cross-product'' $\R^3\times \R^3\to\Rts$ by 
$$\bv\times \w:=vol(\bv,\w,\cdot),$$
or in coordinates, 
$$(\bv\times \w)_i=\epsilon_{ijk}v^jw^k, \quad i=1,2,3.$$

Let $\SLth$ be the group of $3\times 3$ matrices with real entries and determinant 1, acting on  $\Rtt$ by 
\be\label{action}
g\cdot (\q,\p)=(g\q,\p g^{-1})
\ee
 (recall that $\q$ is  a column vector and $\p$ is a row vector). Clearly, this $\SLth$-action leaves the quadratic form $\p\q$ invariant and thus leaves invariant also the quadric $\Q\subset\Rtt$.

Let $e_1,e_2,e_3$ (columns) be the standard basis of $\R^3$ and $e^1, e^2, e^3$ (rows) the dual basis of $\Rts$.

\begin{proposition} \label{prop:sym} 
 {\em (a)} $\SLth$ acts on $\Q$ transitively and effectively. The stabilizer of $(e_3, e^3)$ is the subgroup 
$$H_0=\left\{\left(\begin{array}{c|c}
A&\begin{matrix} 0\\ 0\end{matrix}\\
\hline
\begin{matrix}0& 0\end{matrix}&1
\end{array}
\right)\st  A\in\SLt\right\}.$$

\mn {\em (b)} $\SLth$ acts on $\Q$ by symmetries of  $\D$. 
\end{proposition}

\begin{proof}
 Part (a) is an easy calculation (omitted). For part (b), note that $\SLth$ leaves $vol$ invariant, hence the vector product $\R^3\times \R^3\to\Rts$ is $\SLth$-equivariant: $(g\bv)\times (g\w)=(\bv\times\w)g^{-1}.$ It follows that $\o=d\p-\q\times d\q$ is also $\SLth$-equivariant, $g^*\o=\o g^{-1}$, hence $\D=\Ker(\o)$ is $\SLth$-invariant.\end{proof}

\begin{proof}[Proof of Prop.~\ref{prop235}] Let $\h_0\subset\sl(3,\R)$ 
be the  Lie algebra of the stabilizer  at $(e_3, e^3)\in \Q$. Pick two 
elements $Y_1, Y_2\in\sl(3,\R)$ whose infinitesimal action at $(e_3, e^3)$ generates $\D$.  Then we need to show that 
$$Y_1, Y_2, [Y_1,Y_2],[Y_1,[Y_1,Y_2]],  [Y_2, [Y_1,Y_2]]$$ span $\sl(3,\R)$ mod $\h_0$ (this will show that $\D$ is $(2,3,5)$ at $(e_3,e^3) $, so by homogeneity everywhere.)
Now $\h_0$ consists of matrices of the form 
$$
\left(\begin{array}{c|c}
A&\begin{matrix} 0\\ 0\end{matrix}\\
\hline
\begin{matrix}0& 0\end{matrix}&0
\end{array}
\right),\quad A\in\slt.
$$
Furthermore, $Y\in\sl(3,\R)$ satisfies $Y\cdot(e_3, e^3)\in\D$ if and only if 
$$Y=\left(\begin{array}{cc}A&\bv\\ \bv^*& 0\end{array}\right),\quad \bv={v_1\choose v_2}, \quad \bv^*=(v_2, -v_1),  \quad A\in\slt.$$ 
We can thus take 
$$ 
Y_1=\footnotesize{ 
\left( 
\begin{array}{rrr}
0&0&1\\ 
0&0&0\\ 
0&-1&0
\end{array} 
\right), \quad
Y_2=\left( 
\begin{array}{rrr}
0&0&0\\ 
0&0&1\\ 
1&0&0
\end{array}
\right)},$$
then $[Y_1, Y_2]=Y_3,$ $[Y_1,Y_3]=Y_4$ and $[Y_2,Y_3]=Y_5,$ where 
$$ 
Y_3=\footnotesize{ 
\left( 
\begin{array}{rrr}
1&0&0\\ 
0&1&0\\ 
0&0&-2
\end{array}   
\right),  \quad
Y_4=
\left( 
\begin{array}{rrr}
0&0&-3\\ 
0&0&0\\ 
0&-3&0
\end{array}  
\right), \quad
Y_5=
\left( 
\begin{array}{rrr}
0&0&0\\ 
0&0&-3\\ 
3&0&0
\end{array} 
\right),}
$$
which together with $Y_1, Y_2$ span $\slth/\h_0\simeq T_{(e_3, e^3)}\Q$.
\end{proof}



\subsection{$\g_2$-symmetry via Cartan's submaximality} \label{three_five}Here we show that the symmetry {\em algebra} of our distribution $(\Q,\D)$, given by Eqns.~\eqref{eqnss}, is isomorphic to  $\g_2$, a 14-dimensional simple Lie algebra, the maximum possible for a $(2,3,5)$-distribution. We show this  as an immediate consequence of  a general theorem of Cartan (1910) on $(2,3,5)$-distributions. In the next section this ``hidden symmetry'' is explained and written down explicitly by defining an  embedding  of $(\Q,\D)$ in the standard $\G$-homogeneous model $(\tQ, \tD)$ using split-octonions.

\begin{rmrk} Of course, there is a third way, by ``brute force'', using computer algebra. We do not find it too illuminating but it does produce quickly  a  list of 14 vector fields on $\Rtt$, generating the infinitesimal $\g_2$-action, as given in Cor.~\ref{cor1}  or Cor~.\ref{cor3}  below. 
\end{rmrk}

In a well-known paper of 1910 (the ``5-variable paper"), Cartan proved
the following.

\begin{theorem}[Cartan,  \cite{C_5var}] 

\sn

\begin{enumerate}
\item The symmetry algebra of a $(2,3,5)$-distribution on a connected 5-manifold has dimension at most 14, in which case it is isomorphic to the real split-form of
the simple Lie algebra of type $\g_2$.
\item  All $(2,3,5)$-distributions with 14-dimensional symmetry algebra are
locally diffeomorphic. 
\item If the  symmetry algebra of a $(2,3,5)$-distribution 
has dimension $<14$ then it has  dimension at most 7.
\end{enumerate}
\end{theorem}
The last statement is sometimes referred to as Cartan's
``submaximality" result for $(2,3,5)$-distribution. A $(2,3,5)$-distribution with the maximal symmetry algebra $\g_2$ is called {\em flat}.

Using Prop.~\ref{prop:sym} and the fact
that $\SLth$ is $8$-dimensional we immediately conclude from Cartan's submaximality result for $(2,3,5)$-distributions  

\begin{theorem} \label{sym}
The symmetry algebra of the $(2,3,5)$-distribution defined by Eqns.~\eqref{eqnss}   is 14-dimensional,
isomorphic to the Lie algebra $\g_2$, containing the Lie  subalgebra isomorphic $\slth$ generated by the linear $\SLth$-action given by equation (\ref{action}). 
\end{theorem}

\begin{rmrk} In fact, Cartan \cite{C1} and Engel \cite{Eng} {\em defined} in 1893 the Lie algebra
$\g_2$ as the symmetry algebra of a certain 
$(2,3,5)$-distribution on an open set in $\R^5$, using formulas similar to our Eqns.~\eqref{eqnss}.
For example, Engel considers in \cite{Eng} the $(2,3,5)$-distribution obtained by restricting $ d\p-\q\times d\q$   to the linear subspace in $\Rtt$ given by $q_3=p^3$. 
\end{rmrk}

\section{$\G$-symmetry via split-octonions} \label{Three}In this section we describe  the relation between  the algebra of split-octonions $\tO$ and our equations \eqnsss, thus  explaining  the appearance of the ``hidden symmetries" as in Theorem \ref{sym} of the previous section.  We first review  some well-known facts concerning the algebra of split-octonions $\tO$ and its  automorphism group $\G$.  We then  define the ``standard model"   for the flat   $(2,3,5)$-distribution, 
a compact hypersurface   $\tQ\subset \RP^6$,  the projectivized null cone of imaginary split-octonions,  equipped with a $(2,3,5)$-distribution $\tD\subset T\tQ$. The group  $\G=\Aut(\tO)$  acts naturally on all objects defined in terms  of  the split-octonions, such as $\tQ$ and $\tD$. 

The relation of $(\tQ,\tD)$ with our system $(\Q,\D)$ is seen by finding an  embedding of groups $\SLth\hookrightarrow \G$ and an  $\SLth$-equivariant embedding  $(Q, \D)\hookrightarrow (\tQ, \tD)$. In this way we obtain  an  explicit realization of $\g_2$ as   the 14-dimensional symmetry algebra of $(Q,\D)$, containing the 8-dimensional subalgebra of ``obvious'' $\slth$-symmetries, as defined in  Eqn. (\ref{action})  of Sect. \ref{TwoFour}.  This  construction  explains  also why the infinitesimal $\g_2$-symmetry of $(\Q, \D)$ does not extend to a global $\G$-symmetry.

\subsection{Split-octonions via Zorn's vector matrices}
We begin with a brief review of the algebra of split-octonions, using a somewhat unfamiliar notation due to Max Zorn (of Zorn's Lemma fame in set theory), which we found quite useful in our context. See \cite{Katja} for a similar presentation.

The split-octonions $\tO$ is an $8$-dimensional non-commutative and non-associative real algebra, whose elements can be written  as ``vector matrices"
$$\m=\left(\begin{matrix}x&\q \\ \p&y\end{matrix}\right),\quad x,y\in\R,\quad \q\in\R^3, \quad \p\in\Rts,$$
with the ``vector-matrix-multiplication", denoted here by $\zm$, 
$$\m\zm \m'=\left(\begin{matrix}x&\q \\ \p&y\end{matrix}\right)\zm\left(\begin{matrix}x'&\q' \\ \p'&y'\end{matrix}\right):=
\left(\begin{array}{lr}xx'-\p'\q& x\q'+ y'\q+ \p\times\p'\\ x'\p+ y\p'+ \q\times\q'&yy'-\p\q'\end{array}\right), $$
where, as before,  we use the vector products $\Rt\times\Rt\to\Rts$ and $\Rts\times\Rts\to\Rt$, given by
$$\q\times\q':=vol(\q,\q', \cdot),\quad \p\times\p':=vol^*(\p,\p', \cdot),$$
$vol$ is the standard  volume form on $\R^3,$
$$vol=dq^1\wedge dq^2\wedge dq^3$$
and $vol^*$ is the dual volume form on $\Rts$, 
$$vol^*=dp_1\wedge dp_2\wedge dp_3.$$
In coordinates, 
$$(\q\times\q')_i=\epsilon_{ijk}q^jq'^k,\quad (\p\times\p')^i=\epsilon^{ijk}p_jp'_k.$$

\begin{rmrk}  These ``vector matrices" were  introduced   by Max Zorn in \cite{Z} (p.~144).   There are some minor variations  in the literature in the signs in the multiplication formula, but they are all equivalent to ours  by some simple change of variables (we are using Zorn's original formulas). For example, the formula in Wikipedia's article ``Split-octonion" is 
obtained from  ours by  the change of variable $\p\mapsto -\p.$
A better-known formula for octonion  multiplication uses pairs of quaternions, but we found the above formulas of Zorn more suitable (and it  fits also nicely with the original Cartan and Engel 1894 formulas). 
\end{rmrk}

Conjugation in $\tO$ is given by 
$$\m=\left(\begin{matrix}x&\q \\ \p&y\end{matrix}\right)\mapsto \overline{\m}=\left(\begin{array}{rr}y&-\q \\ -\p&x\end{array}\right),$$
satisfying  
$$\overline{\overline{\m}}=\m, \quad \overline{\m\zm\m'}=\overline{\m'}\zm\overline{\m},
\quad \m\zm\overline{\m}=\<\m,\m\>\rm I,$$
where $\rm I=\left(\begin{smallmatrix}1&0 \\ 0&1\end{smallmatrix}\right)$ and
$$\<\m,\m\>= xy+\p\q$$
is a quadratic form  of signature $(4,4)$ on $\tO$. 

Define as usual $$\Re(\m)=(\m+\overline\m)/2, \quad \Im(\m)=(\m-\overline\m)/2, $$
so that $$\tO=\Re(\tO)\oplus\Im(\tO),$$ where 
$\Re(\tO)=\R\rm I$ and $\Im(\tO)$ are vector matrices of the form $\m=\left(\begin{smallmatrix}x&\q \\ \p&-x\end{smallmatrix}\right)$. 


\subsection{About $\G$}

\begin{definition} $\G$  is the subgroup of $\GL(\tO)\simeq \GL_8(\R)$ satisfying
 $g(\m\zm\m')=g(\m)\zm g(\m')$ for all $\m, \m'\in\tO$. 
 \end{definition}

 \begin{rmrk} There are in fact three essentially distinct groups denoted by $\G$ in the literature: the  complex  Lie group $\G^\C$
  and its two real forms: the  compact form 
 and the  non-compact form, ``our" $\G$. See for example \cite{KN}. 
 \end{rmrk}
 \begin{proposition} Every  $g\in\G$ preserves the splitting $\tO=\Re(\tO)\oplus\Im(\tO)$. The action of $\G$ on $\Re(\tO)$ is trivial. Thus  $\G$ embeds naturally  in $\GL(\Im(\tO))\simeq \GL_7(\R)$. 
 
 \end{proposition}
 
\begin{proof} Let $g\in\G$. Since $\rm I$ is invertible so is $g(\rm I)$. 
 Now $g({\rm I})=g({\rm I}\zm{\rm I} )=g({\rm I})\zm g({\rm I})$, hence $g({\rm I})={\rm I}$. 
 It follows that $g$ acts trivially on $\Re(\tO)=\R\rm I$.

Next,  to show that $\Im(\tO)$ is $g$-invariant, define $S:= \{\m\in\tO|\m\zm\m=-\rm I\}$.
 Then it is enough to show that  (1) $S$ is $g$-invariant, 
 (2)  $S\subset \Im(\tO)$, 
 (3) $S$ spans $\Im(\tO)$. 

(1) is immediate from $g(-\mathrm I)=-\mathrm I.$ For (2),  let  $\m=\left(\begin{smallmatrix}x&\q \\ \p&-y\end{smallmatrix}\right)\in S,$ 
then $\m*\m=-{\rm I}\ent x^2-\p\q=y^2-\p\q=-1,  (x+y)\q=(x+y)\p=0\ent x+y=0\ent \m\in\Im(\tO)$. For (3), it is easy to find a basis of $\Im(\tO)$ in $S$.
\end{proof}

  The Lie algebra of $\G$ is the  sub-algebra  $\g_2\subset\End(\tO)$ of {\em derivations} of $\tO$: the elements  $X\in \End(\tO)$ such that $ X (\m\zm\m')= (X \m)\zm\m'+\m\zm (X \m')$ for all $\m, \m'\in \tO.$ It follows from the last proposition that $\g_2$ embeds as a sub-algebra of $\End(\Im(\tO))$.  \'E.  Cartan gave in his  1894 thesis  explicit  formulas for the image of this embedding, as follows.

For  each  $(A,\b, \bfc)\in \sl_3\oplus \R^3\oplus\Rts$   define  $\rho(A,\b,\bfc)\in\End(\Im(\tO))$, written  as a block matrix, corresponding to the decomposition  $\Im(\tO)\simeq\R^3\oplus\Rts\oplus\R$, 
$\left(\begin{smallmatrix}x&\q \\ \p&-x\end{smallmatrix}\right)\mapsto (\q,\p,x)$, 
by 
$$\rho(A,\b,\bfc)=
    \left(
        \begin{array}{ccc}
           A& R_\bfc &2\b\\
           L_\b &-A^t\,\,\,&2\bfc\\
           \bfc^t &\b^t& 0 
        \end{array}
    \right),$$
where   $L_\b:\Rt\to\Rts$ is given by $\q\mapsto \b\times \q$ and  $R_\bfc:\Rts \to\Rt$ is given by $\p\mapsto\p\times \bfc$.

Now define $\tilde\rho: \sl_3\oplus \R^3\oplus\Rts\to \End(\tO)$ by $$\tilde\rho(A,\b, \bfc)\m=\rho(A,\b, \bfc)\Im(\m).$$

Explicitly, we find
$$\tilde\rho(A,\b, \bfc) \left(\begin{matrix}x&\q \\ \p&y\end{matrix}\right)=
\left(\begin{matrix}\p\b+\bfc\q&A\q +(x-y)\b+\p\times\bfc\\ -\p A+\b\times\q+(x-y)\bfc&-\p\b-\bfc\q\end{matrix}\right)
.$$

\begin{proposition} The image of $\tilde\rho$ in $\End(\tO)$ is  $\g_2$.  That is, for all $(A,\b,\bfc)\in\sl_3\oplus  \R^3\oplus\Rts$,  $\tilde\rho(A,\b,\bfc)$ is a derivation of $\tO$  and all derivations of $\tO$ arise in this way. Thus  $\G$ is a 14-dimensional Lie group. It is a simple Lie group  of type $\g_2$ (the non-compact real form).

\end{proposition}

\begin{proof} (This is a sketch; for more details see for example \cite{Katja}). One  shows first  that $\tilde\rho(A,\a,\b)$ is a derivation by direct calculation. In the other direction, if $X\in\g_2,$ i.e. is a derivation, then its restriction to $\Im(\tO)$ is antisymmetric with respect to  the quadratic form ${\mathrm J}=x^2-\p\q$, i.e. is in the 21-dimensional Lie algebra $\so(4,3)$ of the orthogonal group corresponding to ${\mathrm J}$. One than needs  to show the vanishing of the  projection of $X$ to $\so(4,3)/\Im(\rho)$ (a 21-14=7 dimensional space). The latter  decomposes under $\SLth$ as $\R^3\oplus\Rts\oplus\R$, so by Schur Lemma it is enough to check the claim for one $X$  in each of the three irreducible summands. 

Now one can pick a Cartan subalgebra and root vectors showing that this algebra is of type $\g_2$ (see Cartan's thesis \cite{C_thesis}, p.  146).  
\end{proof}

\begin{rmrk}
Cartan gave the above representation of $\g_2$ in his 1894 thesis \cite{C_thesis} with no reference to octonions (the relation with octonions was published by him later in 1908 \cite{Ca3}). He presented $\g_2$ as the symmetry algebra of a rank 3 distribution on the null cone in $\Im(\tO)$.
 \end{rmrk}

\subsection{The distribution  $(\tQ,\tD)$}\label{dist} 

%

Imaginary split-octonions $\Im(\tO)$ satisfy  $\m=-\overline{\m}$ and are given by vector-matrices of the form 
$$\m=\left(\begin{matrix}x&\q \\ \p&-x\end{matrix}\right)$$ 
where $ (\q,\p,x)\in\Rt\oplus\Rts\oplus\R.$

\begin{definition} Let $\O:=\m \zm d \m$ (an $\tO$-valued 1-form on $\Im(\tO)$).  Explicitly, 
$$\O:=
\left(\begin{matrix}
x\,dx-\q \,d\p & 
x\,d\q - \q \,dx+\p\times d\p\\  
 \p\, dx- x\,d\p +
 \q\times d\q&
x\,dx-\p \,d\q\end{matrix}\right).$$
\end{definition}

\begin{proposition}  Let  $\Ker(\O)$ be the distribution  (with variable rank) on $\Im(\tO)$ annihilated by $\O$ and  let $  C\subset \Im(\tO)$ be the null cone, $C=\{\m\in\Im(\tO)|x^2-\p\q=0 \}$. Then 
$\Ker(\O)$  is
\begin{enumerate}
\item  $\G$-invariant,
\item  $\R^*$-invariant,  under     $\m\mapsto \lambda\m$, $\lambda\in\R^*,$ 
\item tangent to  $C\setminus 0$,
\item a rank 3 distribution when restricted to $C\setminus 0$,
\item the $\R^*$-orbits on $C$ are tangent to $\Ker(\O)$.

\end{enumerate}
\end{proposition}

\begin{proof} 

\sn

\begin{enumerate}
\item $\O$ is $\G$-equivariant, i.e. $g^*\O=g\O$ for all $g\in \G$,  hence $\Ker(\O)$ is 
$g$-invariant. Details: $g^*(\m \zm d\m)=(g\m)\zm d(g\m)=(g\m)\zm [g(d\m)]=g(\m\zm d\m)=g\O.$
\item $\lambda^*\O=\lambda^2\O\ent \Ker(g^*\O)=\Ker(\lambda^2\O)=\Ker(\O).$
\item $C$ is the 0 level set of $f(\m)=\m\zm \bar\m=-\m\zm \m,$ hence the tangent bundle to $C\setminus 0$ is the kernel of $df=-(d\m)\zm\m-\m \zm d \m=-\O-\overline{\O}=-2\Re(\O),$ hence 
$\Ker(\O)\subset\Ker(df).$

\item  Use the fact that $\G\X\R^*$ acts transitively on  $C\setminus 0$, so it is enough to check at say  $\q=e_1,$ $\p=0$,   $x=0$. Then $\Ker(\O)$ at this point is given by $dp_1=dq^2=dq^3=dx=0,$ which define a 3-dimensional subspace of $\Im(\tO))$. 
\item The $\R^*$-action is generated by the Euler vector field 
$$E=
p_i{\partial\over\partial_{p_i}}+q^i{\partial\over\partial_{q^i}}+
x{\partial\over\partial_{x}},$$ hence
$\O(E)=\m \zm d\m(E)=\m\zm\m=0,$ for $\m\in C$. 
\end{enumerate}
\end{proof}

\begin{cor} $\Ker(\O)$ descends to a $\G$-invariant rank-2  distribution $\tD$ on the projectivized null cone $\tQ=(C\setminus 0)/\R^*\subset \P( \Im(\tO))\cong\RP^6.$ 
\end{cor}

Now define an  embedding  $\iota:\Rtt\to \Im(\tO)$ by $(\q,\p)\mapsto (\q,\p,1).$ The pull-back of $\O$ by this map is easily seen to be
\begin{equation}\label{oct}
\iota^*\O=\left(\begin{array}{lr}
-\q \,d\p & 
d\q +\p\times d\p\\  
-d\p+
 \q\times d\q&
-\p \,d\q\end{array}
\right).
\end{equation}

Let $\SLth$ act on $\tO$ by 
$$ \left(\begin{matrix}x&\q \\ \p&y\end{matrix}\right)\mapsto
 \left(\begin{matrix}x&g\q \\ \p g^{-1}&y\end{matrix}\right), \quad g\in\SLth.$$
 This defines an embedding $\SLth\hookrightarrow \Aut(\tO).$  
\begin{theorem} Let $Q=\{\p\q=1\}\subset\Rtt$. Then 
\begin{enumerate}[leftmargin=18pt,label=(\alph*)]\setlength\itemsep{5pt}
\item the composition $$Q\stackrel{\iota}{\longrightarrow} C\setminus 0\stackrel{\R^*}{\longrightarrow}  \tQ,\quad (\q,\p)\mapsto [(\q,\p,1)]\in\tQ\subset \P(\Im(\tO))\cong\RP^6,$$
is an $\SLth$-equivariant embedding of  $(Q,\D)$ in $(\tQ,\tD).$ 

\item The image of $\Q\to\tQ$ is the  open-dense orbit of the $\SLth$-action on the projectivized null cone $\tQ\subset \P(\Im(\tO))\cong\RP^6$; its complement is a closed 4-dimensional  submanifold. 
\end{enumerate}

\end{theorem}

\sn\pf
(a) \,  Under $\SLth$, $\Im(\tO)$ decomposes as $\Rtt\oplus \R$, hence $\Rtt\to\Im(\tO)$, $(\q,\p)\mapsto [\q,\p,1]$, is an $\SLth$-equivariant embedding. Formula (\ref{oct}) for $i^*\O$ shows that $\D$ is mapped to $\tD$. 

\sn (b) \,  From the previous item, the image of $\Q$ in $\tQ$ is a single $\SLth$-orbit,  5-dimensional, hence open. It is dense, since the complement is a  4-dimensional submanifold in $\tQ$, given (in homogeneous coordinates) by the intersection of the hyperplane $x=0$ with the quadric $\p\q-x^2=0$. Restricted to $x=0$ (a 5-dimensional projective subspace in $\RP^6$) the equation $\p\q=0$ defines a smooth 4-dimensional   hypersurface, a projective quadric of signature $(3,3)$. 
\qed

\sn

Now if we consider  the projectivized $\g_2$-action on $[\Im(\tO)\setminus 0]/\R^*$ and pull it back to $\Rtt$ via $(\q,\p)\mapsto [(\q,\p, 1)]$, we obtain a realization of $\g_2$ as a Lie algebra of vector fields on $\Rtt$ tangent to $Q$, whose restriction to $Q$ forms the  symmetry algebra of $(Q,\D)$. 

\begin{cor}\label{cor3} For  each  $(A,\b, \bfc)\in\sl_3\oplus  \R^3\oplus\Rts$  the vector field   on $\Rtt$  
\begin{eqnarray*}
X_{A,\b, \bfc}&=&[ 2\b + A\q  + \p\times\bfc-(\p\b+\bfc \q)\q]\partial_\q 
\\ &&
+[2\bfc-\p A  + \q\times\b-(\p\b+\bfc\q)\p]\partial_{\p }
\end{eqnarray*}
is tangent to $Q\subset\Rtt$. The resulting 14-dimensional vector space of vector fields on  $Q$  forms the  symmetry algebra of $(Q,\D)$. 
\end{cor}

Explicitly, if  $A=(a^i_j), $ $\b=(b^i),$  $\bfc=(c_i), $ then 
\begin{eqnarray*}
X_{A,\b, \bfc}&=&[ 2b^i + a^i_jq^j + \epsilon_{ijk}p^j c^k -(p_jb^j+ c_j q^j)q^i]\partial_{q^i} 
\\ &&
+[2c_i-a_i^jp_j  + \epsilon^{ijk}q_j b_k- (p_j b^j +c_j q^j )p_i]\partial_{p_i }.
\end{eqnarray*}

\begin{proof} Let $\u=(\q,\p)\in\Rtt$, then $\iota(\u)=(\u,1)\in\Im(\tO)$. Any  linear vector field $X$ on $\Im(\tO)$ can be block decomposed as $$(X_{11}\u+ X_{12}x)\partial_\u + (X_{21}\u+ X_{22}x)\partial_x,$$ with 
$$X_{11}\in\End(\Rtt), \quad X_{12}\in\Rt, \quad X_{21}\in(\Rtt)^*, \quad X_{22}\in\R.$$
The induced vector field on $\Rtt$, obtained by  projectivization and pulling-back via $\iota$, is the quadratic vector field 
$$[X_{12} + (X_{11}-X_{22}x)\u-(X_{21}\u)\u]\partial_\u.$$
Now plug-in the formula for $X$ from last corollary. 
\end{proof}

\section{Pseudo-riemannian geometry in signature $(2,2)$}\label{pr}

 In this section we relate the geometry of the (2,3,5)-distribution $(\Q, \D)$ given by Eqns.~(\ref{eqns}) to 4-dimensional  conformal  geometry, by giving $\Q$ the structure of a principal $\R^*$-bundle $\Q\to \M$, the  ``pseudo-Hopf-fibration'', inducing on  $\M$  a  split-signature pseudo-riemannian metric $\met$, which we call the ``dancing metric''; the  name is due to an amusing alternative  definition  of the {\em conformal class} $[\met]$ (see Def.~\ref{def_danc} of Sect.~\ref{sec_danc}).

We then show in Thm.~\ref{eleven} (Sect.~\ref{tang_SD}), using the Maurer-Cartan structure equations of $\SLth$, that the projection $\Q\to \M$ establishes a bijection between  integral curves  in  $(\Q,\D)$   and  (non-degenerate) null-curves in $(\M,[\met])$ with  {\em parallel  self-dual tangent null 2-plane.} 

A more conceptual explanation to Thm.~ \ref{eleven} is given in Thm.~\ref{ident}, where we  show  that  $(\Q, \D)$ can  be naturally embedded in  the  total space of the {\em self-dual twistor fibration} $\RP^1\to \T^+(\M)\to \M$ associated with $(\M,[\met])$, equipped with its canonical {\em twistor distribution} $\D^+$, as introduced in \cite{AN}. The non-integrability of $\D$ is then seen to be due  to the non-vanishing of the self-dual Weyl tensor of $\met$.

\subsection{The pseudo-Hopf-fibration and the dancing metric} 

\subsubsection{First definition of the dancing metric}\label{first}
Recall from Sect.~\ref{TwoOne} that $\Q=\{(\q,\p)|\p\q=1\}\subset \Rtt$ (the ``unit pseudo-sphere"). To each pair $(\q,\p)\in\Q$ we assign the pair $\Pi(\q, \p)=([\q], [\p])=(q, p)\in \RPt\times\RPts$, where  $q\in \RPt$, $p\in\RPts$ are the points with homogeneous coordinates $\q, \p$ (resp.). Let $\I\subset\RPt\times\RPts$ be the subset of pairs $(q,p)$ given in homogeneous coordinates by the equation $\p\q=0$, also called {\em incident pairs}  (the name comes from the geometric interpretation of such a pair  as  a ({\em point, line}) pair, such that the {\em line} passes through  the {\em point}; more on this in  Sect. 5).  It is easy to see from the equation $\p\q=0$ that $\I$ is a 3-dimensional  closed submanifold of $\RPt\times\RPts$.  Its   complement  $$\M:=(\RPt\times\RPts)\setminus \I$$ is  the set of {\em non-incident point-line pairs},   a connected  open dense  subset of $\RPt\times\RPts$. Clearly, if $\p\q=1$ then $([\q], [\p])\not\in\I$, thus $\Pi:\Q\to \M$ is well defined.

 Define  an $\R^*$-action on $\Q$, where $\lambda\in \R^*$  acts   by 
 \be\label{Raction}
 (\q, \p)\mapsto
(\lambda\q, \p/\lambda), \quad \lambda\in\R^*.\ee 
This is a free $\R^*$-action whose  orbits are precisely the fibers of 
$$\Pi:\Q\to \M, \quad (\q, \p)\mapsto([\q], [\p]).$$
That is, $\Pi$ is a principal $\R^*$-fibration.  Now the quadratic form $\p\q$ defines on  $\Rtt$    a flat split-signature metric, whose restriction to $\Q\subset \Rtt$ is a  $(2,3)$-signature metric. Furthermore,   the principal $\R^*$-action on $\Q$ is by isometries, generated by a  negative definite vector field. 
 Combining these we get:

\begin{proposition}[\bf definition of the dancing metric]\label{dancing}
Restrict the flat split-signature metric $-2d\p \, d\q$ on $\Rtt$ to $\Q$. Then there is a unique pseudo-riemannian metric $\met$ on $\M$, of signature $(2,2)$, rendering $\Pi:\Q\to \M$  a pseudo-riemannian  submersion. We call $\met$  the {\em dancing metric}. 
 \end{proposition}
 
\begin{rmrk} The factor $-2$ in the above definition is not essential and is introduced merely for simplifying later explicit formulas for $\met$. 
\end{rmrk}

\begin{rmrk} This definition  is analogous to the definition  of the Fubini-Study metric  on $\C\P^2$ via the (usual)  Hopf fibration $S^1\to S^5\to\C\P^2$. In fact, $\M$ is referred to  by some authors as  the ``para-complex projective plane'' and  $\met$ as the ``para-Fubini-Study metric"  \cite{A,CFG}. 
\end{rmrk}

Using the $\SLth$-invariance of $\met$ it is not difficult to come up with an explicit formula for $\met$ in homogenous coordinates  $\q, \p$ on $\RPt, \RPts$ (resp.).
\begin{proposition}\label{homo}
Let $\widetilde \Pi:\Rtt\setminus \{\p\q=0\}\to \M, $ $(\q, \p)\mapsto ([\q],[\p])$. Then 
\be\label{homog}
\widetilde \Pi^*\met=-2{(\q\times d\q)(\p\times d\p)\over (\p\q)^2}.
\ee
\end{proposition}

\begin{proof} The expression on the right of Eqn.~(\ref{homog}) is a quadratic 2-form, $\R^*\times \R^*$-invariant, $\widetilde\Pi$-horizontal (vanishes on
 $\widetilde\Pi$-vertical vectors) and $\SLth$-invariant. It thus descends to an 
 $\SLth$-invariant quadratic 2-form on $M$. 
By examining the isotropy representation of the stabilizer of a point in 
 $M$ (Eqn.~(\ref{isotropy}) below) we see that $M$ admits a unique $\SLth$ quadratic 2-form, 
 up to a constant multiple. It is thus sufficient to verify  the formula  on a single non-null vector, 
 say $e_1-e^1\in T_{(e_3, e^3)}Q.$ We omit this (easy) verification. 
 \end{proof}
 
\begin{rmrk} Using standard vector identities, formula~(\ref{homog}) can be rewritten also as 
 \be\label{homog_bis}
\widetilde \Pi^*\met=-2{(\p\q)(d\p\,  d \q)-(\p\, d\q) (d\p \,\q)\over (\p\q)^2}.
\ee
\newcommand{\z}{\mathbf z}
An advantage of this formula is that it makes sense in higher dimensions, definining the ``para-Fubini-Study" metric on $\left[\R^{n+1,n+1}\setminus\{\p\q=0\}\right]/(\R^*\times \R^*)$. It also compares nicely  with the usual formula for the (standard) Fubini-Study metric $\met_{FS}$ on $\C P^n=\left[\C^{n+1}\setminus\{0\}\right]/\C^*$, given  in homogenous coordinates $\z=(z_0, \ldots z_n)^t\in\C^{n+1}$, $\z^*:=\bar\z^t$, by
$$\widetilde\Pi^*\met_{FS}={(\z^*\z) (d\z^* \,d\z)-( \z^*\, d\z)( d\z^*\, \z)\over (\z^*\z)^2}.$$
\end{rmrk}

 We give later  three  more explicit formulas for  $\met$: in Prop.~\ref{all}a $\met$ is expressed  in terms of the Maurer-Cartan form  of $\SLth$, analogous to a formula for the Fubini-Study metric on $\C P^n$ in terms of the Maurer-Cartan form of $\mathrm{SU}_{n+1}$. In Prop.~\ref{cr} we give a ``cross-ratio"  formula for $
\met$. In Sect.~\ref{simple} we derive a simple formula  in local coordinates for the {\em conformal class}  $[\met]$, using the ``dancing condition''.   

\subsubsection{Orientation}\label{orientation} We define an orientation on $\M$ via its {\em para-complex structure}. Namely, using the decomposition $T_{(q,p)}\M=T_q\RPt\oplus T_p\RPts$, define  $K: TM\to TM$ by $K(q',p')= (q',-p')$.  A {\em para-complex basis} for $T_{(q,p)}\M$ is then an ordered   basis of the form $(v_1, v_2, K v_1,  K v_2)$. One can check easily that any two such bases are related by a matrix with positive determinant,  hence these bases give a well-defined orientation on $\M$. See Prop.~\ref{all}c below for an alternative definition via a volume form on $\M$, written in terms of the components of the Maurer-Cartan form of $\SLth$.    

\subsubsection{Some properties of the dancing metric}\label{proper}

The dancing metric has remarkable properties. We group in the next theorem some of them.  

\begin{theorem} \label{Ten}

\sn\begin{enumerate}[leftmargin=*]\setlength\itemsep{5pt}
\item  $(\M,\met)$ is the  homogeneous  symmetric space  $\SLth/H$, where   $H\simeq \GL_2(\R)$ (the precise subgroup $H$ is described below in Sect.~\ref{proofs}).  The $\SLth$-action on $\M$ is induced from the standard action on $\Rtt$, $([\q], [\p])\mapsto ([g\q], [\p g^{-1}]).$ The $\GL_2(\R)$-structure endows $\M$ with a structure of a {\em para-Kahler} manifold. 

\item  $(\M,\met)$  is a  complete, Einstein,   irreducible, pseudo-riemannian 4-manifold of signature $(2,2)$. 
It   is self-dual (with respect to the above orientation), i.e. its anti-self-dual Weyl tensor $\mathcal W^-\equiv0$,  but  is not conformally flat; its self-dual Weyl curvature tensor $\mathcal W^+$ is nowhere vanishing, of Petrov type $D$.

\item The splitting $T_{(q,p)}\M=T_q\RPt\oplus T_p\RPts$ equips $\M$ with a pair of complementary  null, self-dual,  parallel,  integrable,   rank 2 distributions. Their integral leaves generate a pair of foliations of $\M$ by  totally geodesic self-dual null surfaces, the fibers of the double  fibration
$$    \xymatrix{
        &\M  \ar[dl]_\pi \ar[dr]^{\bar\pi} & \\
                   \RPt&            & \RPts }
$$

\item $\M$ admits a 3-parameter family of anti-self-dual totally geodesic
 null surfaces, naturally parametrized by the incidence variety  $\I:=\{(\bar q, \bar p)|\bar q\in \bar p\}\subset\RPt\times\RPts.$
For  each incident pair $(\bar q, \bar p)\in\I$, the corresponding surface is the set $\Sigma_{\bar q, \bar p}$ of non-incident pairs $(q,p)$ such that $q\in\bar p$ and $\bar q\in p$. 
\end{enumerate}
\end{theorem}
\begin{figure}[h]\centering
\includegraphics[width=0.3\textwidth]{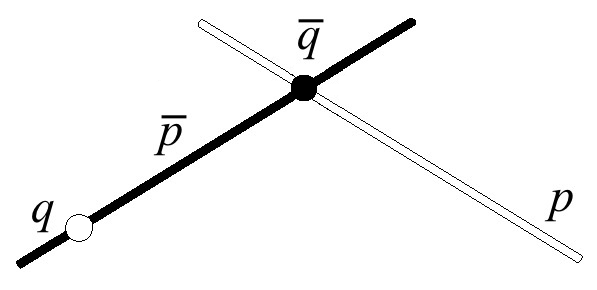}
\caption{The definition of   $\Sigma_{\bar q, \bar p}$}
\end{figure}
\begin{rmrk} The last  point (4) can be reformulated as follows: let $\N\subset  \M\times\I$ be defined via the incidence diagram above, i.e.
$$\N=\{(q,p,\bar q, \bar p)| q\not\in p, \bar q\in\bar  p, q\in\bar p, \bar q\in  p\}\subset\RPt\times \RPts\times\RPt\times \RPts.$$ 
Then $\N$ is a 5-dimensional submanifold of $\M\times \I $, equipped with the  double fibration
\be\label{double}    
\xymatrix{
&\N  \ar[dl]_{\pi_{12}} \ar[dr]^{\pi_{34}} & \\
\M& & \I }
\ee
The right-hand fibration $\pi_{34}:\N\to\I$ foliates $\N$ by 2-dimensional surfaces, each of which projects via $\pi_{12}:\N\to \M$
to one of the surfaces  $\Sigma_{\bar q, \bar p}$. That is, $\Sigma_{\bar q, \bar p}=\pi_{12}\left(\pi_{34}^{-1}(\bar q, \bar p)\right).$ The left-hand fibration $\pi_{12}:\N\to\M$ foliates $\N$  by projective lines and  can be naturally identified with the {\em anti-self-dual twistor fibration} $\T^-\M\to\M$ associated with $(\M,[\met])$ (see Sect.~\ref{twist}). The fibers of $\pi_{34}$ then correspond to the integral  leaves of the anti-self-dual twistor distribution $\D^-$, which is integrable in our case, due to the vanishing of $\mathcal W^-$ (see Cor.~\ref{grancor}, Sect.~\ref{proofs}  below).  
\end{rmrk}

 Most claims of this theorem can be found in  various sources in the literature (see e.g.~\cite{A} and the many references within). Using  the Maurer-Cartan equations of $\SLth$ (Sect.~\ref{proofs}), it is quite straightforward to prove these results. Alternatively, one can write down explicitly the dancing metric  in local coordinates (Sect.~\ref{simple}) and let a computer calculate curvature, symmetries etc. 

\subsection{Rudiments of 4-dimensional geometry in split-signature}
\subsubsection{Linear algebra}\label{linalg}  

Let  $V$ be  an oriented 4-dimensional real vector space equipped  with  a quadratic form $\<\,,\>$
 of signature $(++-\,-)$. It is convenient to introduce  {\em null bases} in such a $V$. This is a basis $\{e_1, e_2, e^1, e^2\}\subset V$  
such that 
$$\<e_a, e_b\>=\<e^a, e^b\>=0, \;  \<e^a, e_b\>=\delta^a_b, \qquad a,b=1,2.$$ 
Note that if $\{\ed^1, \ed^2,  \ed_1,  \ed_2 \}\subset V^*$ is the {\em dual basis} to a null basis, i.e. $\ed^a(e^b)=\ed_a(e_b)=0,  
\ed^a(e_b)=\ed_b(e^a)=\delta^a_b,$ then 
\be\label{metric}\<\,,\>=2(\ed^1 \ed_1+\ed^2 \ed_2).
\ee

\begin{rmrk} Our convention is that the symmetric tensor product $x y\in S^2\,V^*$ of two elements  $x,y\in V^*$   is the symmetric bilinear form 
\be (x y)(v,w):=[x(v)y(w)+y(v)x(w)]/2, \quad v, w\in V. 
\ee
\end{rmrk}

Now let  $vol:=\ed^1\wedge\ed^2\wedge\ed_1\wedge\ed_2\in \Lambda^4\,V^*$  and  $*:\Lambda^2\,V^*\to \Lambda^2\,V^*$ the  corresponding Hodge dual, satisfying  $\alpha\wedge*\beta=\<\alpha,\beta\>vol$, $\alpha, \beta\in \Lambda^2(V^*)$.  Then  $*^2=1$ and one has the splitting 
\be\label{hodge}\Lambda^2\,V^*=\Lambda^2_+\,V^*\oplus \Lambda^2_-\,V^*,\ee
where $\Lambda^2_\pm\,V^*$  are the $\pm 1$ eigenspaces of $*$, called   SD (self-dual) and the  ASD (anti-self-dual) 2-forms (resp.). 

  Let $\SO_{2,2}\subset \GL(V)$ be the corresponding orientation-preserving orthogonal 
group and $\so_{2,2}\subset \End\,V$ its Lie algebra. With respect to a null basis, the matrices of elements in $\so_{2,2}$  are of the form 
\be\label{sott}
\left(\begin{array}{cc}
A&B\\
C&-A^t
\end{array}\right), \quad A,B,C\in Mat_{2\times 2}(\R), \; B^t=-B, \; C^t=-C.
\ee

There is a natural isomorphism (equivalence of $\SO_{2,2}$-representations)
    \be\label{iso2}\so_{2,2}\iso\Lambda^2\,V^*,\quad T\mapsto {1\over 2}\<\,\cdot\, ,T\,\cdot\,\>
\ee
and  a Lie algebra decomposition $\so_{2,2}=\sl^+_2(\R)\oplus\sl^-_2(\R),$ given by
$$
{\footnotesize
 \left(\begin{array}{cc}
A&B\\
C&-A^t
\end{array}\right)=\left(\begin{array}{cc}
A_0&0\\
0&-A_0^t
\end{array}\right)+\left(\begin{array}{cc}
{\tr A\over 2}\II&B\\
C&-{\tr A\over 2}\II
\end{array}\right), \quad A_0=A-{\tr A\over 2}\II\in\sl_2(\R), 
}$$ 
matching  the decomposition of Eqn.~(\ref{hodge}), i.e. $\sl_2^\pm(\R)\iso \Lambda^2_\pm\,V^*.$

 Given a 2-plane $W\subset V$ pick a basis $\theta^1, \theta^2$  of the annihilator $W^0\subset V^*$ and  let $\beta=\theta^1\wedge\theta^2$. If we pick another basis of $W^0$ then $\beta $ is multiplied by a non-zero constant (the determinant of the matrix of change of basis), hence $\R\beta\subset \Lambda^2(V^*)$ is well-defined in terms of $W$ alone. This defines the {\em Pl\"ucker embedding}  of the grassmanian of 2-planes $Gr(2,V)\hookrightarrow\P(\Lambda^2\,V^*)\simeq \RP^5$. Its  image is given in homogeneous coordinates by the quadratic equation $\beta\wedge\beta=0$. We say that a 2-plane $W$ is SD (self-dual) if $\R\beta\subset \Lambda^2_+\,V^*$, and ASD (anti-self-dual)  if $\R\beta\subset \Lambda^2_-\,V^*$. We denote by 
 $$\T^+V:=\{W\subset V\st W\hbox{ is a SD 2-plane}\}.$$
Using the Pl\"ucker embedding,  $\T^+V$   is naturally  identified  with the conic in $\P(\Lambda^2_+\,V^*)\simeq\RPt$  given by the equations $\beta\wedge\beta=0,$ $*\beta=\beta.$ Similarly for the ASD 2-planes $\T^-V$.

A {\em null subspace}  is a subspace of $V$ on which the quadratic form $\<\,, \>$ vanishes. The maximum dimension of a null subspace is 2, in which case we call it a {\em null 2-plane}. It turns out that the null 2-planes are precisely the SD and ASD 2-planes.

\begin{proposition}\label{Nine} Let $V$ be an oriented 4-dimensional vector space  equipped with a quadratic form of signature $(2,2)$. Then 
\begin{enumerate}[leftmargin=*]\setlength\itemsep{5pt}
\item A 2-plane  $W\subset V$ is null if and only if it is  SD or ASD. Thus the space $Gr_0(2, V)$ of null 2-planes in $V$ is naturally identified with 
$$Gr_0(2,V)=(\T^+V)\sqcup(\T^- V),\quad \T^\pm V\simeq \RP^1.$$
\item 
Every  1-dimensional null subspace $N\subset V$  is the intersection of  precisely two null 2-planes, one SD and one ASD, $N=W^+\cap W^-.$  

\end{enumerate}
\end{proposition}

The proof is elementary (omitted). Let us just describe briefly the  picture that emerges from   the last assertion. The set of 1-dimensional null subspaces $N\subset V$  forms  the {\em projectivized null cone}  $\P C$,  a 2-dimensional quadric surface in $\P V\simeq\RP^3$, given in homogeneous coordinates, with respect to a null basis in $V$, by the equation $x^ax_a=0$. 
The statement then is that the SD and ASD null 2-planes in 
$V$ define a  {\em double ruling} of $\P C$. That is, the surface $\P C\subset \P V$, although not flat, contains many lines, forming a pair of foliations, so that through each point $e\in\P C$ pass exactly two lines, one from each foliation. The two lines through $e$ can also be found by intersecting $\P C$ with  the tangent plane to  $\P C$ at $e$. In some  affine chart, if $\P C$ is given by  $z=xy$ and  $e=(x_0,y_0, x_0y_0)$, then  the two null lines through $e$ are given by $z=x_0y$, $z=xy_0$. 
\begin{figure}[H]\centering
\includegraphics[width=0.3\textwidth]{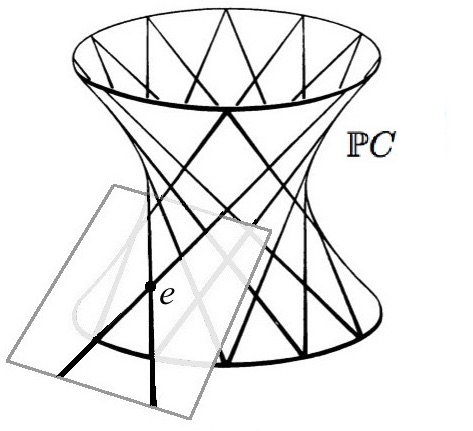}
\caption{The double-ruling of the projectivized null cone  $\P C\subset \P V$.}
\end{figure}

\subsubsection{The Levi-Civita connection and its curvature}\label{LC}

Now let $M$ be an oriented  smooth 4-manifold equipped with a pseudo-riemannian metric  $\met$ of signature $(2,2)$. 
Denote by $\BLA^k:=\Lambda^k(T^*M)$ the bundle of differential $k$-forms on $M$ and by  $\Gamma(\BLA^k)$ 
its space of smooth sections. In a (local) null coframe $\eta=(\eta^1, \eta^2, \eta_1, \eta_2)^t\in \Gamma(\BLA^1\otimes\R^4)$
the metric is  given by $\met=2\eta_a\eta^a$ and the Levi-Civita connection  is given by the unique 
 $\so_{2,2}$-valued  1-form $\Theta$ satisfying $d\eta+\Theta\wedge\eta=0$, i.e.  the connection is {\em torsion-free}. The associated covariant derivative is $\nabla \eta=-\Theta\otimes \eta$ 
and the curvature is the $\so_{2,2}$-valued 2-form  $\Phi=d\Theta+\Theta\wedge\Theta.$ The curvature form  $\Phi$ defines  via the isomorphism $\so_{2,2}\simeq \Lambda^2(T^*_mM)$ of  Eqn.~(\ref{iso2}) the curvature   {\em operator} $\cR\in\Gamma(\End(\BLA^2))$, which is  self-adjoint with respect to $\met$, i.e. $\cR^*=\cR$. Now  we use the  decomposition $\BLA^2=\BLA^2_+\oplus\BLA^2_-$ to    block decompose
\be\label{decompo}
\cR=\left(\begin{array}{cc}
\cA^+&\cB\\ 
\cB^*& \cA^-
\end{array}\right),\ee
where $\cB\in \mathrm{Hom}(\BLA^2_+, \BLA^2_-)$ and
$\cA^\pm\in\End\,\BLA^2_\pm$ are self-adjoint. This can be further refined into an  irreducible decomposition
$$\hbox{$\cR\sim  (\tr\cA^\pm,\,\cB,\,\cA^+-{1\over 3}\tr\cA_+, \,\cA^--{1\over 3}\tr\cA^-)$},$$
where $\tr\cA^+ =\tr\cA^- ={1\over 4}$ scalar curvature, $\cB$ is the traceless Ricci tensor and the last two components are traceless endomorphisms $\cW^\pm\in \Gamma(\End_0(\BLA^2_\pm))$,  defining the conformally invariant {\em Weyl tensor}, $\cW:= \cW^+\oplus\cW^-$  \cite{ST}. Thus the metric is {\em Einstein} iff $\cB = 0,$ {\em conformally flat} iff $\cW=0$, {\em self-dual} iff $\cW=\cW^+$ (i.e. $\cW^-=0$) and {\em anti-self-dual} iff $\cW=\cW^-$ (i.e. $\cW^+=0$). 

\subsubsection{Principal null 2-planes} Associated with  the Weyl tensor $\cW$ are  its {\em principal null 2-planes}, as follows. 
Recall  from Sect.~\ref{linalg} (just before Prop.~\ref{Nine}) that a 2-plane $ W\subset T_mM$  corresponds  to  a unique 1-dimensional space $\R\beta\subset \Lambda^2(T^*_mM)$ satisfying $\beta\wedge\beta=0$;  also,   $W$ is SD iff $\beta\in\BLA^2_+$, ASD iff $\beta\in\BLA^2_-$. 
\begin{definition}\label{def_principal}
A null 2-plane $W\subset  T_mM$  is {\em principal} if the associated  non-zero elements $\beta \in \Lambda^2(T^*_mM)$ satisfy $\beta\wedge\cW\beta=0.$ 

 \end{definition}
 If $\cW^+_m=0$ then  all SD null 2-plane in $T_mM$ are principal (by definition). Otherwise, the quadratic equation  $\beta\wedge\cW^+\beta=0$ defines a conic in $\P\Lambda^2_+(T^*_mM)\simeq \RP^2$, intersecting the conic $\T^+(T_mM)$ given by $\beta\wedge\beta=0$ in at most 4 points, corresponding precisely to  the principal SD 2-planes. The possible patterns of intersection  of these two conics  give rise to an algebraic classification of the SD Weyl tensor $\cW^+$, called the 
 {\em Petrov classification}. A similar classification holds for  $\cW^-$. 
\begin{figure}[h]
\centering
\includegraphics[width=1\textwidth]{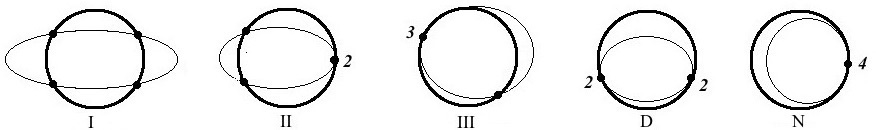}
\caption{The Petrov classification}\label{petrov}
\end{figure} 
 
\begin{rmrk} The above diagram depicts the classification over $\C$.  In the real case (such is ours) there are more sub-cases, as some of the intersection points might be complex. See for example \cite{GHN} for the complete classification. 
\end{rmrk} 

\subsubsection{The twistor fibration and distribution}\label{twist}

(We shall state the  results for the SD twistor fibration, but they apply verbatim  to the ASD case as well). 
Let $M$ be an oriented 4-manifold with a split-signature pseudo-riemannian metric, as in the previous subsection. The  {\em SD (self-dual) twistor fibration} is the fibre bundle  $$\RP^1\to\T^+M\to M$$  
 whose fiber at a point $m\in M$ is the set  $\T^+(T_mM)$ of  SD null 2-planes in $T_m M$ (see Prop.~\ref{Nine} of Sect.~\ref{linalg} above). 
 The total space $\T^+M$  is a 5-manifold  equipped with a natural rank 2 distribution $\D^+\subset T(\T^+M)$, the {\em SD twistor distribution},  defined by the Levi-Civita connection, as follows: a point $\tilde m\in\T^+_mM$  corresponds to a SD 2-plane $W\subset T_mM$; the 2-plane $\D^+_{\tilde m}\subset T_{\tilde m}(\T ^+M)$ is the horizontal lift of $W$  via  the Levi-Civita connection (one can check that $\D^+$ depends only on the conformal class $[\met]$  of the metric on $M$). By construction, the integral curves of $\D^+$ project to null-curves in $M$ with parallel self-dual tangent  2-plane. Conversely, each null curve in $M$ with parallel SD null 2-plane lifts uniquely to an integral curve of $(\T^+M, \D^+).$
 
   This  is the split-signature version  of the famous twistor construction of Roger Penrose \cite{Pen}.
A  standard feature  of  the twistor construction is  the relation between the integrability properties of $\D^+$ and  the vanishing of the SD Weyl  tensor $\mathcal W^+$. Namely,   $\D^+$ is integrable iff $\cW^+\equiv 0$ (i.e. $M$ is ASD). Less standard is the case of non-vanishing $\cW^+$, treated by An-Nurowski in \cite{AN}.

 \begin{theorem}[\cite{AN}]
 Let $(\T^+M,\D^+)$ be the SD twistor space and distribution of a split-signature oriented pseudo-riemannian conformal 4-manifold $(M,[\met])$ with a nowhere-vanishing SD Weyl tensor $\cW^+$. Then $\D^+$ is  $(2,3,5)$  away from the principal locus of  $\T^+M$. That is,  $\D^+$ is  $(2,3,5)$ when restricted to the open subset  $\T^+_*M\subset \T^+M$ obtained by removing the set of points corresponding to the principal SD 2-planes (at most 4 points on each fiber of $\T^+M\to M$; see Def.~\ref{def_principal} above).
 \end{theorem} 
See the theorem   in \cite{AN}, right before Corollary 1.

\subsection{The tangent SD 2-plane along a null curve in the dancing space.}\label{tang_SD}
Now we return to our case of $\M\subset \RPt\times \RPts$ equipped with the dancing metric $\met$, as defined in Prop.~\ref{dancing}.

\begin{definition}\label{ng}
Let $\Gamma$  be a parametrized curve in $\M$,  $\Gamma(t)=(q(t),p(t)) $. 
Then  $\Gamma $ is  {\em non-degenerate} if $q(t),p(t)$ are regular curves in $\RPt, \RPts$ (resp.); i.e., $q'(t)\neq 0$ and $ p'(t)\neq 0$ for all $t$. 
\end{definition}

Note that the non-degeneracy condition  is reparametrization independent, hence it applies to unparametrized curves $\Gamma\subset M$ (1-dimensonal submanifolds). It means that $\Gamma$ 
is nowhere tangent to the leaves of the double fibration $\RPt\leftarrow\M\to\RPts$.   
Equivalently, the projections of $\Gamma$ to $\RPt$ and $\RPts$ are non-singular.

 Now let $\Gamma$ be a null-curve in $ (\M,[\met])$. Then, by Prop. \ref{Nine}, there are two tangent null 2-plane fields defined along $\Gamma$, one SD and the other ASD, whose  intersection is the tangent line field along $\Gamma$. 
\begin{theorem}\label{eleven}
Every  integral curve $\tG$  of $(\Q,\D)$ projects to a non-degenerate  null-curve $\Gamma$  in  $\M$ with a parallel SD tangent 2-plane. Conversely, every non-degenerate null-curve in $(\M, [\met])$ with parallel SD tangent 2-plane lifts uniquely to an  integral curve  of $(\Q,\D)$. \end{theorem}

\begin{theorem}\label{ident} For each $(\q,\p)\in\Q$, the 2-plane 
$$\Pi_*\D_{(\q,\p)}\subset T_{(q,p)}\M,$$
where $(q,p)=\Pi((\q, \p)),$
is a non-principal self-dual 2-plane. The resulting map 
$$ \Q\to \T^+\M,\quad (\q,\p)\mapsto \Pi_*\D_{(\q,\p)},$$ is an  $\SLth$-equivariant  embedding, identifying    $\Q$ with  the non-principal locus of $\D^+$ in $\T^+\M$, and mapping $\D$ over to $\D^+$. 
\end{theorem}

The proofs of these two theorems will be carried out in the next subsection, using  the Maurer-Cartan structure equations of $\SLth$.

\subsection{Proofs of Theorems~\ref{eleven} and \ref{ident}}\label{proofs}

 Let $\{e_1,e_2,e_3\}$ be the standard basis of $\Rt$ and $\{e^1,e^2,e^3\}$ the dual basis of $\Rts$. Recall (Sect.~\ref{TwoFour}) that  $G=\SLth$ acts transitively on $(\Q,\D)$ and $(\M,\met)$ by $g\cdot(\q, \p)=(g\q, \p g^{-1}),$ $g\cdot([\q], [\p])=([g\q], [\p g^{-1}])$,  preserving  $\D$ and $\met$ (resp.). 
Fix $\tilde m_0=(e_3, e^3)\in \Q$ and  $m_0=\Pi(\tilde m_0)=([e_3], [e^3])\in\M$. Define
\be\label{principal}
    \xymatrix{
        G  \ar[dr]_j \ar[r]^\tj & \Q \ar[d]^\Pi\\
                               & \M }
\ee
by $\tj(g)=g\cdot \tilde m_0=(ge_3, e^3g^{-1})$, $ j(g)=g\cdot m_0=([ge_3],[e^3g^{-1}])=(\Pi\circ\tj)(g).$ Then  $j$ is a principal $H$-fibration and $\tj$ a principal $H_0$-fibration, 
where 
$$H=\left\{ \left(
\begin{array}{cc}
A&0\\
0 &
a^{-1}
\end{array}\right)\st \; A\in \GL_2(\R), \; a=\det(A)\right\}
\simeq\GL_2(\R)$$
is the stabilizer subgroup of $m_0$, with Lie algebra
\be\label{lieh}\h=\left\{ \left(
\begin{array}{cc}
X&0\\
0 &
-x
\end{array}\right)\st \; X\in\gl_2(\R), \; x=\tr(X)\right\}\simeq\gl_2(\R)
\ee
and 
$$H_0=\left\{ \left(
\begin{array}{cc}
A&0\\
0 &
1
\end{array}\right)\st \; A\in \SL_2(\R)\right\}
\simeq\SL_2(\R)$$
is the stabilizer subgroup of $\tilde m_0$, with Lie algebra
$$\h_0=\left\{ \left(
\begin{array}{cc}
X&0\\
0 &
0
\end{array}\right)\st \; X\in \sl_2(\R)\right\}\simeq\sl_2(\R).$$

 The left-invariant MC (Maurer-Cartan) form on $G=\SLth$ is the $\g$-valued 1-form  
$\omega=(\om{i}{j}):=g^{-1}dg$, i.e. $\tr(\omega)=\om{i}{i}=0$,  $i,j\in\{1,2,3\}$ (using, as always,  the summation convention on repeated indices). The components of   $\omega$ provide a global coframing on $G$, whose basic properties (immediate from its definition) are
\be\label{MC} \begin{array}{ll}
(a)\quad \omega_e=id_\g&\\
(b)\quad (L_g)^*\omega=\omega & \hbox{(left invariance)}\\
(c)\quad (R_g)^*\omega=g^{-1}\omega g & \hbox{(right $\Ad$-equivariance)}\\
(d)\quad d\omega =-\omega \wedge\omega& \hbox{(the MC structure  equation).}
\end{array}
\ee

Now let us rename  the components of $\omega$: 
\be\label{many}
\eta^a:=\om{a}{3}, \; \eta_b:=\om{3}{b}, \;\phi:=\om{a}{a}=-\om{3}{3},\;\theta^{a}_{\,b}:=\om{a}{b}+\delta^a_{\,b}\phi, \quad a,b\in\{1,2\}.
\ee

Furthermore, introduce the matrix notation
\be\label{compon}
\eta:=
\left(\begin{matrix}\eta^1\\ \eta^2\\ \eta_1\\\eta_2\end{matrix}\right), 
\Theta:=
\left(\begin{array}{cc}
\theta&0\\
0&-\theta^t
\end{array}\right) , 
\theta:=(\theta^{a}_{\,b})=
\left(\begin{array}{cc}
2\om{1}{1}+ \om{2}{2}&\om{1}{2}\\
\om{2}{1}&\om{1}{1}+ 2\om{2}{2}
\end{array}\right).
\ee
With this notation, Eqn.~(\ref{MC}c) now  reads
\be\label{comp}(R_h)^*\eta=\rh^{-1}\eta,\quad (R_h)^*\Theta=\rh^{-1}\Theta\rh,\quad h\in H,
\ee
where $\rho:H\to \SO_{2,2}$ is the {\em isotropy representation},   
\be\label{isotropy}h=\left(
\begin{array}{cl}
A&0\\ 0& a^{-1}
\end{array}
\right)
\mapsto 
\rh=\left(
\begin{array}{cc}
aA&0\\ 0& (aA^t)^{-1}
\end{array}
\right),\quad A\in\GL_2(\R), \quad a=\det A.
\ee
The MC structure  equation (Eqn.~(\ref{MC}d)) also  breaks into two equations,
\be\label{break}
d\eta+\Theta\wedge\eta=0,\qquad
d\Theta+\Theta\wedge\Theta=\left(\begin{array}{cc}\varphi&0\\ 0&-\varphi^t\end{array}\right), 
\ee
where
\be\label{curve}
\varphi:=d\theta+\theta\wedge\theta=
\left(\begin{array}{cc}
2\eta_1\wedge\eta^1+ \eta_2\wedge\eta^2&\eta_2\wedge\eta^1\\[5pt]
\eta_1\wedge\eta^2&\eta_1\wedge\eta^1+ 2\eta_2\wedge\eta^2
\end{array}\right).
\ee

From Formula~(\ref{lieh}) for $\h$, we see that the  four 1-forms $\eta^a, \eta_b\in\Omega^1(G)$
 are pointwise linearly independent and $j$-horizontal, i.e. vanish  on the fibers of 
 $j:G\to M$, hence span $j^*(T^*M)\subset T^*(G).$ Similarly, $\eta^a, \eta_b, \phi$ span  $\tj^*(T^*\Q).$

\begin{proposition}\label{all}    

Consider the principal  fibrations $j, \tj$ of Eqn.~\eqref{principal} and the left-invariant 1-forms $\eta, \phi, \theta, \Theta, \varphi$ on $G$, as defined above in Eqns.~\eqref{many}-\eqref{curve}. 

\begin{enumerate}[leftmargin=18pt,label=(\alph*)]\setlength\itemsep{5pt}
\item  $j^*\met=2\eta_a\eta^a$, where $\met$ is the dancing metric on $\M$, as defined in Prop.~\ref{dancing}. 

\item Let $\nabla$ be the covariant derivative on $T^*M$ associated with   the Levi-Civita connection of   $\met$ and $\wnab=j^*(\nabla)$ its pull-back to $j^*(T^*M)$.Then $\wnab\eta^a=-\theta^a_{\,b}\otimes\eta^b, $ $\wnab\eta_b=\theta^a_{\,b}\otimes\eta_a,$ or in matrix form,
$\wnab\eta=-\Theta\otimes\eta$. The associated  curvature 2-form   is $\Phi:=d\Theta+\Theta\wedge\Theta$,  given in terms of $\eta$ by   Eqns.~\eqref{break})-\eqref{curve} above.

\item  Let $vol\in\Omega^4(M)$ be the positively oriented unit volume form on $\M$ (see  Sect.~\ref{orientation}). Then  $j^*(vol)=\eta^1\wedge\eta^2\wedge\eta_1\wedge\eta_2.$

\item   Let  $\D\subset T\Q$ be the rank 2 distribution given by $d\p=\q\times d\q$ and $\D^0\subset T^*\Q$ its annihilator.  
Then   $\tj^*(\D^0)=\Span\{\eta^2-\eta_1,\;\eta^1+\eta_2, \; \phi\}.$

\end{enumerate}
\end{proposition}

\newcommand{\hT}{\hat\Theta}
\begin{rmrk} We can rephrase the above in terms of coframes on $M$
 and $Q$, as follows: let $\sigma$ be a local section of $j:G\to M$, 
 then (a) $\weta=\sigma^*\eta$ is a null-coframe on $M$, 
 so that $\met=2\weta_a\weta^a$, 
 (b)  $\hT=\sigma^*(\Theta)$ is the connection 1-form of the 
 Levi-Civita connection of $\met$ with respect to the coframe $\weta$, 
 (c) $vol=\weta^1\wedge\weta^2\wedge\weta_1\wedge\weta_2.$ and
  (d) $\D=\Ker\{\tilde\eta^2-\tilde\eta_1,\;\tilde\eta^1+\tilde\eta_2, \; \tilde\phi\},$
   where $\tilde\eta=\tilde\sigma^*\eta,$ $\tilde\phi=\tilde\sigma^*\phi $ and $\tilde\sigma= \sigma\circ\Pi$ (a local section of $\tj:G\to Q$). 
\end{rmrk}
\begin{proof} 
 (a) First,  the formula $(R_h)^*\eta=\rh^{-1}\eta$ of Eqn.~(\ref{comp}) implies  that $\eta_a\eta^a$, a $G$-left-invariant $j$-horizontal symmetric 2-form   
 on $G$, is  $H$-right-invariant, hence descends to a 
well-defined $G$-invariant symmetric 2-form on $M.$ Next, by examining the isotropy
 representation of $H$ (Eqn.~(\ref{isotropy})), one sees that $T_{m_0}M$ admits  a unique 
 $H$-invariant quadratic form, up to a constant multiple, hence $M$ admits a unique $G$-invariant 2-form, up to a constant multiple.  It follows that it is enough to verify the equation $j^*\met=2\eta_a\eta^a$ 
 on a single non-null element $Y\in \g=T_eG$; for example,  $Y=Y_1$ from  the proof of Prop.~\ref{prop235}. We omit this (easy) verification.
  
\mn(b) The relations   $d\eta+\Theta\wedge\eta=0$, $(R_h)^*\Theta=\rh^{-1}\Theta\rh$  and the formula for $\Theta$
(Eqns.~(\ref{compon})-(\ref{break}))  show  that $\Theta$ is an $\so_{2,2}$-valued 1-form on $G$, descending to  a torsion-free $\SO_{2,2}$-connection  on $T^*M$, hence  is in fact the Levi-Civita connection of $\met$.

\mn(c)  First one verifies that $\eta^1\wedge\eta^2\wedge\eta_1\wedge\eta_2$ is  a volume  form of norm 1 with respect to $2\eta_a\eta^a$. Then, to compare to the orientation definition of Sect.~\ref{orientation}, we check that $K^*\eta^a=\eta^a,$ $K^*\eta_b=-\eta_b,$ hence $\eta^1+\eta_1, \eta^2+\eta_2, \eta^1-\eta_1, \eta^2-\eta_2$ is a para-complex coframe. Now one calculates  $(\eta^1+\eta_1)\wedge(\eta^2+\eta_2)\wedge(\eta^1-\eta_1)\wedge(\eta^2-\eta_2)=4\eta^1\wedge\eta^2\wedge\eta_1\wedge\eta_2,$ hence $\eta^1, \eta^2, \eta_1, \eta_2$ is a positively oriented coframe.

\mn(d)  Let $E_j:G\to \Rt$ be the function that assigns  to an element $g\in G$ its   $j$-th  column, $j=1,2,3$. Then $\omega=g^{-1}dg$ is equivalent to $dE_j=E_i\om{i}{j}$. Next let $E^i:G\to \Rts$ be the function  assigning to $g\in G$ 
the $i$-th row of $g^{-1}$. Then clearly $E^iE_j=\delta^i_{\,j}$ (matrix multiplication of a row by column vector), and by taking exterior derivative of the last equation  we obtain $dE^i=-\om{i}{j}E^j.$ Also, $\det(g)= 1$ implies $E_i\times E_j=\epsilon_{ijk}E^k$, $E^i\times E^j=\epsilon^{ijk}E_k$. Next, by definition of $\tj$, $E_3=\q\circ\tj,$ $E^3=\p\circ\tj$. Now we calculate 
$\tj^*(d\p-\q\times d\q)=dE^3-E_3\times dE_3=-\om{3}{j}E^j-(E_3\times E_i)\om{i}{3}=(\eta^2-\eta_1)E^1-(\eta^1+\eta_2)E^2+\phi E^3.$
\end{proof}

\begin{cor}[Proofs of Thms.~\ref{eleven} and \ref{ident}]\label{grancor}
\n\begin{enumerate}[leftmargin=18pt,label=(\alph*)]\setlength\itemsep{5pt}
\item $(\M,\met)$ is Einstein but not Ricci-flat, SD (i.e. $\cW^-\equiv 0$) and $\cW^+$ is nowhere vanishing, of Petrov type D (see 
Fig.~\ref{petrov}). 
More precisely,  at each $(q,p)\in M$ there are exactly two principal SD null 2-planes, 
each of multiplicity 2, given by $T_q\RPt\oplus \{0\}$ and $\{0\}\oplus T_p\RPts. $

\item Every integral curve $\tG$ of $(Q,\D)$ projects to a non-degenerate null curve $\Gamma:=\Pi\circ\tG$ in  $(\M,[\met])$ with parallel SD tangent 2-plane. 

\item Every non-degenerate null curve $\Gamma$ in $(\M,[\met])$ with parallel SD tangent 2-plane lifts uniquely to an integral curve $\tG$ of $(\Q, \D).$

\item For every $\tilde m\in \Q$,  $\Pi_*\D_{\tilde m}\subset T_{\Pi(\tilde m)}\M$ is a non-principal SD null 2-plane.

\item Let $\T^+_*M\subset \T^+M$ be the non-prinicipal locus (the complement of the principal points). Then the map $\nu:\Q\to \T^+_*M$, 
$\tilde m\mapsto \Pi_*\D_{\tilde m}$, 
is  an $\SLth$-equivariant diffeomorphism, mapping $\D$ unto $\D^+$.

\end{enumerate}

\end{cor}

\begin{proof}

\mn (a)   By Prop.~\ref{all}a and \ref{all}c, the coframe $\eta^1, \eta^2,\eta_1, \eta_2$ is null and positively oriented. It follows from the definition of the Hodge dual  that
\begin{subequations}
\begin{align}\label{iso}
\quad j^*( \Lambda^2_+\,M)&= 
\Span\{ \eta_1\wedge\eta^1+\eta_2\wedge\eta^2, \;
\eta^1\wedge\eta^2, \; 
\eta_1\wedge\eta_2
\},\\
\quad j^*( \Lambda^2_-\,M)&= 
\Span\{ \eta_1\wedge\eta^1-\eta_2\wedge\eta^2, \;
\eta_1\wedge\eta^2, \; 
\eta_2\wedge\eta^1
\}.
\end{align}
\end{subequations}

Then using  the formula for the curvature form $\Phi$ (Eqns.~(\ref{break})-(\ref{curve})) and the definition of the curvature operator  $\cR$
(Sect.~\ref{LC}), one finds that $j^*(\cR)$ is diagonal in the above bases,  with matrix 
$$
j^*(\cR)={\footnotesize\left(\begin{array}{ccc|ccc}
-3&&&&&\\ 
&\;0&&&&\\ 
&&\;0&&&\\
\hline
&&&-1&&\\ 
&&&&-1&\\ &&&&&-1
\end{array}\right)}.
$$ 
Comparing this expression with the the decomposition of $\cR$ of Eqn.~(\ref{decompo}), we see that the dancing metric is Einstein ($\cB\equiv 0$), the scalar curvature is  $-12$,  $\cW^-\equiv 0$ and 
$$j^*(\cW^+)={\footnotesize\left(
\begin{matrix}-2&&\\ &1&\\ &&1\end{matrix}\right)}.
$$
Now let  $a,b,c$  be the  coordinates dual to the basis of  $ j^*( \Lambda^2_+\,M)$  of Eqn.~(\ref{iso}). Then $\beta\wedge\beta=0$
 is given by $a^2-bc=0$ and $\beta\wedge\cW^+\beta=0$  by $2a^2+bc=0$. This system of  two homogeneous equations has two non-zero
  solutions (up to a non-zero multiple), $a=b=0$ and $a=c=0$, each with multiplicity 2 (the  pair of conics defined in
  each fiber of $\P\Lambda^2_+M$ by these equations are tangent at their two  intersection points). 
  The corresponding SD  2-forms are $\eta_1\wedge\eta_2, \eta^1\wedge\eta^2$,  corresponding  
  to the principal SD null 2-planes  $T_q\RPt\oplus \{0\}$ and $\{0\}\oplus T_p\RPts $ (resp.), as claimed. 

\mn(b)  Let $\tG(t)=(\q(t),\p(t))$  be  a regular parametrization of an integral curve of $(\Q, \D)$, i.e.  $\tG'=(\q', \p')$  is nowhere vanishing and  $\p'=\q\times \q'$. We  first show that $\Gamma=\Pi\circ\tG$ 
  is non-degenerate. Let $\Gamma(t)=\Pi(\tG(t))=(q(t), p(t)),$ 
  where $q(t)=[\q(t)],$ 
  $p(t)=[\p(t)].$ We need to show that 
 $q',  p'$ are nowhere vanishing.

\begin{lemma}\label{dual} The distribution $\D\subset T\Q$, given by $d\p=\q\times d\q$, is also given by $d\q=-\p\times d\p.$ 
\end{lemma}
\begin{proof} Let  $\D'=\Ker(d\q+\p\times d\p)\subset T\Q$. Then both $\D,\D'$ are $\SLth$-invariant, hence it is enough to compare them at say $(e_3, e^3)\in\Q$. At this point $\D$ is given by $dp_1+dq^2=dp_2-dq^1=dp_3=dp^3+dq_3=0$, and $\D'$ by $dq^1-dp_2=dq^2+dp_1=dq^3=dp^3+dq_3=0.$ These obviously have the same 2-dimensional space of  solutions. 
\end{proof}

Now   $q'=\q'\,(\mod \q)$, hence  $q'=0\ent  \q'\times\q=0\ent \p'=0,$ so by Lemma \ref{dual}, $\q'=-\p\times\p'=0$. Similarly, $p'=0\ent \q'=\p'=0,$ hence $\Gamma$ is non-degenerate.

 Next we show that $\Gamma$ is null. Let $\sigma$ be a lift of $\tG$ (hence of  $\Gamma$) to $G=\SLth$. Let 
$\sigma^*\eta^a=s^adt, $ $\sigma^*\eta_b=s_bdt,$ $a,b=1,2, $ for some  real-valued functions (of $t$) $s^1, s^2, s_1, s_2$.  Then, by Prop.~\ref{all}a and \ref{all}d, $\met(\Gamma', \Gamma')=2s_as^a=
2(s_1 (-s_2)+s_2s_1)=0,$ hence $\Gamma$ is a null curve. 

 Next we show that  the SD null 2-plane along $\Gamma$  is parallel.
Let  $\weta=\eta\circ\sigma$ be the  coframing of $\Gamma^*(TM)$ determined by the lift $\sigma$ of $\Gamma$ (a ``moving coframe"  along 
$\Gamma$).

\begin{rmrk} $\weta$  should not be  confused with
 $\sigma^*\eta=(s^1, s^2, s_1, s_2)^tdt$, the restriction of $\weta$ to $T\Gamma$.   
 \end{rmrk}
 
 Let  $W$ be  the 2-plane field along 
 $\Gamma$ defined by  $\weta^1+\weta_2=\weta^2-\weta_1=0$.  By Prop.~\ref{all}d, 
 $\sigma^*(\eta^1+\eta_2)=\sigma^*(\eta^2-\eta_1)=0,$ hence $W$ is tangent to
 $\Gamma$. The  2-form corresponding to $W$ is 
 $\beta=(\weta^1+\weta_2)\wedge(\weta^2-\weta_1)=
 \weta^1\wedge\weta^2+\weta_1\wedge\weta_2+\weta_1\wedge\weta^1+\weta_2\wedge\weta^2$, 
 which is SD by  Formula~(\ref{iso}),  hence $W$ is the  SD tangent 2-plane field along $\Gamma$.
  Now a short calculation, using Prop.~\ref{all}b, shows that
  $$\nabla \beta =3(\sigma^*\phi)\otimes (\weta_1\wedge\weta_2-\weta^1\wedge\weta^2).$$ By Prop.~\ref{all}d, $\sigma^*\phi=0$, hence $\nabla \beta=0$, so $W$ is parallel. 

\mn(c)  Let $\sigma$ be a lift of $\Gamma$ to $G$, with 
$\sigma^*\eta^a=s^adt, $ $\sigma^*\eta_b=s_bdt.$ 

\begin{lemma}\label{adapted} Given a non-degenerate parametrized null curve $\Gamma:\R\to M$, there exists a lift $\sigma$ of $\Gamma$ to $G$ such that  $s^1=s_2=0,$ $s_1=s^2=1$. In other words, 
$$\sigma^*\omega=
\left(\begin{array}{ccc} *&*&0\\ *&*&1\\ 1&0&*\end{array}\right)dt.
$$
\end{lemma}
 
\begin{rmrk} We call such a lift $\sigma$ {\em adapted} to $\Gamma$.  
\end{rmrk}
 
\begin{proof} Starting with an  arbitrary lift $\sigma$, any other lift is of the form 
$\bar\sigma=\sigma h$, where $h:\R\to H$ is an arbitrary $H$-valued smooth function, i.e. 
$$h= 
\left(
\begin{array}{cl}
A&0\\
0 &
a^{-1}
\end{array}\right), \quad  A:\R\to \GL_2(\R), \; a=\det(A).
$$
 Now a short calculation shows that 
 $$\bar\sigma^*\omega=
 h^{-1}(\sigma^*\omega)h+h^{-1}dh=
\left(\begin{array}{ccc} *&*&0\\ *&*&1\\ 1&0&*\end{array}\right)dt
$$ 
provided 
\be\label{B}
aA \left(\begin{matrix} 0\\  1\end{matrix}\right)=
 \left(\begin{matrix}  s^1\\ s^2\end{matrix}\right), \quad 
(1, \; 0)=
( s_1, \;  s_2)aA.
\ee
Now one checks that the  last system of equations can be solved for $A$  iff  $(s^1,  s^2)\neq 0,$ $(s_1, s_2)\neq 0$ 
and $s_as^a=0$. These are precisely the non-degeneracy and nullity conditions on $\Gamma$. From  $A$ we obtain  $h$ and  the desired $\bar\sigma$.
\end{proof}
Once we have an adapted lift $\sigma$ of $\Gamma$, with associated moving coframe $\weta:=\eta\circ\sigma$, we define a 2-plane field $W$  along $\Gamma$ 
  by   
 $\weta^1+\weta_2=\weta^2-\weta_1=0$. Then 
 $\sigma^*(\eta^1+\eta_2)=(s^1+s_2)dt=0,$ $\sigma^*(\eta^2-\eta_1)=(s^2-s_1)dt=0,$ hence $W$ is tangent to $\Gamma$. Let $\beta:=(\weta^1+\weta_2)\wedge(\weta^2-\weta_1).$ 
 Then $\beta$ is SD, so  $W$ is the SD tangent 2-plane along 
 $\Gamma$. Now  $W$ is  parallel  
 $\ent \nabla \beta =3(\sigma^*\phi)\otimes (\weta_1\wedge\weta_2-\weta^1\wedge\weta^2)\equiv 0\,(\mod \beta)\ent \sigma^*\phi=0$, since $\weta_1\wedge\weta_2-\weta^1\wedge\weta^2$ is a non-zero ASD  form, hence  $\not\equiv 0\,(\mod\beta).$ It follows that $\sigma$ satisfies $\sigma^*(\eta^1+\eta_2)=
\sigma^*(\eta^2-\eta_1)= \sigma^*(\phi)=0$, hence, by Prop.~\ref{all}d, $\tG:=\tj\circ\sigma$ is a lift of $\Gamma$ to an integral curve of $(Q, \D)$.

To show uniqueness, if $\tG(t)=(\q(t), \p(t))$ then  any other lift of $\Gamma$ to $\Q$ is of the form $(\lambda(t)\q(t), \p(t)/\lambda(t))$ for some non-vanishing real function $\lambda(t)$. If this other lift is also an integral curve of $(\Q, \D)$ then 
$(\p/\lambda)'-(\lambda\q)\times( \lambda\q)'=-(\lambda'/\lambda^2)\p+(1/\lambda-\lambda^2)\p'=0.$ 
Multiplying the last equation by $\q$ and using $\p\q=1, \p'\q=0$, we get $\lambda'=0\ent (1-\lambda^3)\p'=0$. Now $\p'\neq 0$ since $\Gamma$ is non-degenerate $\ent \lambda^3=1\ent \lambda= 1.$

\mn(d)  Let $m=\Pi(\tilde m)$,  $W=\Pi_*(\D_{\tilde m})\subset T_mM$, $g\in G$ such that $\tj(g)=\tilde m$  and  $\weta=\eta(g)$ the corresponding coframing of $T_mM$. Then, by Prop.~\ref{all}d,  $W=\Ker\{\weta^1+\weta_2, \weta^2-\weta_1\}.$ As before (item (b)),  one checks that $\beta:=(\weta^1+\weta_2)\wedge(\weta^2-\weta_1)=\weta^1\wedge\weta^2+\weta_1\wedge\weta_2+\weta_1\wedge\weta^1+\weta_2\wedge\weta^2$ is SD $\ent W$ is SD, but not principal (the SD 2-planes are given by $\weta^1\wedge\weta^2$ and $\weta_1\wedge\weta_2$; see the proof of item (a) above).

\mn(e)  One checks that $\nu$ is  $\SLth$-equivariant, $\Q$ and $\T^+_*M$ are $\SLth$-homogeneous manifolds, with the same stabilizer at $\tilde m_0=(e_3, e^3)$ and $\nu(\tilde m_0)$, hence $\nu$  is a diffeomorphism. It remains to show that   $\nu_*\D=\D^+$. This is just a reformulation  of items (b) and (c) above. 
\end{proof}

\section{Projective geometry:  dancing pairs and projective rolling}

We give here two related projective geometric interpretations of the Cartan-Engel distribution $(\Q,\D)$: ``dancing pairs" and ``projective rolling".   We start in Sect.~\ref{sec_danc} with the {\em dancing condition},  characterizing   null curves in $(\M, [\met])$. Next  in Sect.~\ref{simple} we use this characterization for an elementary derivation of an explicit  coordinate formula for $[\met]$. In Sect.~\ref{ss_cr} we give yet another formula for the dancing metric $\met$, this time in terms of the {\em cross-ratio} (a classical projective invariant of 4 colinear points). This is followed in Sect.~\ref{pj} by a study of the relation between the {\em projective structures} of the members of a dancing pair (the structures  happily match up), which we use in Sect.~\ref{sec:mate} for deriving the ``dancing mate equation". 
To illustrate all these concepts we study two examples: the  ``dancing mates of the circle" (Sect.~\ref{sec:circle}) and ``dancing pairs with constant projective curvature"  (Sect.~\ref{const}). 

We mention also in Sect.~\ref{sec:centro} a curious geometric interpretation for Eqns.~\eqref{eqns} that we found  during the proof of Prop.~\ref{TwentyOne}: curves in $\R^3$ with constant ``centro-affine torsion''.

The rest of the section (Sect.~\ref{RT}) is dedicated to {\em projective rolling}. Our motivation comes from the intrinsic geometric formulation of ordinary (riemannian) rolling, as appears in \cite{BrHs}. After making the appropriate definitions,  the nullity condition for curves on $(\M, [\met])$ becomes  the ``no-slip" condition for the projective rolling of $\RPt$ along $\RPts$, self-dual null 2-planes become ``projective contact elements" of the two surfaces and the condition of ``parallel self-dual tangent 2-plane" is the ``no-twist" condition of projective rolling, expressed in terms of the osculating conic of a plane curve and its developments, as appear in \'E.~Cartan's book \cite{Cbook}.

\subsection{Projective duality and the dancing condition}\label{sec_danc}

Let $\RP^2:=\P(\R^3)$ be the  {\em  real projective plane}, i.e. the space of 1-dimensional linear subspaces in $\R^3$,  with   $$\pi:\R^3\setminus\{ 0\}\to \RP^2,\quad \q\mapsto \R\q,$$
the canonical projection. If $q=\pi(\q)\in\RPt$, where $\q=(q^1,q^2,q^3)^t\in\R^3\setminus\{0\}$, we write $q=[\q]$ and say that  $q^1,q^2,q^3$ are the {\em homogeneous coordinates} of $q$. 
Similarly, $\RPts:=\P(\Rts)$ is the {\em dual projective plane}, with $$ \quad \bar\pi:\Rts\setminus 0\to \RPts, \quad \p\mapsto \R\p,$$ the canonical projection,  $\bar\pi(\p)=[\p]$. 

 A {\em projective line} in $\RPt$ is the projectivization of a 2-dimensional linear subspace in $\R^3$, 
i.e.~the set of 1-dimensional subspaces of $\R^3$ contained in a fixed 2-dimensional linear 
subspace of $\R^3$. The space of  projective lines in $\RPt$ is 
naturally identified with $\RPts$;
to each $p=[\p]\in\RPts$ corresponds the {\em dual projective line}
 $\hat p\subset \RP^2$, the projectivization of  
$\p^0=\{\q\in\R^3|\p\q=0\}$, and each projective line in 
$\RPt$ is of this form. 
Similarly, $\RPt$ is 
naturally identified with the space of projective lines in $\RPts$; to each 
point $q=[\q]\in\RPt$ corresponds the dual  projective line $\hat q\subset \RPts$, 
the projectivization of the 2-dimensional subspace $\q^0=\{\p\in\Rts|\p\q=0\}$. 

We say that $(q,p)\in \RPt\times\RPts$ 
are {\em incident} if $q\in\hat p$ (same as $p\in \hat q).$ We  also write this condition as $q\in p$. 
In homogeneous coordinates this is simply $\p\q=0$.

 Given a smooth  curve $\gamma\subset \RP^2$ (a 1-dimensional submanifold), the {\em duality map}  $*:\gamma\to \RPts$ assigns to each point $q\in \gamma$ its tangent line $q^*\in\RPts$. The image of $\gamma$ under the duality map is the {\em  dual curve}  $\gamma^*\subset \RPts$. In homogeneous coordinates, if $\gamma$ is parametrized by $q(t)=[\q(t)]$, where  $\q(t)\in\R^3\setminus \{0\}$, then $\gamma^*$ is parametrized by $q^*(t)=[\q^*(t)]$, where $\q^*(t):=\q(t)\times \q'(t)\in\Rts\setminus \{0\}.$ If $\gamma$ is a smooth curve without inflection points (points where $q''\equiv 0 \, \mod q'$, see  Def.~\ref{def:inflect} below) then $\gamma^*$ is smooth as well. More generally,  inflection points of $\gamma$  map to singular (or ``cusp") points 
of $\gamma^*$, where $(q^*)'=0$. 

 Similarly, given a curve $\bar \gamma\subset \RPts$, the duality map $\bar \gamma\to \RP^2$
 assigns to each line $p\in\bar \gamma$  its {\em turning point} $p^*$.
 In homogeneous coordinates:  if $\bar \gamma$ is parametrized by $p(t)=[\p(t)]$, 
then its dual   $\bar \gamma^*\subset \RPt$ is parametrized by $p^*(t)=[\p^*(t)]$, 
where $\p^*(t)=\p(t)\times\p'(t).$

Geometrically, $\bar \gamma$ is a 1-parameter
 family of lines  in $\RPt$, and its dual $\bar \gamma^*$ is the {\em envelope} of the family. Using the above formulas 
 for the duality map,  it is easy to verify that, away from inflection points, $(\gamma^*)^*=\gamma$ ; that is, $p(t)$ is the tangent line to $\bar \gamma^*$ at $p^*(t)$.

\begin{figure}[h!]\centering
\includegraphics[width=0.5\textwidth]{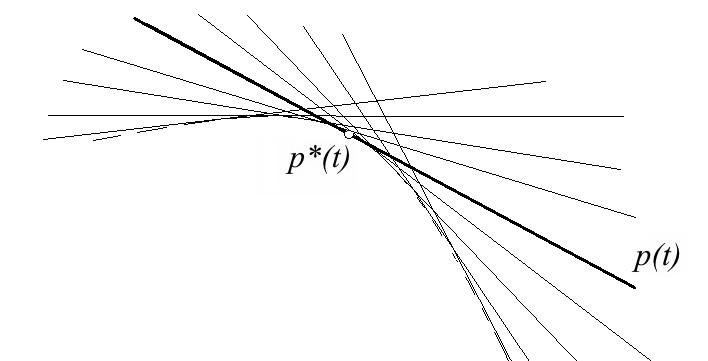}
\caption{The envelope  of a family  of lines}
\end{figure}

\begin{definition}\label{def_danc} A pair of  parametrized curves $q(t), p(t)$ in $ \RPt,  \RPts$ (resp.) satisfies the {\em  dancing condition}  if for each $t$ 

\begin{enumerate}
\item $(q(t), p(t))$ is non-incident;  

\item  if $q'(t)\neq 0$ and $p'(t)\neq 0$ then the tangent line $q^*(t)$ at $q(t)$  is incident to the turning point $p^*(t)$ of $p(t)$.

 \end{enumerate}

\end{definition}
See Fig.~\ref{fig:dance} of Sect.~\ref{ProjGeom}. 

\begin{rmrk} In condition (2), if either $q'$ or $p'$ vanish, then $q^*$ or $p^*$ is not well-defined, in which case, by definition, the curves satisfy the dancing condition. 
 \end{rmrk}

\begin{proposition}\label{prop_danc} The following conditions on a parametrized curve $\Gamma(t)=(q(t), p(t))$ in $\M$ are equivalent:

\begin{enumerate}
\item  $\Gamma(t)$ is a null curve in $(\M,[\met])$;

\item  the pair of curves $q(t), p(t)$ satisfies the dancing condition.
\end{enumerate}
\end{proposition}

\begin{proof} We use the notation of the proof of Prop.~\ref{all}d. Let  $\sigma$ be a lift of $\Gamma$ to $G=\SLth$ (see diagram~(\ref{principal}) of  Sect.~\ref{proofs}). Then 
$(\q(t),\p(t)):=(E_3(\sigma(t)), E^3(\sigma(t))$, $\sigma^*\eta^a=s^adt,$ $\sigma^*\eta_b=s_b dt.$
If  either $q'$ or $p'$ vanish, then either $(s^1, s^2)=(0,0)$ or $(s_1, s_2)=(0,0)$, hence, by Prop.~\ref{all}a,  $\met(\Gamma', \Gamma')=2s_as^a=0$ so  (1) and (2) are both satisfied. If  neither $q'$ nor $p'$ vanish then  $q^*=[\q^*],$ $p^*=[\p^*]$,
where  
$$\begin{array}{l}
\q^*=\q\times \q'=E_3\times E_3'
=E_3\times E_a s^a
=s^1E^2-s^2 E^1,\\
\p^*=\p\times \p'=E^3\times (E^3)'
 =-E^3\times E^bs_b=
-s_1E_2+s_2E_1.
\end{array}
$$
The dancing condition is  then $\q^*\p^*=0$, i.e. 
 $(s^1E^2-s^2 E^1)(-s_1E_2+s_2E_1)=-s_as^a=0,$
 which is the  nullity condition on $\Gamma$.\end{proof}
 
 \begin{rmrk} It is clear  that both the dancing metric and the dancing condition are $\SLth$-invariant and homogeneous in the velocity   $\Gamma'$ of a parametrized curve $\Gamma$ in $\M$, thus defining each a field of tangent cones on $\M$. It is also clear from the formula of the isotropy representation (Eqn.~\eqref{isotropy}) that $\M$ admits a unique $\SLth$-invariant  conformal metric (of whatever signature). The main point of the last proposition, perhaps less obvious, is then  that the dancing condition  is {\em quadratic} in the velocities $\Gamma'$,  thus defining {\em some} conformal metric on $\M$. This point can be proved   in an elementary fashion, as we now proceed to show, and thus gives an alternative  proof of the last proposition. 
  \end{rmrk}
  
\subsection{A coordinate formula for the   conformal class of the dancing metric} \label{simple}
Let us use Cartesian coordinates $(x,y)$  for  a point   $q\in\RPt$ (in some affine chart) and the coordinates $(a,b)$ for a  line 
$y=ax+b$ (a point  $p\in\RPts$). If $q(t)$ is given by $(x(t), y(t))$ then its tangent line $y=Ax+B$ at time $t$ satisfies 
\be\label{d1}
y(t)=Ax(t)+B, \quad y'(t)=Ax'(t).
\ee 
Likewise, if $p(t)$ is  a curve in $\RPts$ given by $y=a(t)x+b(t)$ then its  ``turning point"    $(X,Y)$ at time $t$  satisfies
\be\label{d2}
Y=a(t)X+b(t), \quad 0=a'(t)X+b'(t).
\ee 

The dancing condition (``the turning point  lies on the tangent line") is then 
$$
Y=AX+B.
$$
Expressing $A,B, X,Y$ in the last equation in terms of $x,y,a,b$ and their derivatives via Eqns.~(\ref{d1})-(\ref{d2}), we obtain $a'[(y - b)x' - xy'] + b'[ax' - y']=0.$
Combining this calculation with Prop.~\ref{prop_danc},  we have shown
\begin{proposition}
The dancing metric $\met$ on  $M$ is given in the above local coordinates $x,y,a, b$ by 
$$\met\sim da[(y - b)dx - xdy] + db[adx - dy],$$
where $\sim$ denotes conformal equivalence (equality up to multiplication by some  non-vanishing function on $M$). 

\end{proposition}

\begin{rmrk} In fact, although somewhat less elementary, it is not hard   to show that the missing conformal factor on the right hand side of   the above formula  is  of the form $const./(y-ax-b)^2.$ 
 \end{rmrk}

\subsection{A cross-ratio  formula for the dancing metric}\label{ss_cr}
\newcommand{\tq}{\bar q}

\begin{definition} The cross-ratio of 4 distinct points $a_1,a_2,a_3,a_4$ on a line $\ell\subset \RPt$  is  $$ [a_1,a_2,a_3,a_4]:= \frac {x_1-x_3 }{x_1-x_4 }\cdot\frac{x_4-x_2} {x_3-x_2},$$ where $x_i$ is the coordinate of $a_i$ with respect   to some affine coordinate  $x$ on  $\ell$.\end{definition}

It is well-known (and not hard  to verify) that this definition is independent of the affine coordinate chosen on $\ell$  and that it is $\SLth$-invariant. 

\sn 

Now consider a non-degenerate parametrized curve $\Gamma$  in $\M$ and two points on it, $\Gamma(t)=(q(t), p(t))$ and $ \Gamma(t+\epsilon)=(q(t+\epsilon), p(t+\epsilon))$. These determine 4 colinear points 
$q,  q_\epsilon,  \tq,  \tq_\epsilon,$
 where  $q:=q(t),$ $q_\epsilon=q(t+\epsilon),$ and $ \tq,   \tq_\epsilon$   are the intersection points of the two lines $p:=p(t), p_\epsilon:=p(t+\epsilon)$ with the line $\ell$ through $q, q_\epsilon$ (resp.), as  in the picture. (The line $\ell$ is well-defined, for small enough $\epsilon$, by the non-degeneracy assumption on $\Gamma$). 
\begin{figure}[h]\centering
\includegraphics[width=0.35\textwidth]{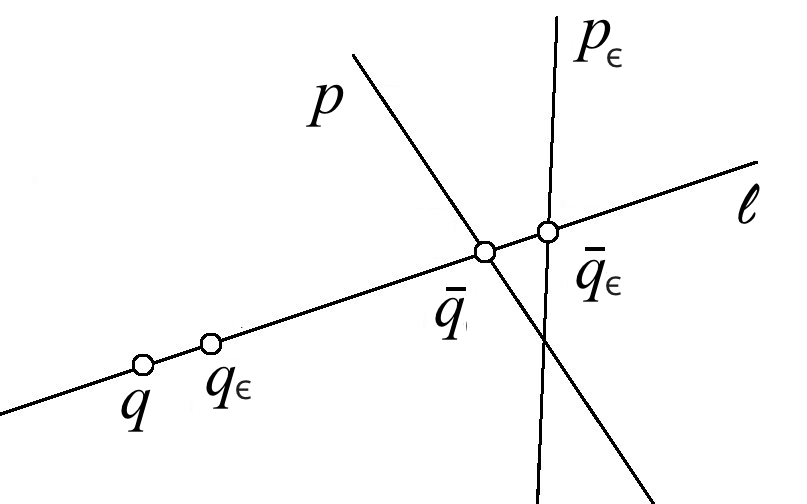}
\caption{\small The cross-ratio definition  of the dancing metric }
\end{figure}

\sn 

Let us expand the cross-ratio  of $q, q_\epsilon, \tq, \tq_\epsilon$ in powers of $\epsilon$. 

\begin{proposition}\label{cr} Let $\Gamma(t)$ be a non-degenerate parametrized curve in $\M$, $q_0, q_\epsilon, \tq_0, \tq_\epsilon$ as defined above and  $v=\Gamma'(t)$. Then 
$$[q, q_\epsilon, \tq, \tq_\epsilon]={1\over 2}\epsilon^2\met (v, v)+O(\epsilon^3),$$
where $\met$ is the dancing metric on $M$, as defined in Prop.~\ref{dancing}. 
\end{proposition}


\begin{proof}  Lift $\Gamma(t)$ to a curve $\tG(t)=(\q(t), \p(t))$ in $Q$. Then $\ell=[\q\times\q_\epsilon]$, 
$\tq=[\mathbf\tq]$ and  $\tq_\epsilon=[\mathbf\tq_\epsilon],$ where 
\begin{align*}
&\mathbf\tq=( \q\times\q_\epsilon)\times \p=\q_\epsilon-(\p\q_\epsilon)\q,\\ 
&\mathbf\tq_\epsilon=(\q\times\q_\epsilon)\times \p_\epsilon=(\p_\epsilon\q)\q_\epsilon-\q.
\end{align*}

Now  it easy to show that if 4 colinear points $a_1,\ldots, a_4\in\RP^2$ are given by  homogeneous coordinates  $\a_i\in\R^3\setminus 0$, such that  $\a_3=\a_1+\a_2,$ $\a_4=k\a_1+\a_2,$ then 
  $ [a_1,a_2,a_3,a_4]=k$ (see for example \cite{K}). 
 Using this  formula and the above expressions for $\mathbf\tq, \mathbf\tq_\epsilon$, we obtain, after some manipulations, 
$$[q,q_\epsilon, \tq, \tq_\epsilon]=1-{1\over(\p\q_\epsilon)(\p_\epsilon\q)}=
-\epsilon^2(\q\times\q')(\p\times\p')+O(\epsilon^3).$$
Now we use  the expression for $\met$ of Prop.~\ref{homo}.
\end{proof}

\subsection{Dancing pairs and their projective structure}\label{pj}
\begin{definition}A  {\em dancing pair} is a pair of parametrized curves $q(t), p(t)$ in  $ \RP^2, \RPts$ (resp.)  
obtained from the projections $q(t)=[\q(t)], p(t)=[\p(t)]$ of an integral curve $(\q(t), \p(t))$ of $(\Q,\D)$. If $q(t), p(t)$ is a dancing pair we say that $p(t)$ is a {\em dancing mate} of $q(t)$. 
\end{definition}
Equivalently, by Thm.~\ref{eleven}, this is the pair of curves one obtains from a non-degenerate  null curve in $\M$ with parallel self-dual tangent plane, when projecting it to $\RPt$ and $\RPts$. 

We already know that dancing pairs  satisfy the dancing condition.  
 We now want to study further the projective geometry of 
such pairs of curves, using the classical notions of projective differential geometry, such as the projective 
structure of a plane curve, projective curvature and projective arc length. We will  derive a 4th 
order ODE whose solutions give the {\em dancing mates} $p(t)$ of a given 
(\lc) curve $q(t)$.  We will give several  examples of dancing pairs, including the surprisingly non-trivial  case of the  dancing mates  
associated with a point moving on a circle.

 Differential projective geometry  is not so well-known nowadays, so  we  begin with a brief review of the pertinent notions. Our favorite   references are \'E.~Cartan's book \cite{Cbook} and the more modern references of Ovsienko-Tabachnikov \cite{OT} and Konovenko-Lychagin \cite{KL}. 

A {\em projective structure} on a curve $\gamma$ (a 1-dimensional manifold) is an atlas of charts $(U_\alpha, f_\alpha)$, where $\{U_\alpha\}$ is an open cover of $\gamma$ and the $f_\alpha:U_\alpha\to \RP^1$, called projective coordinates, are embeddings whose transitions functions $f_\alpha\circ f_\beta^{-1}$ are given by (restrictions of) M\"obius transformation in $\rm PGL_2(\R)$. 

For example, stereographic projection from any point $q$ on a conic $\CC\subset \RPt$ to some line $\ell$ (non-incident to $q$) gives $\CC$ a projective structure, independent of the point $q$ and line $\ell$ chosen (a theorem attributed to  Steiner, see  \cite[p.~7]{OT}).

 An embedded  curve $\gamma\subset\RPt$  is {\em\lc} if it has no inflection points (points where the tangent line has a 2nd order contact with the curve). Every \lc\  curve $\gamma\subset\RPt$ inherits a canonical projective structure. There are various equivalent ways to define this projective structure, but we will give the most classical one, using the {\em tautological ODE} associated with a plane curve (we follow here closely  Cartan's book \cite{Cbook}). 

Let $q(t)$ be a regular parametrization of $\gamma$, i.e. $q'\neq 0$, and $\q(t)$ a lift of $q(t)$ to $\R^3\setminus \{0\},$ i.e. $q(t)=[\q(t)].$ Then local-convexity (absence of inflection points) is equivalent to  $\det(\q(t),\q'(t),\q''(t))\neq 0$, so there are unique   $a_0,  a_1, a_2$ 
 (functions of $t$) such that  
\be\label{taut}
 \q'''+a_2\q''+a_1\q'+a_0\q=0.
 \ee
The last equation is  called the {\em tautological ODE} associated with $\gamma$ (or rather  its parametrized lift $\q(t)$). Solving   for  the unknowns $a_0,  a_1, a_2$ (by Kramer's rule), we get  
$$
a_0=-{J\over I},\quad
a_1={K\over I},\quad
a_2=-{I'\over I},
$$
where
$$I=\det(\q,\q', \q''), \quad J=\det(\q',\q'',\q'''), \quad 
K=\det(\q, \q'', \q''') .
$$

The tautological ODE Eqn.~\eqref{taut} depends on the  choice of parametrized lift  $\q(t)$. One can modify  $\q(t) $ in two ways:
\begin{itemize}

\item {\em Re-scaling}:  $\q(t) \mapsto \bar\q(t)=\lambda(t) \q(t) ,$ where $ \lambda(t) \in \R^*.$ This changes 
$I\mapsto\lambda^3I,$ so if $I\neq 0$ (no inflection points) one  can rescale (uniquely)  to say $I=1$, 
then obtain $a_2= 0$.  
\item {\em Re-parametrization}: $ t\mapsto \bar  t=f(t),$ for some diffeomorphism  $f$. This changes $I\mapsto (f')^3I,$ so again, if $I\neq 0$ then one can reparametrize (uniquely up to an additive constant) to $I=1$, so as to obtain $a_2=0$. 

 \end{itemize}
 So one can achieve a tautological ODE for $\gamma$ with $a_2=0$ by either rescaling or reparametrization. 
Can we combine reparametrization and rescaling so as to reduce the tautological ODE to  $\q'''+a_0\q=0$?

The answer is ``yes" and the resulting OED is called  the  {\em Laguerre-Forsyth  form}   (LF) of the tautological 
ODE for $\gamma$. A straightforward calculation (\cite{Cbook}, p.~48) 
shows 

\begin{proposition}\label{TwentyOne} Given a \lc\ curve $\gamma\subset\RPt$ with  a parametrized lift $\q(t)$ satisfying $\q'''+a_1\q'+a_0\q=0$,

\sn\begin{enumerate}[leftmargin=30pt,label=(\arabic*)]\setlength\itemsep{5pt}

\item one  can achieve the LF form by modifying  $\q(t)$ to  $\bar \q(\bar  t)= f'(t)\q(t),$ where $\bar  t=f(t)$ solves
 $$S(f)={a_1\over 4},$$
and  where $$S(f)={1\over 2}{f'''\over f'}-{3\over 4}\left({f''\over f'}\right)^2$$
 is the {\em Schwarzian derivative} of $f$. 
 
\item The LF form is unique  up to the  change $\q(t)\mapsto \bar \q(\bar  t)= f'(t)\q(t),$ where $\bar  t=f(t)$ is a M\"obius transformation. 

\item Given an LF form $\q'''+a_0\q=0$  for $\gamma$, the  one form $d\sigma=(a_0)^{1/3}dt$ is a well-defined 1-form on $\gamma$ (independent of the particular LF form chosen), called the {\em projective arc length} \cite[p.~50]{Cbook}.

\end{enumerate} 

\end{proposition}
It follows from  item (2) that the LF form defines a local coordinate $t$ on $\gamma$, well-defined up to a M\"obius transformation, hence a {\em projective structure}  on  $\gamma$.

\begin{rmrk} It is possible to extend the definition of the projective structure  to all curves, not necessarily \lc\ (see  \cite{DZ}). 
\end{rmrk}

\begin{example} A classical application of the last proposition is to show that 
{\em  the duality map $\gamma\to \gamma^*$ preserves the projective structure but reverses the projective arc length} 
(provided both $\gamma$ and $\gamma^*$ are \lc): 
parametrize $\gamma$ by $[\q(t)]$ in the LF form, i.e. $\q'''+a_0\q=0$, 
then $\gamma^*$ is parametrized by $[\p(t)]$, where $\p(t)=\q(t)\times \q'(t)$. 
Then one can calculate easily that $\p(t)$ satisfies $\p'''-a_0\p=0,$ which is also in the LF form, hence 
$t$ is a common projective parameter on $\gamma,\gamma^*$, so 
$[\q(t)]\mapsto [\p(t)]$ preserves the projective structure, but reverses the projective arc length.
\end{example}

\begin{proposition}\label{TwentyOne}
Let $\gamma, \bar \gamma$ be a pair of non-degenerate curves in $\RPt,\RPts$ (resp.), parametrized by a dancing pair $q(t), p(t)$ (i.e.~$q(t)=[\q(t)], p(t)=[\p(t)]$, where $\p'=\q\times \q'$). Then the map $\gamma\to\bar \gamma$, $q(t)\mapsto p(t)$, is projective, i.e. preserves the natural projective structures on $\gamma, \bar \gamma$ induced by their embedding in $\RPt, \RPts$ (resp.). 
\end{proposition}

\begin{proof} According to the last proposition, it is enough to show that $(\q(t), \p(t))$ can be reparametrized in such a way that    $\q(t), \p(t)$ satisfy each a tautological ODE with $a_2=0$ and the same $a_1$. 

 \begin{lemma}\label{mucho} Let $(\q(t),\p(t))\in\Rtt$ be a solution of $\p'=\q\times\q',$ $   \p\q=1,$ with $I(\q)=\det(\q,\q', \q'')\neq 0$ and $ I(\p)=\det(\p,\p', \p'')\neq 0$. Let $I=I(\q), \bar  I=I(\p), J=J(\q), \bar J=J(\p),$ etc. Then 
\begin{enumerate}
\item $\p'\q=\p\q'=\p'\q'=\p\q''=\p''\q=0.$
\item $I\p=\q'\times\q'',$ $\bar I\q=\p'\times\p''.$
\item $\q'=-\p\times\p'.$

\item $I^2+J=\bar I^2- \bar J=0$. 
\item $\bar I=I,$ $\bar J=-J,$ $\bar  K=K.$
\item  $\bar a_2=a_2,$  $\bar a_1=a_1, $  $\bar a_0=-a_0$. 
\end{enumerate}
\end{lemma}

\begin{proof}

\begin{enumerate}[leftmargin=18pt,label=(\arabic*)]\setlength\itemsep{5pt}
\item From $\p'=\q\times \q'\ent \p'\q=\p'\q'=0.$ From $\p\q=1\ent \p'\q+\p\q'=0\ent \p\q'=0\ent 0=(\p\q')'=\p'\q'+ \p\q''=\p\q''.$ Similarly, $0=(\p'\q)'=\p''\q+\p'\q'=\p''\q.$

\item From (1),  $\p\q'=\p\q''=0\ent c\p=\q'\times\q''$ for some function $c$ (we assume $I\neq 0$, hence $\q'\times\q''\neq 0$). Taking dot product of last equation with $\q$ and using $\p\q=1$ we get $c=I\ent I\p=\q'\times\q''.$

 Next, from (1), $ \p'\q=\p''\q=0\ent \bar c\q=\p'\times\p''$
  for some function $\bar c$ (here we assume 
  $\bar I\neq 0$). Take dot product with $\p$ and get $\bar c=\bar I \ent \bar I\q=\p'\times \p''.$
 
 \item (This was already shown in Lemma~\ref{dual} but we give another proof here). From (1), $\p\q'=\p'\q'=0\ent \q'=f\p\times\p'$ for some function $f$. Cross product with $\q$, use the vector identity $$(\p_1\times\p_2)\times \q =(\p_1\q)\p_2 - (\p_2\q)\p_1,$$  and get $-\p'=\q'\times \q=f(\p\times\p')\times \q=f[(\p\q)\p'-(\p'\q)\p]=f\p'\ent f=-1\ent \q'=-\p\times\p'.$
   
 \item  $I\p=\q'\times\q'', \p'=\q\times \q'\ent I'\p+I(\q\times\q')=\q'\times\q'''.$ Now dot product with $\q''$, use $\p\q''=0$ and get $I^2+J=0.$
   Very similarly, get
   $(\bar I)^2-\bar J=0.$
   
   \item Use the vector identity
   $$\det(\q_1\times\q_2, \q_2\times \q_3,\q_3\times \q_1)=[\det(\q_1,\q_2, \q_3)]^2,$$ to get 
   $I\bar I=\det(I\p, \p',\p'')=\det(\q'\times\q'', \q\times\q', \q\times\q'')=I^2,$ hence
   $I=\bar I.$ 
   
   \mn From (4), $\bar J=\bar I^2=I^2=-J.$
   
   \mn From $\p'=\q\times \q\ent \p''=\q\times \q''\ent \p''\q''=0\ent \p'''\q''+\p''\q'''=0.$
   Now $K=\det(\q, \q'', \q''')=(\q\times \q'')\q'''=\p''\q''',$ $\bar K=\det(\p,\p'', \p''')=-\p'''\q'',$ hence 
   $K-\bar K=\p''\q'''+\p'''\q''=(\p''\q'')'=0.$
   
   \item Immediate from item  (5) and the definition of $a_0, a_1, a_2$. \qedhere 
  \end{enumerate}
\end{proof}

Now $\gamma$ is \lc\ so we can reparametrize $\q(t)$ to achieve $I(\q)=1.$ The equation $\p'=\q\times\q'$ is reparametrization invariant so it still holds. It follows from item (5) of the lemma that $I(\p)=1$ as well, hence both $a_2=\bar a_2=0.$ From item (6) of the lemma we have that $a_1=\bar a_1$. Hence the  equation for projective parameter    $S(f)=a_1/4$ is the same equation for both curves $\q(t)$ and $\p(t)$. \end{proof}

 \subsection{An aside: space curves with constant ``centro-affine torsion"}\label{sec:centro}
We mention here in passing a curious geometric interpretation of a formula that appeared  
during the proof of Prop.~\ref{TwentyOne} (see part (4) of Lemma~\ref{mucho}):
\begin{equation}\label{peqn}
 J(\q)+I^2(\q)=0,
 \end{equation} 
where $$I(\q)=\det(\q,\q',\q''),\quad J(\q)=\det(\q,\q'', \q''').$$

\newcommand{\peqn}{(\ref{peqn})}

Effectively, this formula means that it is possible to eliminate the $\p$ variable from our system of Eqns.~   \eqref{eqnsss}, reducing them  to  a {\em single}  3rd order ODE   for a  space curve $\q(t)$.

In fact, it is not hard to show that Eqn.~\peqn\ is {\em equivalent} to Eqns.~\eqnsss; given a nondegenerate ($I(\q)\neq 0$)
solution  $\q(t)$ to Eqn.~\peqn, use the ``moving  frame"  $\q(t), \q'(t), \q''(t)$ to define $\p(t)$ by  
\begin{equation}\label{peq}
\p(t)\q(t)=1, \;\p(t)\q'(t)=0, \;\p(t)\q''(t)=0,
\end{equation}
then check that Eqn.~\peqn\ implies  that 
 $(\q(t), \p(t))$ is a solution to Eqns.~\eqnsss. 
 
The curve $\p(t)$ associated to a non-degenerate curve $\q(t)$ via Eqns.~(\ref{peq}) represents the   
  {\em osculating plane} $H_t$ along   $\q(t)$,  via the formula $H_t=\{\q|\p(t)\q=1\}.$
  
For any space curve (with $I\neq 0$) the quantity $\mathcal J=J/I^2$ is parametrization-independent  and $\SLth$-invariant, called by some authors  the {\em (unimodular) centro-affine torsion} \cite{O}. Hence Eqns.~\eqnsss\  can be also interpreted as the equations  for space curves with $\mathcal J=-1$.

\subsection{Projective involutes and the  dancing mate equation}\label{sec:mate} The reader may suspect now  that the necessary condition of Prop.~\ref{TwentyOne} is also sufficient for a pair of curves to be a dancing pair. This is not so, as the following example shows.

\begin{example} Let $\gamma, \bar \gamma$ be the pair consisting of a  circle $q(t)=[\cos(t), \sin(t), 1]$ 
 and the dual of the concentric circle 
$p^*(t)=[\sqrt{2}\cos(t+\pi/4), \sqrt{2}\sin(t+\pi/4), 1].$ 
One can check easily that $(q(t), p(t))$ satisfies the dancing condition (i.e. defines  a null curve in $\M$) and that the map $q(t)\mapsto p(t)$ is projective (as the restriction to $\gamma$ of an element in $\SLth$: a  dilation followed by a rotation). Nevertheless,  the pair of curves  $q(t), p(t)$ is not a dancing pair  
(there is no way to lift $(q(t), p(t))$ to a solution $(\q(t), \p(t))$ of $\p\q=1, \p'=\q\times \q'$). 
\end{example}

\begin{figure}[h!]\centering
\includegraphics[width=0.4\textwidth]{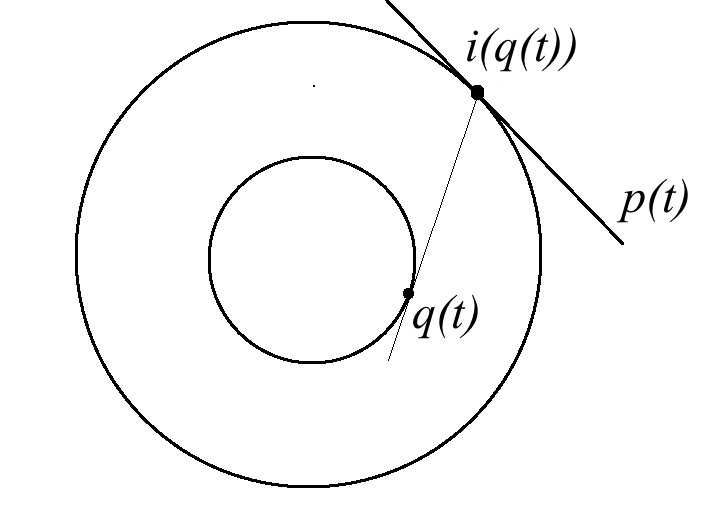}
\caption{A  projective involute which is not a dancing pair}
\end{figure}

 We are going to study carefully the situation now and find an extra condition that the map $q(t)\mapsto p(t)$ should satisfy, for the pair of curves $q(t), p(t)$ to be  a dancing pair.

\begin{definition}  Let $\gamma\subset\RPt$ be a \lc\ curve.   A {\em projective  involute}  of $\gamma$ is a smooth map  $i:\gamma\to \RP^2$ such that  

\begin{itemize}\item  for all $q\in \gamma$,   $i(q)\in q^*$ (the tangent line to $\gamma$ at $q$).

\item  $i$ is a projective immersion. 
\end{itemize}
\end{definition}

The last phrase means that $i$ is an immersion and  the resulting  local diffeomorphism between $\gamma$ and its image  is projective with respect to the natural projective structures on $\gamma$ and  $i(\gamma)$, defined by their embedding in $\RPt$.

 \begin{proposition} Near a non-inflection point  of a curve $\gamma\subset \RPt$ there is a 4-parameter family of
   projective involutes, given by the solutions of the following 4th order ODE:
 if $\gamma$ is given by a tautological ODE in the  LF form $A'''+rA=0,$ then its projective involutes  are  given  by $[A(t) ]\mapsto [B(t) ], $ where  $B(t) =(C-y'(t) )A(t) +y(t) A'(t) $, 
$C$ is an  arbitrary  constant and  $y(t) $ is a solution of  the  ODE 
      $$y^{(4)}+2{y'''(y'-C)\over y}+3ry'+r'y=0.$$
          \end{proposition}

 \begin{proof}  Calculate, using $A'''+rA=0$: 
 \begin{eqnarray*} 
B&=&xA+ yA'\\
B'&=&x'A+ (x+y')A' + yA''\\
B''&=&(x''-ry)A+ (2x'+y'')A'+ (x+2y')A''\\
B'''&=&(x'''-r'y-r(x+3y'))A+ (3x''-ry+y''')A'+ 3(x'+y'')A''
 \end{eqnarray*}
Hence $B\times B'''=0\ent y(x'+y'')=x(x'+y'')=x(3x''+y''')-y(x'''-3ry'-r'y)=0,$ then $y\neq 0\ent x'+y''=0\ent x+y'=C, $ for some constant $C$, hence  
 $$y^{(4)}= 2{y'''(C-y')\over y}-3ry'-r'y.$$
This gives a 5 parameter family of solutions. Then imposing say $\det(B,B',B'')=1$ reduces it to a 4 parameter family (removing the scaling ambiguity on $B$). \end{proof}

 Now given a nondegenerate $\gamma\subset \RPt$, parametrized by $q(t)$, we 
know, by Prop.~\ref{TwentyOne} and the preceding example, that each of its  dancing mates $p(t)$ gives rise to a projective involute
$q(t)\mapsto p(t)\mapsto p^*(t)$. The dancing mates of $\gamma$ form a  3 parameter sub-family of 
the projective involutes,  as they are obtained by lifting 
$\gamma$ horizontally via $\Q\to \RP^2$, followed by the projection $\Q\to\RPts$. 
We are thus looking for a single equation characterizing projective involutes of $\gamma$ that correspond to dancing mates. 
\begin{proposition}\label{TwentyThree} 
Let  $\gamma\subset \RPt$ be a nondegenerate curve with a tautological ODE in LF form $A'''+rA=0$. Let $i:\gamma\to\RPt$ be a projective involute given in homogeneous coordinates by $B=xA+yA'.$ Then $b=B\times B'$ is a dancing mate of $A$  if and only if $x+y'=0$. That is, $C=0$ in the previous proposition so  $y$ satisfies the ODE
$$y^{(4)}+2{y'''y'\over y}+3ry'+r'y=0.$$ 
\end{proposition}

\begin{proof} Let $(\q(t), \p(t))$ be an integral curve of $(\Q,\D)$. Then $\p^*=\p\times \p'=-\q'$. Then, to bring both $\q,\p^*$ to LF form we
  need the same projective parameter $\bar t=f(t)$,  so that $A(\bar t)=f'(t)\q(t) ,$ $B(\bar t)=f'(t)\p^*(t)=-f'(t)\q'(t).$ 
  Taking derivative of $A(\bar t)=f'(t)\q(t)$  with respect $t$, get $f'(dA/d\bar t)=f''\q+f'\q'=(f''/f')A-B$, 
  hence $B=xA+y(dA/d\bar t),$ with $x=f''/f',$ $y=-f'\ent x+dy/d\bar t=f''/f'-f''/f'=0.$\end{proof}
  
\begin{rmrk}  The geometric meaning of the condition $x+y'=0$ is the following.  Since $B=xA+yA'$ then $B'=x'A+(x+y')A'+yA''$. Hence the condition $x+y'=0$ means that $B'$ is the intersection point of the line $b$ and the line $a'=A\times A''$ (the line connecting $A$ and $A''$). In Sect.~\ref{RT} below, we will further interpret this condition in terms of the osculating conic and Cartan's developments. 
 \end{rmrk}

\subsection{Example: dancing around a circle}\label{sec:circle}Take $\gamma$ to be a conic, e.g. a circle, $A=(1+t^2,2t,1-t^2)$. Then $A'''=0$, so $A(t) $ is in the LF form with 
 $r=0$. Then the dancing mate equation in this case is $$y^{(4)}+2{y'''y'\over y}=0.$$
  Any quadratic polynomial solves this (since  $y'''=0$),  and the corresponding involute $B=-y'A+yA$ is a straight line. To show this, take $y=at^2+bt+c,$ then 
   $$B=( bt^2+2(c-a)t-b,-2at^2+2c, -bt^2-2(a+c)t-b),$$ which  is contained in the 2-plane 
   $$(a+c)x+by+(c-a)z=0,$$ so projects into a straight line in $\RPt$.

   If $y$ is not a quadratic polynomial, then  in a neighborhood of $t$ where $y'''(t) \neq 0$,
  $$0={y^{(4)}\over y'''}+2{y'\over y}=[\log (y''' y^2)]'=0\ent y''' y^2=const.$$
  
 Now we can assume, without loss of generality, that $const.=1$ (multiplying $y$ by a constant does not affect $[B(t)]$, so we end up with the ODE
$$y'''y^2=1.$$ 

We do not know how to solve this equation explicitly, so we do it numerically. The result is Fig. 2 of Sect. \ref{ProjGeom} above.  

 A few words about this drawing: we make the drawing in the $XY$ plane, where the circle is $X^2+Y^2=1$ and dancing curves around it  are obtained from solutions of $y'''y^2=1 $ via the formulas
\begin{eqnarray*}
B&=&-y'A+yA'-y'\left(1+t^2-2tz,2(t-z), 1-t^2-2tz\right),\\
(X,Y)&=&\left({B_2\over B_1}, {B_3\over B_1}\right)={(2(t-z), 1-t^2-2tz)\over 1-t^2-2tz}, \quad z=y/y'.
\end{eqnarray*}

The projective coordinate $t$ on a conic  $\CC$ misses a point (the point at ``infinity"), so when integrating this equation  numerically, one needs a second coordinate, $\bar  t=f(t) =1/t,$ and the transformation formulas $\bar  y=f'y$, etc.

\subsection{Example: dancing  pairs of constant projective curvature}\label{const}
 
\n{\bf The idea:}  fix a point 
$(\q_0,\p_0)\in Q$ and an element $Y\in\mathfrak{sl}(3,\R)$  such that $Y\cdot(\q_0,\p_0)$ 
is $\D$-horizontal. That is, $\p_0d\q$ and $d\p-\q_0\times d\q$ both vanish on $Y\cdot (\q_0, \p_0)=(Y\q_0, -\p_0Y).$  The subspace of such $Y$ has codimension 3 in $\mathfrak{sl}(3,\R)$, i.e.~is 5-dimensional, since $\SL(3,\R)$ acts transitively on $Q$ and $\D$ has corank 3.

Then the orbit of $(\q_0, \p_0)$ under the flow of $Y$, 
 $$(\q(t), \p(t))=exp(tY)\cdot(\q_0, \p_0)=(exp(tY)\q_0, \p_0exp(-tY))$$ 
 is an integral curve of $\D$
 (this follows from the $\SL(3,\R)$-invariance of $\D$). The projected curves $q(t)=[\q(t)], p(t)=[\p(t)]$ are then a dancing pair. Each of the curves is an orbit of the 1-parameter subgroup $exp(tY)$ of $\SLth$. Such curves are called $W$-curves or ``pathcurves".  They are very interesting curves, studied by Klein and Lie in 1871 \cite{KLie}. They are: straight lines and conics,   exponential curves, logarithmic spirals  and ``generalized parabolas" (see below for explicit formulas).  This class of curves (except lines and conics, considered degenerate) coincides with the class of curves with {\em constant projective curvature} $\kappa$. 
 
\sn {\bf About the projective curvature.} The projective curvature $\kappa$ of a curve $\gamma\subset \RPt$ is a function defined along $\gamma$, away from {\em sextactic}  points,  where the osculating conic has order of contact with $\gamma$  higher then expected (5th or higher).  
The  sextactic points are also given by the vanishing of the {\em projective arc length}  element $d\sigma$, so away from such points one can use $\sigma$ as a natural parameter on $\gamma$ and compare it to a projective parameter $t$, given by the projective structure (see Prop.~\ref{TwentyOne}).  More precisely, when $d\sigma\neq 0$, it defines a local diffeomorphism $\RP^1\to \gamma$, whose  Schwarzian derivative is  the  quadratic form $\kappa(\sigma)(d\sigma)^2$. The pair $\{ d\sigma, \kappa\}$  forms a complete set of projective invariants for curves in $\RPt$ (analogues to the arc length element and curvature for regular curves in the euclidean plane). For  curves with constant projective curvature, the constant $\kappa$ by  itself is a complete invariant (also same as in the euclidean case). Along a conic $d\sigma\equiv 0$ and so $\kappa$ is not defined.

 \begin{rmrk}  In the book of Ovsienko-Tabachnikov \cite{OT} (a beautiful book, we highly recommend it) the  term ``projective curvature" is used to denote what we call here the projective structure and it is stated  that ``the projective curvature is, by no means, a function on the curve" (\cite[p.~14]{OT},  the online version). This can be somewhat   confusing if one does not realize the difference in usage of terminology. We adhere to the classical terminology, as in  Cartan's book \cite{Cbook}. 
 \end{rmrk}

The classification of curves with constant projective curvature $\kappa$, up to projective equivalence,  is as follows. There are two generic cases, divided (strangely enough) by the borderline value $\kappa_0=-3(32)^{-1/3}\approx-0.94$: 
\begin{itemize}
\item $\kappa>\kappa_0$: logarithmic spirals, $r=e^{a\theta}$, $a>0$ (in polar coordinates); 
\item $\kappa=\kappa_0$: the exponential curve $y=e^x$;
\item $\kappa<\kappa_0$: generalized parabolas, $y=x^m$, $m>0$, $m\neq 2, 1, 1/2.$ 
\end{itemize}

\mn{\bf The curves.} Take $\q_0=(0, 0, 1)^t$, $\p_0=(1,0,0).$ Then $Y\cdot(\q_0,\p_0)\in\D_{(\q_0,\p_0)}$, $Y\in\slth$,  implies 
$$Y=
\left(\begin{array}{c|c}A&\bv\\ \hline\bv^*& 0\end{array}\right),\quad 
\bv={v_1\choose v_2}, \quad \bv^*=(v_2, -v_1), \quad A\in\slt.$$
Let $H_0\cong \SLt$ be the stabilizer subgroup of $(\q_0, \p_0)$. It  acts on $Y$ by the adjoint representation, $$(A,\bv)\mapsto (hAh^{-1}, h\bv), \quad h\in\SLt.$$ Then, reducing by this $\SLt$-action as well as by rescaling, $Y\mapsto \lambda Y, $ $ \lambda\in\R^*$ (this just reparametrizes the orbit), and  removing orbits which are fixed points and straight lines, we are left with a list of ``normal forms''  of $Y$ (two one-parameter families and one isolated case):
\begin{subequations}\label{Y}

\begin{align}  Y_1&:=
 \left( 
\begin{array}{rrr} 
1& 0 &1\\ 
0&-1&a\\ 
a&-1&0
\end {array} 
\right), 
\quad a> 0.\\  
Y_2&:=
\left( \begin{array}{rrr} 
0&1&b\\
-1&0&0\\
0&-b&0
\end{array} 
\right), 
\quad b>0,\\ 
Y_3&:=
\left(
\begin{array}{rrr}
0 &1&0\\ 
0&0&1\\ 
1&0&0
\end{array} 
\right).
\end{align}
\end{subequations}

\begin{proposition} The pair of curves $[\q(t)], [\p(t)],$ in $\RPt, \RPts$ (resp.), where $\q(t)=exp(tY)\q_0$, $\p(t)=\p_0exp(-tY)$,
  $\q_0=(0,0,1)^t$, $\p_0=(0,0,1)$ and $Y$ is any of the matrices in Eqns.~\eqref{Y} above, is a dancing pair  of curves with constant projective curvature $\kappa$ (same value of $\kappa$ for each member of the pair). All values of $\kappa\in\R$ can be obtained in such a way.

\end{proposition}

\begin{proof} A matrix $Y$ with $\tr(Y)=0$  has characteristic  polynomial  of the form $\det(\lambda I-Y)={\lambda}^{3}+a_1\lambda+a_0.$ 
Then, using $Y^3+a_1Y+a_0I=0$ (Cayley-Hamilton), we have that  $\q(t):=exp(tY)\q_0$ satisfies the tautological ODE
$$\q'''+a_1\q'+a_0\q=0.$$ 
From Cartan's formulas (\cite{Cbook}, p.~69 and  p.~71), we then find easily
$$\kappa=a_1a_0^{-2/3}/2.$$
Now in our case, the characteristic  polynomials are 
\begin{align*}
 (a)\quad&{\lambda}^{3}-\lambda-2a, &
 (b)\quad  & {\lambda}^{3}+\lambda-\,b^2, &
  (c)\quad & {\lambda}^{3}-1,
  \end{align*}   
hence we get projective curvatures 
\begin{align*}
 (a)\quad&  \kappa=-(32a^2)^{-1/3}&
 (b)\quad  &\kappa=b^{-4/3}/2,&
  (c)\quad & \kappa=0.
  \end{align*}   
We thus get all possible values of $\kappa$. \end{proof}
To visualize such a pair, we draw the pair $(q(t), p^*(t))$, where $p^*(t)=[\p^*(t)]$ is the curve dual  to $p(t)$, given by $\p^*(t)=\p(t)\times\p'(t)=-Y\q(t)$.

\begin{figure}[h]\centering
\includegraphics[width=0.6\textwidth]{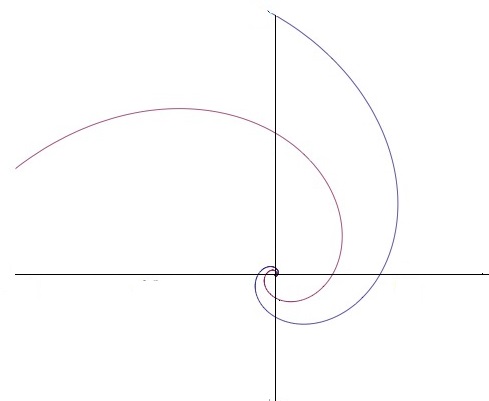}
\caption{\small A dancing pair of logarithmic spirals with $\kappa=0$.}
\end{figure}

\subsection{Projective rolling without slipping and twisting}\label{RT}

\newcommand{\RC}{\mathcal{RC}}

\subsubsection{About riemannian rolling}
Let us describe  first   ordinary (riemannian) rolling, following \cite[p.~456]{BrHs}.   Let $(\Sigma_i, \met_i)$, $i=1,2,$ be two   riemannian surfaces. 
The configuration space for  the rolling of the two  surfaces along each other  is the space $\RC$ of {\em riemannian contact elements} $(u_1, u_2, \psi)$, where $u_i\in\Sigma_i$ and $$\psi:T_{u_1}\Sigma_1\to T_{u_2}\Sigma_2$$ is an isometry.  
$\RC$ is a 5-manifold, and if $\Sigma_i$ are oriented then $\RC$ is the disjoint union $\RC=\RC^+\sqcup\RC^-$, where each $\RC^\pm$ is a circle bundle over $\Sigma_1\times\Sigma_2$ in an obvious way, so that  $\RC^+$ consists of the  orientation preserving riemannian contact elements and   $\RC^-$ are the orientation  reversing.

A parametrized curve $(u_1(t), u_2(t), \psi(t))$ in $\RC$ satisfies the {\em non-slip condition} if 
$$u'_2(t)=\psi(t)u_1'(t)$$ 
for all $t$. It satisfies also the {\em no-twist condition}   if for every  parallel vector field $v_1(t)$   along $u_1(t)$, 
$$v_2(t)=\psi(t)v_1(t)$$
 is parallel along $u_2(t)$ (``parallel" is with respect to   the Levi-Civita connection of the corresponding metric). 

It is easy to show that these two conditions define together a rank 2 distribution $\D\subset T\RC$ which is $(2,3,5)$ unless the surfaces are isometric (\cite{BrHs}, p. 458). For some special pairs of surfaces (e.g. a pair of round spheres  of radius ratio $3:1$) $\D$ is maximally symmetric, i.e. admits $\g_2$-symmetry (maximum possible). Recently \cite{AN}, some new pairs of surfaces were found where the corresponding {\em rolling distribution} $(\RC, \D)$ admits $\g_2$-symmetry, but the general case is not settled yet.  

Now  in \cite{AN} it was  noticed that  riemannian rolling can be reformulated as follows. 
Let $\M:=\Sigma_1\times\Sigma_2$, equipped with the difference metric $\met=\met_1\ominus \met_2$. 
This is a pseudo-riemannian metric of signature $(2,2)$. Then one can check easily that $\psi:T_{u_1}\Sigma_1\to T_{u_2}\Sigma_2$ is an isometry 
if and only if its graph 
$$W_\psi=\{(v,\psi v)|v \in T_{u_1}\Sigma_1\}\subset 
T_{u_1}\Sigma_1\oplus T_{u_2}\Sigma_2\simeq T_{(u_1,u_2)}\M$$ 
is a {\em non-principal null 2-plane}; i.e.~a null 2-plane not of the form  $T_{u_1}\Sigma_1\oplus \{0\}$ or  $\{0\}\oplus T_{u_2}\Sigma_2$ (compare with Cor.~\ref{grancor}a). This defines an embedding  of $\RC$ in the total space of the twistor fibration $\T \M\to\M$ of $(\M, [\met])$ (see Sect.~\ref{twist}). Furthermore, if $\Sigma_i$ are oriented then  one can orient $\M$ so that $\RC^+$ (orientation preserving $\psi$'s) is mapped to the self-dual twistor space  $\T^+\M$ and $\RC^-$ to the anti-self-dual twistor space $\T^-\M$. Finally, it is shown in \cite{AN}, that under the embedding $\RC\hookrightarrow\T\M$, the {\em rolling distribution} $\D$ on $\RC$ goes over to the {\em twistor distribution} associated with  Levi-Civita connection of $(\M, \met)$. 

 In what follows, we give a similar ``rolling interpretation" of the self-dual twistor space of  the dancing space $(\M, [\met])$, and thus, via the identification $\Q\hookrightarrow \T^+\M$ of Thm.~\ref{ident}, a ``rolling interpretation" of our Eqns.~\eqnsss. The novelty here is that the dancing metric $(\M, \met)$ is irreducible, i.e.   not a difference metric as in the case of riemannian rolling. And yet,  it can be given a rolling interpretation of some sort and in addition   admits $\g_2$-symmetry. 

We try to keep our   terminology as close as possible to the above terminology of riemannian rolling, in order to make the analogy transparent.

\subsubsection{A natural isomorphism of projective spaces} 

 A projective isomorphism of two projective spaces $\P(V),$ $\P(W)$  is the projectivization $[T]$ of a linear isomorphism $T:V\eto W$ of the underlying vector spaces, $[T]:[v]\mapsto [Tv]$. Two linear isomorphisms $T, T':V\to W$ induce the same projective isomorphism if and only if $T'=\lambda T$ for some $\lambda \in \R^*$.

For each non-incident pair $(q,p)\in\M$ we define a projective isomorphism  
\be\label{eq:pi}\Psi_{q,p}: \P(T_q\RPt)\eto\P(T_p\RPts) 
\ee
by first identifying  $\P(T_q\RPt)$ with the pencil of   lines through $q$ and $ \P(T_p\RPts)$ with the  points on the line  $p$. We then send a line $\ell$ through $q$ to its   intersection point $\ell^*$ with $p$.
One can verify easily that $\Psi_{q,p}$  is a projective isomorphism. 

\begin{figure}[h]\centering
\includegraphics[width=0.3\textwidth]{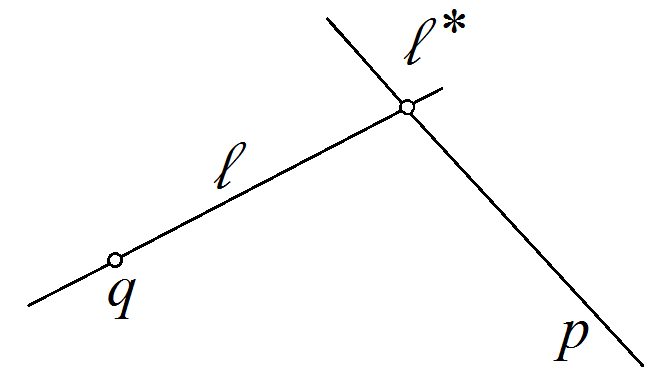}
\caption{\small The natural isomorphism $\Psi_{q,p}:\ell\mapsto \ell^*$}
\end{figure}

\subsubsection{Projective contact}

\begin{definition} 
 A {\em projective contact element} between $\RPt$ and  $\RPts$ is a triple  $(q,p,\psi)$ where $(q,p)\in\M$ and  
 $$\psi:T_q\RPt\to T_p\RPts$$
  is a  linear isomorphism covering the natural projective  isomorphism $\Psi_{q,p}$ of Eq.~\eqref{eq:pi}; that is, $[\psi]=\Psi_{q,p}.$ 
 The set of projective contact elements forms a principal $\R^*$-fibration  $\PC\to \M,$ $(q,p,\psi)\mapsto (q,p)$. 
\end{definition}

\begin{rmrk} We only allow projective contacts between $\RPt$ and $\RPts$  at a {\em non-incident} pair $(q,p)\in\M.$ 
 \end{rmrk}
 
Let us look at  the projective contact condition on $\psi$. We take a non-zero $v\in T_q\RPt$ and let $w=\psi(v)$.  To $v$ corresponds a line $\ell$ through $q$, tangent to $v$ at $q$. Likewise, to $w$ corresponds a point  $\ell^*\in p$, whose dual line in $\RPts$ (the pencil of lines through $\ell^*$) is tangent to $w$ at $p$. The projective contact condition on $\psi$ is then the incidence relation  $\ell^*\in \ell$. But this is precisely the {\em dancing condition}, i.e. $(v,w)\in T_{(q,p)}\M$ is a null vector. In other words, the graph of $\psi$,
$$W_\psi=\{(v, \psi (v))\st v\in T_q\RPt\}\subset T_q\RPt\oplus T_p\RPts=T_{(q,p)}\M,$$
is a null 2-plane. More precisely,

\begin{proposition}
A linear isomorphism $\psi:T_q\RPt\to T_p\RPts$ is a projective contact if and only if its graph $W_\psi\subset T_{(q,p)}\M$ is a non-principal self-dual null 2-plane (see Cor.~\ref{grancor}). 
\end{proposition}

\begin{proof} We recall from Sect.~\ref{proofs}: a local section $\sigma$ of   $j:\SLth\to M$ around $(q,p)\in M$ provides  a null coframing
$\weta:=\sigma^*\eta=(\weta^1, \weta^2, \weta_1, \weta_2)^t$ such  that $T_q\RPt\oplus\{0\}=\{ \weta_1= \weta_2=0\}$,  
$\{0\}\oplus T_p\RPts=\{ \weta^1= \weta^2=0\}$  and
$\met=2\weta_a\,\weta^a.$ 
Let  
$f=(f_1, f_2, f^1, f^2)$ be the dual framing, $\psi(f_a)=\psi_{ab}f^b.$ Now   the projective contact condition is $\psi(v)(v)=0\ent \psi_{ab}=-\psi_{ba}.$ Say $\psi(f_1)=\lambda f^2,$ $\psi(f_2)=-\lambda f^1$ for some $\lambda\in\R^*\ent  W_\psi=\Span\{f_1+\lambda f^2, f_2-\lambda f^1\}=\Ker\{\lambda\weta^1-\weta_2, \lambda\weta^2-\weta_1\}$. The 2-form corresponding to $W_\psi$  is thus $\beta=(\lambda\weta^1-\weta_2)\wedge(\lambda\weta^2+\weta_1).$ Using formula (\ref{iso}), this is easily seen to be the general form of a SD non-principal null 2-plane. \end{proof}

\begin{cor} The map $\psi\mapsto W_\psi$ defines an $\SLth$-equivariant embedding 
$$\PC\hookrightarrow \T^+\M$$ 
whose image is the set $\T^+_*\M$ of non-principal SD 2-planes in $T\M$ (the non-integrable locus of the twistor distribution $\D^+$).  
\end{cor}

Now combining this last Corollary  with  Prop.~\ref{ident}, we obtain the 
identifications
$$\Q \simeq\T^+_*\M\simeq\PC.$$
Tracing through our definitions, we find

\begin{proposition}
There is an isomorphism of principal $\R^*$-bundles over $\M$
$$\Q\eto\PC$$
sending $(\q,\p)\in\Q$  to the projective contact element $(q,p,\psi)$, where $q=[\q], p=[\p]$ and $\psi:T_q\RPt\to T_p\RPts$ is given in homogeneous coordinates by $$\psi([\bv])=[\q\times\bv].$$ 
That is, if $v=d\pi_\q(\bv),$ then $\psi(v)=d\bar \pi_\p(\q\times\bv)$. 
\end{proposition}

\begin{definition}
A parametrized curve $(q(t), p(t), \psi(t))$ in $\PC$ satisfies the {\em no-slip condition} if  $$\psi(t)q'(t)=p'(t)$$ for all $t$.
\end{definition}

\begin{proposition}The projection $\PC\to \M$ defines a bijection between curves in $\PC$ satisfying the no-slip condition and non-degenerate null-curves in $\M$. 
\end{proposition}

\begin{proof} If $(q(t), p(t), \psi(t))$ satisfies the no-slip condition then $(q'(t),p'(t))\in W_{\psi(t)}$, which is a null plane, hence $(q'(t),p'(t))$ is a null vector. Conversely, if $(q(t), p(t))$ is  null and non-degenerate then for all $t$ there is a unique non-principal  SD null 2-plane $W_t$ containing the null vector $(q'(t), p'(t))$. By the previous proposition, there is a unique $\psi(t)$ such that $W_t=W_{\psi(t)}$, hence $\psi(t)q'(t)=p'(t)$ and so 
$(q(t), p(t), \psi(t))$ satisfies the no-slip condition. \end{proof}


\subsubsection{The normal acceleration}

\begin{definition}Given a parametrized regular curve $q(t)$ in $\RP^2$, i.e. $q'(t)\neq 0$, 
its {\em normal acceleration}, denoted by  $q''$, is  a section  of the normal line bundle of the curve, defined as follows:  
lift the curve to $\q(t)$ in $\R^3\setminus 0$, then  let 
$$q'':=  d\pi_{\q}(\q'')\;(\mod q'), $$ where $\pi:\R^3\setminus 0\to \RPt$ is the canonical projection, $\q\mapsto [\q].$ 
\end{definition}

\mn {\bf Claim:} {\em  this definition is independent  of the lift $\q(t)$ chosen.} 

\begin{proof} Note first  that $\R\q=\Ker(  d\pi_\q)$ and that $  d\pi_\q\q'=q'$. Now if  we modify the lift by $\q\mapsto \lambda\q$, where $\lambda$ is some non-vanishing real function of $t$, 
then $\q''\mapsto (\lambda\q)''\equiv\lambda\q''\;(\mod \q, \q')\ent  d\pi_{\lambda\q}(\lambda\q)''
\equiv  d\pi_{ \lambda\q}(\lambda\q'')=
  d\pi_{\lambda\q}(d\lambda_\q(\q''))=
d(\pi\circ\lambda)_{\q}(\q'')=
d\pi _\q(\q'')\;(\mod q')$. \end{proof}

\begin{rmrk}  If we write $q(t)$ in some affine coordinate chart, $q(t)=(x(t), y(t))$, then    the above definition implies that $q''=(x'', y'')\; \mod (x',y')$. The disadvantage of this simple formula is that it is not so easy to show directly that this definition is independent of the affine coordinates chosen (the reader is invited to try). 
\end{rmrk}

\begin{definition}\label{def:inflect} An inflection point of a regular curve in $\RP^2$ is a point where the normal acceleration vanishes. 
\end{definition}
\begin{figure}[H]\centering
\includegraphics[width=0.18\textwidth]{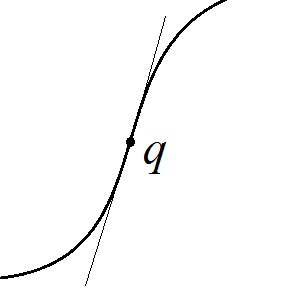}
\caption{An inflection point}
\end{figure}
\begin{rmrk} It is easy to check that the definition is parametrization independent. In fact, it is equivalent to the following, perhaps better-known, definition:  an inflection point   is a point where the tangent line has a higher order of contact  with the curve than expected  (second order or higher).  
 \end{rmrk}

 Now given a curve $(q(t), p(t), \psi(t))$ in $\PC$, 
if it  satisfies the no-slip condition,  $\psi(t)q'(t)=p'(t)$,  
then $\psi(t)$ induces a bundle map, denoted also by $\psi(t)$,  
between the normal line bundles along $q(t)$ and $p(t)$.

 \begin{proposition}For any  curve $(q(t), p(t), \psi(t))$ in $\PC$ satisfying the no-slip condition $\psi(t)q'(t)=p'(t)$, $$\psi(t) q''(t)=p''(t).$$
 \end{proposition}
 
\begin{proof} First note that for the normal accelerations $q'', p''$
 to be well-defined, both $q', p'$ must be non-vanishing, 
 i.e.  $\Gamma(t):=(q(t), p(t))$ is a non-degenerate null curve 
 in $\M$ (see Def.~\ref{ng}). It follows (see Lemma~\ref{adapted}), 
 that we can choose  an {\em adapted} lift $\sigma$ of $\Gamma$ to $\SLth$ with associated coframing  
 $\weta=\sigma^*\eta$ and dual framing $\{f_1, f_2, f^1, f^2\}$ such that  
 $\Gamma'=f_2+f^1$. Let $\tG=\tj\circ\sigma$, with $\tG(t)=(\q(t), \p(t))$, $s^i_{\,j}=\om{i}{j}(\tG'),$ $E_i(t)=E_i(\sigma(t)).$ 
Then  $\q'=E_3'=E_2+s^3_{\,3}E_3\ent \q''\equiv E_2'\equiv s^1_{\,2} E_1\;(\mod \q, \q')\ent
q''\equiv s^1_{\,2} f_1\;(\mod q'),$ and similarly $p''\equiv s^1_{\,2} f^2\;(\mod p').$ Now $W^+=\Span\{f_2+f^1, f_1-f^2\}\ent \psi f_2=f^1\ent  \psi q''=\psi(s^2_{\,1}f_2)=s^2_{\,1}f^1=p'' \;(\mod p').$
\end{proof}

\begin{cor}For a pair of regular curves $(q(t), p(t))$ satisfying the dancing condition (equivalently, 
$\Gamma(t)=(q(t), p(t))$ is a non-degenerate null-curve in $\M$), $q(t)$ is an inflection point if and only if $p(t)$ is an inflection point. 
\end{cor}

\subsubsection{Osculating conics and Cartan's developments}\label{osc}
 To complete the ``projective rolling" interpretation of   $(\Q,\D)$ we introduce a  projective connection associated with a plane curve $\gamma\subset\RPt$, defined on  its fibration of osculating conics $\CC_\gamma$; the  associated horizontal curves of this connection project to plane curves  which are the ``Cartan's developments" of $\gamma$. The ``no twist" condition for projective rolling is then   expressed in terms of this   connection, in analogy with the rolling of riemannian surfaces.

\newcommand{\LL}{\mathcal L}

 Let    $\gamma\subset\RPt$ be a smooth \lc\ curve (i.e. without   inflection points).  
For each $q\in\gamma$ there is a unique conic $\CC_q\subset \RPt$ which is tangent to 
$\gamma$ to order 4 at $q$ (this is the projective analog of the osculating circle to a curve 
in euclidean differential geometry). Define $$\CC_\gamma=\{(q,x)\st q\in \gamma, x\in \CC_q\}\subset \gamma\times\RPt.$$ We get a fibration 
\renewcommand{\II}{\mathbb I}
$$\CC_\gamma\to \gamma, \quad (q,x)\to q.$$

\begin{rmrk}  The fibration $\CC_\gamma\to \gamma$ has some remarkable properties. We refer the reader to the beautiful article \cite{GTT}, from which the following figure is taken. 
 \end{rmrk}
\begin{figure}[h]\centering
\includegraphics[width=0.6\textwidth]{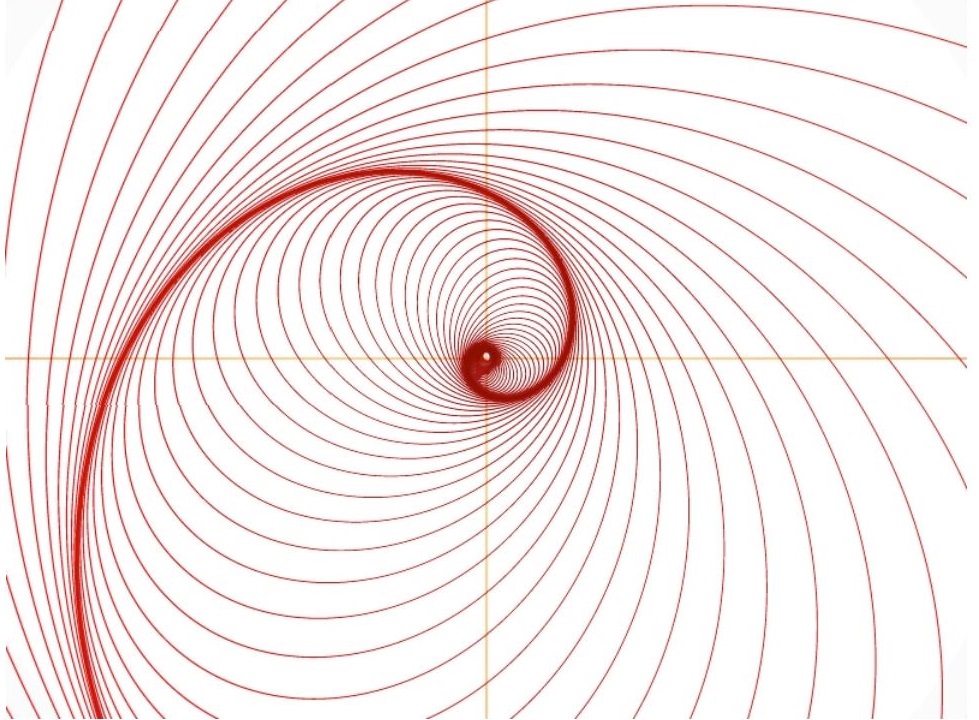}
\caption{\small The osculating conics of a spiral }
\end{figure}
There is a {\em projective connection}  defined on $\CC_\gamma\to \gamma$, i.e. a line field on $\CC_\gamma$, transverse to the fibers, whose associated parallel transport identifies the fibers of $\CC_\gamma$ projectively; its integral curves (the horizontal lifts of $\gamma$ to $\CC_\gamma$)   are defined   as follows: if we parametrize $\gamma$ by $q(t)$, then its horizontal lifts are parametrized curves $(q(t),x(t))\in\CC_\gamma$ such that  $x(t)\in\CC_{q(t)}$ is tangent to the line $\ell(t)$ passing through  $x(t)$ and  $q(t)$.  The projections $x(t)$ of such  horizontal curves on $\RPt$  are   {\em Cartan's developments} of 
$\gamma$ (see \cite[p.~58]{Cbook}). 

\begin{figure}[H]\centering
\includegraphics[width=0.4\textwidth]{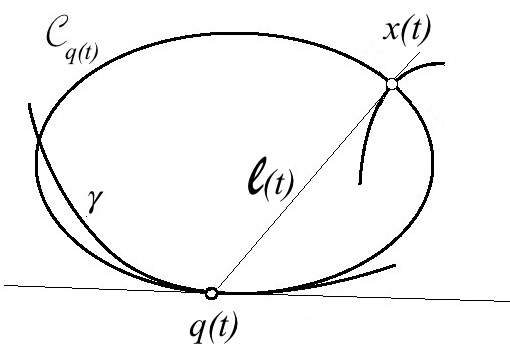}
\caption{\small Cartan's development  $x(t)$ of $\gamma$}
\end{figure}

Next we consider another  fibration of projective lines along $\gamma$
$$\LL_\gamma:= \P(T\RPt)|_\gamma\to \gamma.$$
The fiber over $q\in\gamma$ is the projectivized tangent space $\P(T_q\RPt)$, which we can also identify with $\hat q\subset\RPts$ (the pencil of lines through $q$).

 We    identify the  fiber bundles $\CC_\gamma\simeq \LL_\gamma$ using the usual ``stereographic projection": a point $x\in \CC_{q}$, $x\neq q$, is mapped to  the line $\ell$ joining $x$ with $q$, while $q$ itself is mapped to the tangent line to $\gamma$ at $q$. 
 Thus the projective connection on $\CC_\gamma$ defines, via the identification $\CC_\gamma\simeq \LL_\gamma,$
a projective connection on $\LL_\gamma.$

\subsubsection{Examples of developments.} \label{examples}
 
These are important examples and will be used later.

\begin{enumerate}[leftmargin=18pt,label=(\arabic*)]\setlength\itemsep{5pt}
\item    Parametrize  a \lc\ curve $\gamma\subset\RPt$ by $A(t)$   in LF form, i.e. $A'''+rA=0$ (see Prop.~\ref{TwentyOne}). Using  homogeneous coordinate $(x,y,z)$ on $\RPt$  with respect to the frame $A(t),A'(t),A''(t)$, the osculating conic  $\CC_t$ at  $[A(t)]$    is given by 
$y^2=2xz$  (see \cite{Cbook}, p.~55). 

In particular, taking  $x=y=0$, we get that  
$[A''(t)]$  is on the osculating conic at $[A(t)]$. In   fact:  {\em $ x(t):=[A''(t)]$
is  a  development of $\gamma$}. 

\begin{proof}  $(A'')' =-r A,$ so the tangent line to $A''(t)$ passes through $A(t)$. \end{proof}

The associated parallel line $\ell(t)$ along $\gamma$ is given by $a'=A\times A''$ where  $a=A\times A'$ is the dual curve.

\item  In fact, the development $[A''(t)]$ of the previous item is not so special. It is easy to see that for any point $x\in\CC_q$ (other then $x=q$), one can   pick a parametrization   $A(t)$ of $\gamma$ in LF form such that $x=[A''(0)]$. 

\begin{proof} (Sketch).  Start with any $A(t)$ such that $q=[A(0)]$, than find a M\"obius transformation $\bar t=f(t)$ such that $f(0)=0$ and  $\bar A(\bar t)=f'(t)A(t)$ satisfies $x=[\bar A''(0)].$\end{proof}

\item    Another way to get all  developments of $\gamma$, using the notation of the 1st example,  is to parametrize $\CC_t$ by $P(u)=A(t)+uA'(t)+(u^2/2)A''(t)$, than the developments are given by $x(t)=[P(u(t))], $ where  $u(t)$ satisfies $u'+1=0,$ i.e. $P_c(t)=A(t)+(c-t)A'(t)+[(c-t)^2/2]A''(t) $ is a  development of $\gamma$ for every  constant $c$. (Note that these developments miss exactly the first example $x(t)=[A''(t)]$ above). 

Using this formula, Cartan shows that every development curve $P_c(t)$ is tangent to  $\gamma$ as $t\to c$, with a cusp at $t=c$.

\item    Consider the  curve $\gamma^*\subset \RPts$ dual to a curve $\gamma\subset\RPt$ with a parametrization $A(t)$ in LF form. Parametrize $\gamma^*$  by $a =A\times A'$. One can check easily that $a(t)$ satisfies $a'''-ra=0$, so is also in LF form. It follows, as  in the last example,  that $a''(t)$ is a  development  of $\gamma^*$. The associated   ``parallel line"  along $a(t)$ (a point on $a(t)$) is $A'=a\times a''$. 

\begin{figure}[h]\centering
\includegraphics[width=0.6\textwidth]{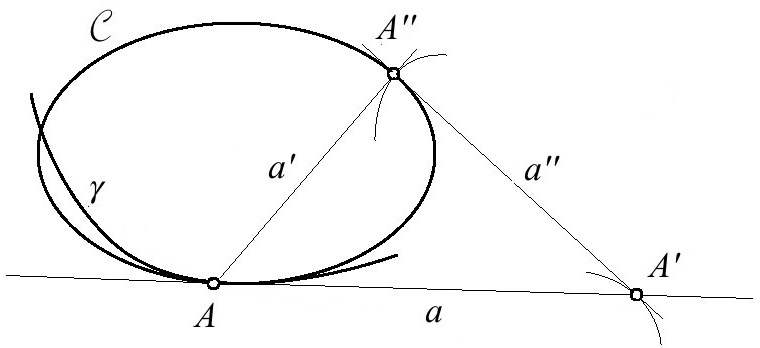}
\caption{Development of the 2nd kind (development of the dual curve)}
\end{figure}

\begin{rmrk} Cartan calls the curve $A'(t)$ a {\em development of the 2nd kind} of $\gamma$.  It can be also characterized as the envelope (or dual) of the family of tangents to osculating conics along the development  $A''(t)$ (of the 1st kind). 
\end{rmrk}

\item    When $\gamma$ is itself a conic $\CC$, then the osculating conic is obviously $\CC$ itself for all $q\in \CC$, hence  the  development curves $(q(t),x(t))$ satisfy  $x(t)=const$. It follows that  if we parallel transport a line along a conic, we get a family of {\em concurrent} lines $\ell(t). $
\begin{figure}[h]\centering
\includegraphics[width=0.35\textwidth]{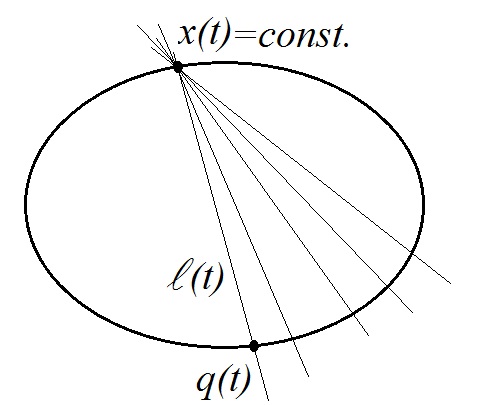}
\caption{Parallel transport of a line along a conic}
\end{figure}

\end{enumerate}

\subsubsection{The no-twist condition} 

\begin{definition}
A {\em projective rolling without slipping or twisting} of $\RPt$ along $\RPts$ is a parametrized curve $(q(t), p(t), \psi(t))$ in $\PC$, satisfying for all $t$

\begin{itemize}
\item the no-slip condition: $\psi(t)q'(t)=p'(t)$;

\item the no-twist condition: if $u(t)$ is a parallel section of $\P(T\RPt)$ along $q(t)$,  then $\psi(t)u(t)$ is a parallel section of 
 $\P(T\RPts)$ along  $p(t)$. 
\end{itemize}
\end{definition}

\begin{proposition}\label{TwentyNine}
Under the identification $ \Q\simeq\PC$, integral curves of the Cartan-Engel distribution $(\Q,\D)$ correspond to projective rolling  curves in $\PC$ satisfying the no-slip and no-twist condition.

\end{proposition}

\begin{proof}  Let $(\q(t), \p(t))$ be an integral curve of $(\Q, \D)$ and $(q(t), p(t), \psi(t))$  the corresponding projective rolling curve in $\PC$. Then $(q(t), p(t))$ is a null curve in $\M$ hence $(q(t), p(t), \psi(t))$ satisfies the no-slip condition.  Let $\ell(t)$ be a parallel line along $q(t)$. We need to show that $\ell^*(t):=\psi(t)\ell(t)=\ell(t)\cap p(t)$ is parallel along $p(t)$. Pick a projective parameter $t$   for  $q(t)$ and    a lift $A(t)$ of $q(t)$ to $\R^3\setminus 0$ such that 
  (1) $A'''+rA=0$ (the LF form) and (2) $\ell(t)$ is the line $a'=A\times A''$ connecting $A(t)$ and $A''(t)$ (see Example (2)  in Sect. \ref{examples}). Now $p(t)$ is a dancing mate of $q(t)$ hence  its dual $p^*(t)$  is  given   by $B =xA+yA'$, where $B'''\times B=0, x+y'=0$. It follows that  $\ell^*(t)=[B'(t)]$ (see the remark following Proposition \ref{TwentyThree}),  which is parallel along $p(t)=[b(t)],$ by  Example (4)  in Sect. \ref{examples}. \end{proof} 
  
\begin{figure}[H]\centering
\includegraphics[width=.9\textwidth]{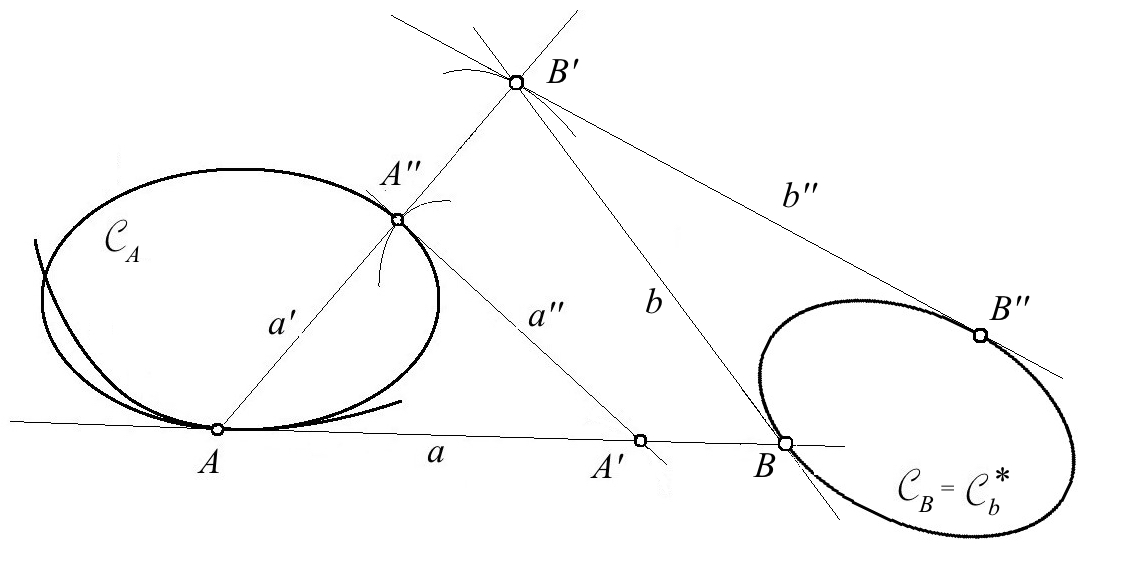}
\caption{\small The proof of Proposition \ref{TwentyNine}}
\end{figure}

\end{document}